\documentclass[review,3p]{elsarticle}

\usepackage{framed,multirow}

\usepackage{amssymb}
\usepackage{latexsym}

\usepackage{url}
\usepackage[table,xcdraw]{xcolor}
\usepackage{bm} 
\usepackage{physics}
\usepackage{amsmath}
\usepackage{mathtools}
\usepackage{booktabs}
\usepackage{yhmath}
\newcommand{\volume}{{\ooalign{\hfil$V$\hfil\cr\kern0.08em--\hfil\cr}}}
\usepackage{graphicx}
\usepackage[labelfont=it]{caption}
\usepackage{subcaption}
\usepackage{float}

\usepackage{cancel}
\usepackage[normalem]{ulem}

\usepackage{algorithm}
\usepackage{algpseudocode}

\usepackage{tikz}

\usepackage{pifont}
\usepackage{hyperref}

\journal{Computer \& Fluids}

\begin{document}

\begin{frontmatter}

\title{Investigation of Shock-Capturing with Bound-Preserving Limiters for the Nonlinearly Stable Flux Reconstruction Method}%
\author[1]{Sai Shruthi Srinivasan\corref{cor1}}
\cortext[cor1]{Corresponding author.}

\ead{sai.srinivasan@mail.mcgill.ca}
\author[1]{Siva Nadarajah}
\address[1]{Department of Mechanical Engineering, McGill University, Montreal, Quebec H3A OC3, Canada}

\begin{abstract}
Nonlinearly stable flux reconstruction (NSFR) combines the key properties of provable nonlinear stability with the increased time step from energy-stable flux reconstruction. The NSFR scheme has been successfully applied to unsteady compressible flows in arbitrary curvilinear coordinates, utilizing low-storage, weight-adjusted approaches to scale efficiently with low memory consumption. Through the use of a bound-preserving limiter, positivity of thermodynamic quantities is preserved, and this enables the extension of this scheme to hyperbolic conservation laws. We extend the bound-preserving limiter of~\citet{ZHANG20103091} to ensure robustness of the limiter for the proposed scheme. The limiter is modified to consider the minimum density and pressure values at the solution nodes when determining the value to scale the solution. The modifications are thoroughly tested with a suite of test cases. In addition to the limiter modifications, this paper conducts a thorough investigation into the shock-capturing capabilities of the NSFR scheme and the advantages it presents over standard discontinuous Galerkin (DG) methods, where, on select variants of the flux reconstruction scheme, essentially oscillation-free solutions are demonstrated. Various parameters of the scheme are extensively tested and analyzed through several 1D and 2D compressible Euler tests that verify the high-order accuracy, entropy stability, time step advantage and shock-capturing capabilities of the NSFR scheme. These parameters of interest include the choice of two-point flux, the choice of quadrature nodes and the strength of the flux reconstruction parameter. In addition to investigating the impact of the various two-point fluxes, this paper also presents numerical studies to determine the CFL condition required to maintain positivity for the two-point flux of choice. The investigation yields insightful results for all parameters, with the results pertaining to the type of flux reconstruction scheme being of special interest. The tests showcase increased robustness, time step advantages and oscillation/overshoot mitigation when employing a stronger flux reconstruction parameter.

\end{abstract}


\end{frontmatter}

\section{Introduction}\label{section: Introduction}

The discontinuous Galerkin (DG) method, proposed by \citet{reed1973triangular}, combines both the key properties of finite volume and finite element schemes. As explained in \citet{NDG}, the high-order scheme provides stability through a numerical flux function and utilizes high-order shape functions to represent the solution. An attractive extension of the DG approach leads to the computationally efficient and robust Summation-by-Parts (SBP) approach \cite{gassner2013skew,fernandez2014generalized,fernandez2019entropy,crean2018entropy,chan2018discretely,ranocha2016summation, fisher2013discretely,montoya2021unifying}. The SBP property allows for the implementation of split forms, initially proposed in \citet{tadmor1984skew}, and in \citet{gassner2013skew} for high-order schemes, which ensures entropy stability on extremely coarse meshes for unsteady problems.

Another variation of the DG method is Flux Reconstruction (FR), initially proposed by \citet{huynh_flux_2007} and later presented as a class of energy stable flux reconstruction (ESFR) schemes \cite{vincent_new_2011,jameson_proof_2010,wang2009unifying}. Through the introduction of correction functions, which are equivalently viewed as a filtered DG correction field \cite{zwanenburg_equivalence_2016,allaneau_connections_2011, aalund2019encapsulated}, ESFR recovers various high-order schemes such as the spectral difference from \citet{liu_spectral_2006} and Huynh's lumped-Lobatto scheme \cite{huynh_flux_2007}. \citet{ranocha2016summation} was the first paper to merge a collocated DG split form and ESFR, although the authors only proved stability for a collocated DG scheme. For the Euler Equations, \citet{abe2018stable} took a different approach and were able to demonstrate entropy stability for the specific ``{$\text{g}_2$ lumped-Lobatto}" ESFR split form, which corresponds to a collocated DG scheme on Gauss-Lobatto-Legendre (GLL) quadrature.

Recent developments by \citet{CicchinoNonlinearlyStableFluxReconstruction2021,cicchino2022provably} proved that nonlinear stability for ESFR schemes can only be satisfied if the FR correction functions are applied to the nonlinear volume terms. This led to the development of nonlinearly stable flux reconstruction (NSFR), where provable nonlinear stability was accomplished for the general FR case. In \citet{cicchino2025discretely}, the application of NSFR was extended for the Euler equations, and a low-storage, weight-adjusted formulation was applied for compressible flows in curvilinear coordinates. 

Physically relevant entropic solutions satisfy the strict maximum principle \cite{CMDAFERMOS}, the TVD property \cite{harten1997highresolution}, and positivity preservation for hyperbolic conservation laws. Unfortunately, numerically satisfying the second law of thermodynamics through entropy stability does not ensure that the numerical approximation also satisfies these additional properties. Higher-order entropy-stable methods improve the robustness of solvers but don’t address the problems of spurious oscillations or nonphysical density and pressure in the vicinity of a shock. Additionally, two-point fluxes, which are used for FR methods, don't have an established CFL condition under which they are expected to preserve positivity for density and pressure.  Positivity preservation strategies are essential to the success of entropy-stable methods and have been a key focus within the research community. There are a few common strategies for controlling oscillations and preventing nonphysical properties for entropy-stable methods, including artificial viscosity, subcell limiting, and limiter-based strategies. 

Artificial viscosity (AV) methods are a classical technique that dates as far back as the 1940s as seen in \cite{richtmyer1950method}. The strategy involves introducing an additional term in the PDE to control the amount of dissipation added based on shock detectors that distinguish between smooth regions and shock waves. AV methods have been successfully used to develop shock-capturing methods that are compatible with discretely entropy stable split forms as seen in \citet{wintermeyer2017entropy, wintermeyer2018entropy}. While this is a popular method, there are critical disadvantages to this strategy, such as potentially changing the order of the PDE or time-step limitations introduced as a result of enforcing boundary conditions for the extra terms. 

Subcell limiting, introduced by \citet{dumbser2014posteriori}, is a novel idea wherein the DG cell is broken down into subcells and recomputed within the troubled cells from the old time level using a different and more robust numerical scheme on a subgrid level. Subcell limiting has seen much development in the context of entropy-stable methods wherein the strategy has been modified from a hard switch to the seamless blending of high-order and low-order approaches to eliminate oscillations while maintaining entropy stability as seen in \citet{hennemann2021provably, rueda2021subcell, lin2024high, sonntag2014shock, vilar2019posteriori}. While this approach is effective, it has a rigorous implementation compared to other limiting strategies. Several components are involved, such as projecting from the cell to subcells, a detection procedure, recomputation of the solution using a different method for all troubled subcells, and a gathering procedure to establish a cell-centred polynomial on the main grid with the corrected subcells.

The class of bound-preserving limiters established by \citet{zhang2010maximumpreservinglimiter} are designed to preserve strict bounds for DGSEM and FV schemes while maintaining the order of accuracy of the original schemes. These limiters are easily extended to multi-dimensions and unstructured meshes. Bound-preserving limiters achieve strict bounds by scaling the values locally in cells that do not fall within the prescribed bounds. By using the bounds defined by the strict maximum principle or positivity preservation properties to determine whether a cell needs to be limited, the need for a troubled cell indicator is eliminated. Additionally, the cost of the limiter is minimal as it only involves the evaluation of the unlimited solution polynomial at a few predetermined quadrature points, and the limiter preserves the original cell average, thereby ensuring conservation. Bound-preserving limiters have been shown to work well with entropy-stable methods by \citet{chen2017entropy}. The limiters are simple, efficient and don’t rely on a shock detection method, unlike the aforementioned strategies; thus, we choose to incorporate bound-preserving limiters to ensure positivity preservation and enable the implementation of NSFR in a shock-capturing context. Despite all the aforementioned strengths of the positivity-preserving limiter, it still lacks robustness guarantees. To address this, modifications are proposed for the higher dimension ($>1$) implementations of the limiter.

The paper first introduces an advancement of the Zhang-Shu bound-preserving limiter~\cite{ZHANG20103091} for the current state-of-the-art NSFR schemes; second, it demonstrates the robustness of the approach for a large suite of cases without the need for TVD limiters; third, we numerically demonstrate that the two-point fluxes abide by a CFL condition; and lastly, we show that the scheme provides essentially oscillation-free solutions with select flux reconstruction schemes. 

In Sec.~\ref{sec: Prelim} we outline the DG scheme, the corresponding ESFR scheme, and then present the NSFR scheme for compressible flows. In Sec.~\ref{sec: methods}, we describe the implementation of the positivity-preserving limiter, including modifications made to enhance robustness. In addition to the limiter, Sec.~\ref{sec: methods} also outlines other critical implementation details including the two-point flux, temporal discretization method, and the flux reconstruction parameter. Lastly, in Sec.~\ref{sec:numericalresults}, we numerically verify that the proposed limiter modifications provide enhanced robustness and that the limiter preserves positivity and maintains a high-order of accuracy. This is achieved through a grid refinement order of accuracy study. We also assess the impact of the choice of two-point flux, quadrature nodes, and flux reconstruction parameters on problems involving strong shocks using various 1D and 2D numerical experiments. In addition to this extensive study, certain test cases are also used to numerically verify and establish a CFL condition for two-point fluxes that preserves positivity for density and pressure. The findings of this assessment are summarized in Sec.~\ref{sec: conclusion}.

\section{Preliminaries}\label{sec: Prelim}
\subsection{Nonlinearly Stable Flux Reconstruction Scheme}
In this section, the nonlinearly stable FR (NSFR) scheme will be presented. The preliminaries, including the discontinuous Galerkin (DG) formulation and energy stable flux reconstruction (ESFR) scheme, are first discussed, followed by the NSFR discretization scheme used in this investigation. In this paper, row vector notation will be used. We will also use boldface for vectors and italic for scalar values, and the $c$-superscript refers to the use of Cartesian coordinates.

\subsubsection{DG Formulation}\label{sec: DG}

Consider the scalar 3D conservation law,

\begin{equation}
\begin{aligned}
    \frac{\partial}{\partial t}u(\bm{x}^c,t) +\nabla\cdot\bm{f}(u(\bm{x}^c,t) )&=0,\:t\geq 0,\:\bm{x}^c\coloneqq[x,y,z]\in\Omega,\\
    u(\bm{x}^c,0)&=u_0(\bm{x}^c),
\end{aligned}\label{eq: consLaw}
\end{equation}

\noindent where $\bm{f}(u(\bm{x}^c,t) )$ represents the flux. The computational domain, $\Omega^h$ is partitioned into $M$ non-overlapping elements, denoted by $\Omega_m$,

\begin{equation}
    \Omega \simeq \Omega^h \coloneqq \bigcup_{m=1}^M \Omega_m.
\end{equation}

The global solution can then be taken as the direct sum of each approximation within each element,
\begin{equation}
    u(\bm{x}^c,t) \simeq u^h(\bm{x}^c,t)=\bigoplus_{m=1}^M u_m^h(\bm{x}^c,t).
\end{equation}

On each element, we represent the solution with $N_p$ linearly independent modal or nodal polynomial basis of a maximum order of $p$; where, $N_p=(p+1)^d$, $d$ is the dimension of the problem. The solution can be represented as

\begin{equation}
     u_m^h(\bm{x}^c,t)  = \sum_{i=1}^{N_p}{{\chi}_{m,i}(\bm{x}^c)\hat{u}_{m,i}(t)},
\end{equation}

\noindent where 

\begin{equation}
     \bm{\chi}_m(\bm{x}^c) \coloneqq [\chi_{m,1}(\bm{x}^c)\text{ }\chi_{m,2}(\bm{x}^c) \dots \chi_{m,N_p}(\bm{x}^c)]=\bm{\chi}(x)\otimes \bm{\chi}(y)\otimes \bm{\chi}(z)
\end{equation}

\noindent holds the polynomial basis functions for the element.

The physical coordinates are mapped to the standard reference space, $    \bm{\xi}^r \coloneqq \{ [\xi \text{, } \eta \text{, }\zeta]:-1\leq \xi,\eta,\zeta\leq1 \}
$, on each element with the transformation mapping 
$\bm{x}_m^c(\bm{\xi}^r)\coloneqq \bm{\Theta}_m(\bm{\xi}^r)=\Theta_{m,i}(\bm{\xi}^r)\hat{\bm{x}}_{m,i}^c,
    \nonumber
$ where $\Theta_{m,i}(\bm{\xi}^r)$ are the mapping shape functions, and $\hat{\bm{x}}_{m,i}^c$ are the physical mapping control points. In this work, the reference coordinates, $\bm{\xi_r}$, are chosen as Gauss-Lobatto-Legendre (GLL) nodes for the solution bases and either GLL or Gauss-Legendre (GL) nodes for the flux bases. The impact of the choice of flux nodes for shock problems is discussed in Section~\ref{sec: 1DShuOsherProblem}.

Continuing with the DG formulation, the elementwise reference residual for Eq.~\ref{eq: consLaw} is,
\begin{equation}
     R_m^h(\bm{x}^c,t)=\frac{\partial}{\partial t}u_m^h(\bm{x}^c,t) + \nabla \cdot \bm{f}(u_m^h(\bm{x}^c,t)).\label{eq: residual}
\end{equation}

Following a DG framework, we integrate the residual over the computational domain, utilizing orthogonality of the residual with the test function, chosen to be the basis function. After integrating by parts twice we arrive at the strong form, 

\begin{equation}
    \bm{M}_m\frac{d}{dt}\hat{\bm{u}}_m(t)^T+\bm{S}_{\xi}\hat{\bm{f}}_{1m}^r(t)^T 
    +\bm{S}_{\eta}\hat{\bm{f}}_{2m}^r(t)^T+\bm{S}_{\zeta}\hat{\bm{f}}_{3m}^r(t)^T
    +\sum_{f=1}^{N_f}\sum_{k=1}^{N_{fp}} \bm{\chi}(\bm{\xi}_{f,k}^r)^T{W}_{f,k}[\hat{\bm{n}}^r\cdot {\bm{f}}_m^{{C,r}}]=\bm{0}^T\label{eq:dg strong}
\end{equation}

\noindent where $\bm{\xi_v^r}$ represents the reference coordinate at the volume quadrature nodes and $\bm{\xi^r_{f,k}}$ represents the reference coordinate on the face $f$ at the facet cubature node $k$, where $k\in[1,N_{fp}]$ with $N_{fp}$ being the number of facet cubature nodes on the surface. The mass and stiffness matrices are defined as,

\begin{equation}
    \begin{aligned}
      (\bm{M}_m)_{ij}\approx\int_{\Omega_r}J_m^\Omega {\chi}_i(\bm{\xi}^r){\chi}_j(\bm{\xi}^r)d\Omega_r \to \bm{M}_m= \bm{\chi}(\bm{\xi}_v^r)^T\textbf{w}\bm{J}_m\bm{\chi}(\bm{\xi}_v^r),\\
      (\bm{S}_\xi)_{ij}=\int_{\Omega_r} {\chi}_i(\bm{\xi}^r)\frac{\partial}{\partial\xi}{\chi}_j(\bm{\xi}^r)d\Omega_r \rightarrow \bm{S}_\xi= \bm{\chi}(\bm{\xi}_v^r)^T\textbf{w}\bm{\chi}_\xi(\bm{\xi}_v^r),
    \end{aligned}\label{eq: mass}
\end{equation}

\noindent and similarly for the other directions. 

Here, $\textbf{w}$ is a diagonal operator storing the quadrature weights of integration at the volume cubature nodes, and $\bm{J}_m$ is a diagonal operator storing the determinant of the metric Jacobian evaluated at the volume cubature nodes. We let $\hat{\bm{n}}^r$ represent the outward pointing normal on the face at the facet cubature node in the reference element. Here, $\hat{\bm{f}}_{m}^r(t)=\bm{P}\bm{f}_{m}^r(t)$, where $\bm{P}=\bm{M}^{-1}\bm{\chi}(\bm{\xi}_{v}^r)^T\textbf{w}$ is the ${L}_2$ projection operator. That is, the modal reference flux is the modal ${L}_2$ projection of the nodal reference flux. Lastly, ${\bm{f}}_m^{{C,r}}\coloneqq({\bm{f}}_m^{*,r}) -({\bm{f}}_m^r)$ is the numerical flux minus the reference flux across the face.

\subsubsection{Corresponding ESFR Scheme}\label{sec: FR}

Following an ESFR framework~\cite{huynh_flux_2007,vincent_new_2011}, the reference flux is composed of a discontinuous and a corrected component,
\begin{equation}
    \bm{f}^r(u_m^h(\bm{\Theta}_m(\bm{\xi}^r),t))\coloneqq 
    \bm{f}^{D,r}(u_m^h(\bm{\Theta}_m(\bm{\xi}^r),t)) + \sum_{f=1}^{N_f}\sum_{k=1}^{N_{f_{p}}} \bm{g}^{f,k}(\bm{\xi}^r)[\hat{\bm{n}}^r\cdot \bm{f}_m^{C,r}].\label{eq: ESFR corr flux}
\end{equation}

As illustrated in ~\cite{allaneau_connections_2011,zwanenburg_equivalence_2016,Cicchino2020NewNorm,CicchinoNonlinearlyStableFluxReconstruction2021,cicchino2022provably} the corresponding ESFR strong form as a filtered DG scheme is given by,
\begin{equation}
    (\bm{M}_m+\bm{K}_m)\frac{d}{dt}\hat{\bm{u}}_m(t)^T
     +\bm{S}_{\xi}\hat{\bm{f}}_{1m}^r(t)^T
  +\bm{S}_{\eta}\hat{\bm{f}}_{2m}^r(t)^T+\bm{S}_{\zeta}\hat{\bm{f}}_{3m}^r(t)^T
    +\sum_{f=1}^{N_f}\sum_{k=1}^{N_{fp}} \bm{\chi}(\bm{\xi}_{f,k}^r)^T{W}_{f,k}[\hat{\bm{n}}^r\cdot {\bm{f}}_m^{{C,r}}]=\bm{0}^T.\label{eq:ESFR strong}
\end{equation}

The ESFR strong form as a filtered DG scheme is also equivalently expressed as,
\begin{equation}\label{eq: esfrStrongForm}
\begin{split}
    (\bm{M}_m+\bm{K}_m)\frac{d}{dt}\hat{\bm{u}}_m(t)^T &+ \bm{\chi}(\bm{\xi}^r_v)^T\textbf{w}\nabla^r\phi(\bm{\xi}^r_v)\cdot\hat{f}_m^r(t)^T\\
    &+\sum_{f=1}^{N_f}\sum_{k=1}^{N_{fp}} \bm{\chi}(\bm{\xi}_{f,k}^r)^T{W}_{f,k}\hat{\bm{n}}^r\cdot [\bm{f}_m^{*,r} -\phi(\bm{\xi}_{fk}^r){\bm{\hat{f}}}_m^r(t)^T]=\bm{0}^T 
\end{split}
\end{equation}

The entire influence of the ESFR correction functions are stored in $\bm{K}_m$, which is represented as,

\begin{align}
\begin{split}
    (\bm{K}_m)_{ij} &\approx \sum_{s,v,w } c_{(s,v,w)} \int_{ {\Omega}_r} J_m^\Omega \partial^{(s,v,w)} \chi_i(\bm{\xi}^r) 
    \partial^{(s,v,w)}\chi_j(\bm{\xi}^r) d {\Omega_r}\\
    \to
    \bm{K}_m &= 
     \sum_{s,v,w} c_{(s,v,w)}
    \partial^{(s,v,w)}\bm{\chi}(\bm{\xi}_v^r) ^T \textbf{w}\bm{J}_m
   \partial^{(s,v,w)}\bm{\chi}(\bm{\xi}_v^r)\\
   &=\sum_{s,v,w } c_{(s,v,w)}\Big(\bm{D}_\xi^s \bm{D}_\eta^v\bm{D}_\zeta^w \Big)^T\bm{M}_m\Big(\bm{D}_\xi^s \bm{D}_\eta^v\bm{D}_\zeta^w \Big),
\label{eq:Km}
\end{split}
\end{align}

\noindent where $\partial^{(s,v,w)}$ represents a multidimensional index spanning the $p$-th order broken Sobolev-space~\cite{sheshadri2016stability}, and its corresponding correction parameter,

\begin{equation}
\begin{split}
&\text{2D:} \indent c_{(s,v)}=c_{1D}^{(\frac{s}{p}+\frac{v}{p})},\: \text{such that } s=\{ 0,p\},\: v=\{ 0,p\},\: s+v\geq p ,\\
&\text{3D:}\indent     c_{(s,v,w)}=c_{1D}^{(\frac{s}{p}+\frac{v}{p}+\frac{w}{p})},\:\text{such that } s=\{ 0,p\},\: v=\{ 0,p\},\: w=\{ 0,p\},\: s+v+w\geq p.
    \end{split}
\end{equation}
The impact of the FR correction parameter is discussed in Section.~\ref{sec: c_param}.

\subsubsection{Nonlinearly Stable FR}
Following the derivations in Chan \cite{chan2018discretely, chan2019discretely, chan2019skew}, the general differential operator is used to recast Eqn.\ref{eq: esfrStrongForm} in a skew-symmetric two-point flux differencing form,
\begin{equation}\label{eq: nsfrStrong}
    (\bm{M}_m+\bm{K}_m)\frac{d}{dt}\hat{\bm{u}}_m(t)^T + [\bm{\chi}(\bm{\xi_v^r})^T\bm{\chi}(\bm{\xi_f^r})^T]\left[\left(\tilde{\bm{Q}} - \tilde{\bm{Q}^T}\right)\odot \bm{F_m^r}\right]\bm{1}^T+\sum_{f=1}^{N_f}\sum_{k=1}^{N_{fp}} \bm{\chi}(\bm{\xi}_{f,k}^r)^T{W}_{f,k}\hat{\bm{n}}^r\cdot \bm{f}_m^{*,r}=\bm{0}^T 
\end{equation}
where, $\tilde{\bm{Q}} - \tilde{\bm{Q}^T}$ is the general hybridized skew-symmetric stiffness operator, $\odot$ denotes a Hadamard product, and $\bm{F_m^r}$ is the matrix storing the reference two-point flux values,
\begin{equation}
    (\bm{F_m^r})_{ij} = \bm{f}_s(\tilde{u_m}(\xi_i^r),\tilde{u_m}(\xi_j^r))\left(\frac{1}{2}(\bm{C_m}(\bm{\xi_i^r}) + \bm{C_m}(\bm{\xi_j^r})\right), \forall 1 \leq i, j \leq N_v + N_{fp}
\end{equation}
where $\bm{C_m}$ is the metric Jacobian cofactor matrix formulated as described in \citet{cicchino2022provably} and $\bm{f}_s$ denotes the two-point flux which is described in Section.~\ref{sec: ec_flux}. This procedure, first proposed by Chan \cite{chan2018discretely}, is essential for entropy stability on uncollocated solution and flux nodes. For the detailed derivation and proofs for the NSFR discretization, we direct the reader to \cite{cicchino2024weight}.

\section{Methodology}\label{sec: methods}
This section presents the critical parameters in this investigation. The choice in positivity-preserving strategy, two-point flux, and FR correction parameter are key factors considered in extending the NSFR scheme to problems involving strong shocks. The impact of these factors is discussed in the following subsections. 

\subsection{A Modified Implementation of the Zhang-Shu Positivity-Preserving Limiter~\cite{ZHANG20108918}} \label{sec: limiters}
Consider the 2D Euler equations for perfect gas,
\begin{equation}
    \begin{split}\label{eq: 1DEuler}
        &\textbf{w}_t + \textbf{f(w)}_x + \textbf{g(w)}_y= 0, \: t\geq0, \: (x,y)\in\mathbb{R}^2,\\
        &\textbf{w} = \begin{pmatrix}
            \rho\\
            m\\
            n\\
            E
        \end{pmatrix}, \: \textbf{f(w)} = \begin{pmatrix}
            m\\
            \rho u^2 + p\\
            \rho uv\\
            (E+p)u
        \end{pmatrix}, \: \textbf{g(w)} = \begin{pmatrix}
            m\\
            \rho uv\\
            \rho v^2 + p\\
            (E+p)v
        \end{pmatrix},\\
        \text{with } \qquad     &m = \rho u, \quad 
        n = \rho v, \quad
        E = \rho e +\frac{1}{2}\rho \left(u^2 +v^2\right), \quad p = (\gamma-1)\rho e,
    \end{split}
\end{equation}
where $\rho$ is density, $u$ and $v$ are the velocities in the $x$-and $y$-directions, $m$ and $n$ are the momentum in their respective directions, $E$ is total energy, $p$ is pressure, $e$ is internal energy and $\gamma=1.4$ (for air). Physically, it is required that density and pressure remain positive. Failure to preserve the positivity of these two properties results in the numerical algorithm diverging. High-order numerical schemes commonly used to solve systems of hyperbolic conservation laws, such as Eq.~(\ref{eq: 1DEuler}), do not provably preserve positivity without additional tools. In ~\citet{ZHANG20108918}, techniques used to establish a maximum-principle preserving limiter \cite{ZHANG20103091} were generalized to derive a limiter designed to preserve positivity of density and pressure, provided density and pressure were positive at the previous timestep. 

The set of admissible states for Eq.~(\ref{eq: 1DEuler}) is given by,
\begin{equation}\label{eq: admissibleStates}
    G = \left\{\textbf{w}=\begin{pmatrix}
        \rho\\
        m\\
        n\\
        E
    \end{pmatrix}\: \middle| \: \rho>0 \text{ and } p=(\gamma-1)\left(E - \frac{1}{2}\rho\left(u^2 +v^2 \right)\right) > 0\right\}.
\end{equation}

Given the vector of approximation polynomials, $\bm{u}^h(\bm{x})=(\bm{\rho}(\bm{x}),\bm{m}(\bm{x}),\bm{n}(\bm{x}),\bm{E}(\bm{x}))^T$, and cell averages, $\Bar{\textbf{w}}_m=(\Bar{\rho}_m,\Bar{m}_m, \Bar{n}_m,\Bar{E}_m) \in G$, the positivity-preserving limiter is applied to transform $\bm{u}^h(\bm{x})$ to $\tilde{\bm{u}}^h(\bm{x})$ while maintaining high-order accuracy, positivity (ie. $\tilde{\bm{u}}^h(\bm{x}) \in G$) and conservativity. 

For dimensions greater than one, the limiter utilizes a mixed set of quadrature rules to determine the cell averages and minimum values for density and pressure. Using the polynomial representation of the local solution, $\bm{u}^h(\bm{x})$, we interpolate the solution to the new sets of quadrature rules:
\begin{equation}\label{eq: mixed_nodes}
    \begin{split}
        \bm{\xi}^{r,1} = \xi^\alpha \otimes \eta^\beta,\\
        \bm{\xi}^{r,2} = \xi^\beta \otimes \eta^\alpha,
    \end{split}
\end{equation}
where the sets $\{\xi^\beta,\eta^\beta:-1\leq \xi^{\beta}\text{, }\eta^{\beta}\leq1 \}$ and $\{\xi^\alpha,\eta^\alpha:-1\leq \xi^{\alpha}\text{, }\eta^{\alpha}\leq1 \}$ denote GL and GLL quadrature points, respectively.

Using these two solution spaces, the cell averaged solution vector,
\begin{equation}\label{eq: cell_avg_solution_vector_ppl}
    \Bar{\textbf{w}}_m=(\Bar{\rho}_m,(\Bar{\rho \text{v}_1})_m, (\Bar{\rho \text{v}_2})_m, \Bar{E}_m)^T
\end{equation} is computed using the quadrature rules as follows: 
\begin{align}\label{eq: cell_avg_solution_vector_ppl_2d}
    \Bar{\textbf{w}}_m &= \frac{a_1\lambda_1}{\mu}\Bar{\textbf{w}}_m + \frac{a_2\lambda_2}{\mu}\Bar{\textbf{w}}_m\\
    &= \frac{a_1\lambda_1}{\mu}\sum^L_{\beta=1}\sum^N_{\alpha=1}\omega_\alpha\omega_\beta\textbf{w}(\xi^{r,1}) + \frac{a_2\lambda_2}{\mu}\sum^L_{\beta=1}\sum^N_{\alpha=1}\omega_\alpha\omega_\beta\textbf{w}(\xi^{r,2}),
\end{align}
where $\omega$ denotes the quadrature weight, $a$ denotes maximum wavespeed, $\lambda_1 = \frac{\Delta t}{\Delta x}$, $\lambda_2 = \frac{\Delta t}{\Delta y}$ and $\mu = a_1\lambda_1 + a_2\lambda_2$.

The positivity-preserving limiter designed by Zhang and Shu, later improved by Wang and Shu \cite{WANG2012653}, enforces positivity in two steps. In the first step, the density is limited using a scaling value, $\theta_1$, that is determined with the cell average density, $\Bar{\rho}_m$, and minimum density, $\rho_{min}$. In the original implementation of the limiter, $\rho_{min}$ is determined by taking the lowest density value between the two quadrature rules. The density values in the vector of approximation polynomials are replaced with the limited density values, and the new vector of approximation polynomials is denoted as $\bm{\hat{u}}^h(\bm{x})$.

The second step is to enforce the positivity of pressure through a scaling value, $\theta_2$, that is determined by calculating the cell average pressure and pressure at the quadrature nodes. The cell average pressure is obtained using the cell averaged solution vector,r and the pressure at the quadrature nodes utilizes both mixed quadrature rules given in Eqn.~\ref{eq: mixed_nodes}. The values in the vector of approximation polynomials are scaled using $\theta_2$ and the new vector of approximation polynomials is denoted as $\bm{\Tilde{u}}^h(\bm{x})$.

From here, in the NSFR scheme with Euler forward in time, $\tilde{\bm{u}}^h(\bm{x})$ is used in the place of $\textbf{u}^h(\textbf{x})$. This limiter is applied at each stage of the RK method. It is also important to note that while the bound-preserving class of limiters provably ensures that density and pressure remain positive, the solutions are still prone to oscillations since it is not provably non-oscillatory. 

The implementation of the positivity-preserving limiter described in \citet{WANG2012653} improves the original limiter through modification to the $\theta_2$ scaling parameter but still lacks robustness guarantees as the values of density and pressure at the original solution nodes are not taken into account. This could result in negative density and pressure values at the solution nodes while these properties remain positive at the tensored quadrature points. Thus, to further ensure robustness, in this paper modifications are made to account for the values of density and pressure at the solution nodes, $\bm{\xi}^{r,3} = \xi^\alpha \otimes \eta^\alpha$, in addition to the other two sets of nodes when determining minimum density and pressure at the quadrature nodes. It should be noted that the cell averaged solution vector is still calculated as per Eq.~(\ref{eq: cell_avg_solution_vector_ppl_2d}) using the tensored quadrature points in Eqn.~\ref{eq: mixed_nodes}.

In the original implementation of the limiter, the $\theta$ values which are affected by these changes are used to scale the solution nodes directly, and as such, this modification preserves the properties of the limiter. This modification has been verified through an extensive suite of test cases and also recovers the expected orders as highlighted in Sec.~\ref{sec:numericalresults}.

\subsubsection{Implementation}
The implementation of the positivity-preserving limiter in 2D with the suggested modifications is briefly outlined in this section, and an algorithm is provided for clarity.

To achieve positivity preservation, first we define $\varepsilon$ to be a small number greater than zero that serves as a lower bound for the average density and average pressure for all cells (ie. $\Bar{\rho}_m > \varepsilon$ and $\Bar{P}_m > \varepsilon$). In this paper we use $\varepsilon = 10^{-13}$ throughout. Next, we interpolate the solution nodes to the two sets of tensored quadrature points ($\bm{u}^h(\xi^{r,1})$ and $\bm{u}^h(\xi^{r,2})$) to determine the cell averaged solution vector as per Eq.~(\ref{eq: cell_avg_solution_vector_ppl_2d}). With the cell-averaged solution vector, the average pressure in the cell can be defined as:
\begin{equation}\label{eq: cell_avg_pressure}
    \Bar{P}_m = (\gamma -1) \left(\Bar{E}_m -\frac{1}{2}\frac{\Bar{m}_m^2}{\Bar{\rho}_m}-\frac{1}{2}\frac{\Bar{n}_m^2}{\Bar{\rho}_m}\right).
\end{equation}
With the cell averages determined, the solution vector can be scaled to ensure positivity. The first scaling value is determined using the cell average density ($\Bar{\rho}_m$) and the minimum density in the cell ($\rho_{min}$); where the latter is determined using three sets of nodes, that is, the solution nodes and the two tensored quadrature sets. 
\begin{equation*}
        \rho_{\min}=\min\:\left\{\bm{\rho}(\bm{\xi}^{r,1}), \bm{\rho}(\bm{\xi}^{r,2}), \bm{\rho}(\bm{\xi}^{r,3})\right\}.
\end{equation*}
Using the cell average density and cell minimum density, the first scaling value is obtained:
\begin{equation*}
    \theta_1 = \min\left\{\frac{\Bar{\rho}_m-\varepsilon}{\Bar{\rho}_m-\rho_{\min}}, 1\right\},
\end{equation*}
which is employed to scale only density:
\begin{equation*}
     \hat{\bm{\rho}}(\bm{x}) = \theta_1(\bm{\rho}(\bm{x}) - \Bar{\rho}_m) + \Bar{\rho}_m.
\end{equation*}
The original density values of the vector of approximation polynomials, $\bm{u}^h(\bm{x})=(\bm{\rho}(\bm{x}),\bm{m}(\bm{x}),\bm{n}(\bm{x}),\bm{E}(\bm{x}))^T$ are replaced with the scaled density values and the new vector of approximation polynomials is denoted as $\bm{\hat{u}}^h(\bm{x})=(\bm{\hat{\rho}}(\bm{x}),\bm{m}(\bm{x}),\bm{n}(\bm{x}),\bm{E}(\bm{x}))^T$.

The second scaling value is then determined using the cell average pressure and the cell minimum pressure across the three sets of nodes. The pressure at each node is determined, and if the pressure violates the preservation of positivity, the scaling value at that node is set to:
\begin{equation*}
    t = \frac{P(\Bar{\textbf{w}}_m)}{P(\Bar{\textbf{w}}_m)-P(\hat{\textbf{u}}^h)}
\end{equation*}
or otherwise, it is set to 1. Once the scaling value at all the nodes for the three sets has been determined, the minimum scaling value across all sets of nodes is chosen as the second scaling value:
\begin{equation*}
    \theta_2 = \min \left\{\bm{t}_{1}, \bm{t}_{2}, \bm{t}_{3}\right\}.
\end{equation*}
The second scaling value is then used to scale the vector of approximation polynomials with limited density thereby giving us:
\begin{equation*}
    \Tilde{\textbf{u}}^h(\bm{x}) = \theta_2(\hat{\textbf{u}}^h(\bm{x}) -\Bar{\textbf{w}}_m) + \Bar{\textbf{w}}_m,
\end{equation*}
which is a vector of approximation of polynomials that preserves accuracy, positivity, and conservativity. The presented positivity-preserving limiter based on Zhang and Shu~\cite{ZHANG20108918} with the modifications presented above is detailed in Algorithm~\ref{alg:ppl}.

\begin{algorithm}
\caption{2D Positivity-Preserving Limiter Implementation with Modifications}\label{alg:ppl}
\begin{algorithmic}[1]
\Require $\bm{\rho}(\bm{\xi}^{r,3})>0 \text{ and } \bm{P}(\bm{\xi}^{r,3}) > 0$ at previous timestep, where pressure is $P=(\gamma-1)\left(E - \frac{1}{2}\frac{m^2}{\rho}- \frac{1}{2}\frac{n^2}{\rho}\right)$

\State Project solution to tensored quadrature nodes ($\bm{\xi}^{r,1}$ and $\bm{\xi}^{r,2}$)
\State Determine cell averages $\Bar{\textbf{w}}$ and $\Bar{P}$ using Eq.~\ref{eq: cell_avg_solution_vector_ppl_2d} and Eq.~\ref{eq: cell_avg_pressure}

\State $\rho_{min} \gets \min\{\bm{\rho}(\bm{\xi}^{r, }),\bm{\rho}(\bm{\xi}^{r,2}),\bm{\rho}(\bm{\xi}^{r,3})\}$

\If{$\rho_{min} < 0$}
    \State $\theta \gets \min\{\frac{\Bar{\rho}-\varepsilon}{\Bar{\rho}-\rho_{min}},1\}$
    \State $\hat{\bm{\rho}} \gets \theta_1(\bm{\rho}(\bm{\xi}^{r,3}) - \Bar{\rho})+\Bar{\rho}$
    \State Replace $\bm{\rho}$ with $\hat{\bm{\rho}}$ in $\bm{u}^h(\bm{x})$
\EndIf

\ForAll{$\bm{\xi}^r$ in $\{\bm{\xi}^{r,1},\bm{\xi}^{r,2},\bm{\xi}^{r,3}\}$}
    \If{$P(\hat{\textbf{u}}^h(\bm{\xi^{r}})) < 0$}
        \State $t(\bm{\xi}^r) \gets \frac{P(\Bar{w})}{P(\Bar{w}) - P(\hat{\textbf{u}}^h(\bm{\xi^{r}}))}$
        \State $\theta_2 \gets \min\{\theta_2,t\}$
    \EndIf
\EndFor

\State $\Tilde{\bm{u}} \gets \theta_2(\bm{u}^h(\bm{\xi^{r,3}}) - \Bar{\textbf{w}})+\Bar{\textbf{w}}$

\Ensure  $\Tilde{\bm{\rho}}(\bm{\xi}^{r,3})>0 \text{ and } \Tilde{\bm{P}}(\bm{\xi}^{r,3}) > 0$
    
\end{algorithmic}
\end{algorithm}

\subsection{Two-Point Flux}\label{sec: ec_flux}

In this paper, the split form convective two-point flux of Chandrashekar modified by Ranocha $\left(\text{CH}_{\text{RA}}\right)$ \cite{ranocha2021preventing} is used. The $\left(\text{CH}_{\text{RA}}\right)$ flux is entropy conserving (EC), kinetic energy preserving (KEP), and pressure-equilibrium-preserving (PEP). The $\text{CH}_{\text{RA}}$ numerical two-point flux in the $x$ and $y$ directions is given by:
\begin{equation}\label{eq:two-point-generic-chra}
    \bm{f}_{c} = \left[F_{\text{CH}_{\text{RA}}}^{\#},~G_{\text{CH}_{\text{RA}}}^{\#}\right],
\end{equation}
where $F_{\text{CH}_{\text{RA}}}^{\#}\equiv F_{\text{CH}_{\text{RA}}}^{\#}(\bm{U}_{L},\bm{U}_{R})$ and $G_{\text{CH}_{\text{RA}}}^{\#}\equiv G_{\text{CH}_{\text{RA}}}^{\#}(\bm{U}_{L},\bm{U}_{R})$ are computed using the left state $\bm{U}_{L}$ and right state $\bm{U}_{R}$ as:
\begin{equation}\label{eq:two-point-chra}
    F_{\text{CH}_{\text{RA}}}^{\#} = 
    \begin{bmatrix}
    \hat{\rho}\hat{u}\\
    \hat{\rho}\hat{u}^{2}+\hat{p}_{1}\\
    \hat{\rho}\hat{u}\hat{v}\\
    \hat{\rho}\hat{u}\hat{h}-\{\!\!\{ \text{v}_{1}P\}\!\!\}
    \end{bmatrix}, \quad G_{\text{CH}_{\text{RA}}}^{\#} = 
    \begin{bmatrix}
    \hat{\rho}\hat{v}\\
    \hat{\rho}\hat{v}\hat{u}\\
    \hat{\rho}\hat{v}^{2}+\hat{p}_{1}\\
    \hat{\rho}\hat{v}\hat{h}-\{\!\!\{ \text{v}_{2}P\}\!\!\}
    \end{bmatrix}
\end{equation}

\begin{equation}
    \hat{\rho}=\rho^{\ln},\quad \hat{u}=\{\!\!\{\text{v}_{1}\}\!\!\},\quad \hat{v}=\{\!\!\{\text{v}_{2}\}\!\!\},\quad \quad \hat{p}_{1}=\{\!\!\{ P\}\!\!\},\quad 
\end{equation}
\begin{equation}
    \hat{h}=\frac{1}{\hat{p}_{2}\left(\gamma-1\right)}+\frac{1}{2}\sum_{i=1}^{2}\left(2\{\!\!\{\text{v}_{i}\}\!\!\}^{2}-\{\!\!\{\text{v}_{i}^{2}\}\!\!\}\right)+\frac{2\hat{p}_{1}}{\rho^{\ln}},\quad \hat{p}_{2}=\left(\frac{\rho}{P}\right)^{\ln},
\end{equation}
where $\{\!\!\{\cdot \}\!\!\}$ denotes the mean:
\begin{equation}\label{eq:mean_operator}
    \left\{\!\!\{\cdot\right\}\!\!\} = \frac{\left(\cdot\right)_{L}+\left(\cdot\right)_{R}}{2},
\end{equation}
and $\left(\cdot\right)^{\ln}$ denotes the logarithmic mean:
\begin{equation}\label{eq:log_mean}
    \left(\cdot\right)^{\ln} = \frac{\left(\cdot\right)_{L}-\left(\cdot\right)_{R}}{\ln\left(\left(\cdot\right)_{L}\right)-\ln\left(\left(\cdot\right)_{R}\right)}.
\end{equation}

 Due to the nature of the logarithmic mean, determining a CFL condition for the two-point flux is a challenge. This paper aims to numerically establish a CFL condition instead. The two-point flux is expected to preserve positivity for pressure if it follows the condition provided in \citet{ZHANG20108918}, which states that for the Lax-Friedrichs flux, the $CFL$ condition can be relaxed to,
\begin{equation}\label{eq: cfl-condition-1}
   \frac{\Delta t}{\Delta x}\: ||(|u|+c)||_\infty \leq 1.
\end{equation}
This condition can be used for all polynomial degrees as the $\Delta x$ is determined using the minimum distance between quadrature nodes within a given cell.

In addition to the condition given above, the two-point flux is expected to preserve positivity for density if it follows the condition in \citet{ranocha2018comparison},
\begin{equation}\label{eq: cfl-condition-2}
    \frac{\Delta t}{\Delta x}\: ||(|u|+c)||_\infty \leq 0.5.
\end{equation} 
It should be noted however that this condition uses cell width for $\Delta x$ which wouldn't reflect the change in the minimum $\Delta x$ when the polynomial degree is changed. Thus, this condition is only expected to hold for $p = 2$. 

A numerical study of the CFL conditions is given in Sections \ref{sec:tpf_sod}, \ref{sec:tpf_shuosher}, and \ref{sec:tpf_svsw}.

\subsection{FR Correction Parameter}\label{sec: c_param}
The modified mass matrix seen in Eq.\ref{eq: nsfrStrong} has been shown to effectively act as a linear filter operator on the DG residual \cite{zwanenburg_equivalence_2016,allaneau_connections_2011}. The ESFR correction functions stored in $\bm{K}_m$ are scaled by the correction parameter, $c$, which can be used to tune the linear filtering operation. The correction parameter is generally chosen within the range of $c\in[c_-,c_+]$, where $c_-$ and $c_+$ denote the minimum and maximum value, respectively, for which the scheme is stable \cite{allaneau_connections_2011,castonguay2012high}. For values of $c\in[c_-,0)$, the filter will amplify the highest mode of the residual and for values of $c\in(0,c_+]$, the filter dampens the highest mode. When the correction parameter is chosen to be $c=0$, it recovers the unfiltered DG method and we denote this value as $c_{DG}$. Other notable values of $c$ include $c_{SD}$ which recovers a stable spectral difference scheme \cite{allaneau_connections_2011} and $c_{HU}$ which recovers Hyunh's $g_2$ scheme \cite{huynh_flux_2007}.

\subsection{Temporal Discretization}
In this work, the solution is advanced in time using the strong stability preserving third-order accurate Runge-Kutta (SSPRK3) explicit time advancement method \cite{shu1988efficient}. The method is presented below,
\begin{align*}
    w^{(1)} &= u^n + \Delta t \mathcal{L}(t^n, u^n)\:,\\
    w^{(2)} &= \frac{1}{4}(3u^n + w^{(1)} + \Delta t \mathcal{L}(w^{(1)},t^n + \Delta t)\:,\\
    u^{n+1} &= \frac{1}{3}(u^n + 2w^{(2)} + 2\Delta t\mathcal{L}\left(w^{(2)}, t^n + \frac{1}{2}\Delta t\right)
\end{align*}

The convective-based CFL condition used in this work is determined as follows:
\begin{align}
    \Delta t = CFL \: \frac{\Tilde{\Delta x}}{\lambda_{max}}, \quad \Tilde{\Delta x} = \frac{x_{max}-x_{min}}{(DOF)^{1/dim}}, \quad \lambda_{max} = max(|v|+c),
\end{align}
where $\Tilde{\Delta x}$ is the approximate grid spacing, $dim$ is the dimension of the problem, and $\lambda_{max}$ is the maximum wavespeed throughout the entire domain with speed of sound $c = \sqrt{\gamma R T}$.

\section{Numerical Results}
\label{sec:numericalresults}
In this section, the accuracy and robustness of the NSFR method in the context of shock-capturing are demonstrated through a variety of 1D and 2D numerical tests. The positivity-preserving limiter (PPL) enables the extension of the NSFR method to problems involving strong shocks and complex features through the guarantee of positivity preservation, but as stated in Section \ref{sec: limiters}, while positivity is preserved, the limiter does not eliminate oscillations in the solution. We investigate the stability of the NSFR scheme in the absence of either a TVD limiter or shock-capturing mechanisms such as slope limiters or artificial viscosity. 

All tests use GLL quadrature points, the Chandrashekar flux with the Ranocha pressure fix $\left(\text{CH}_{\text{RA}}\right)$, and Roe dissipation unless otherwise specified. The polynomial degree is varied for some tests to demonstrate the effect of increasing polynomial order on problems involving shocks. The polynomial degree is denoted with a $p$ followed by the degree, such as $p3$ to denote an approximation polynomial degree of 3. We also observe numerically that the two-point flux abides by the $CFL$ conditions, Eq.\ref{eq: cfl-condition-1} and \ref{eq: cfl-condition-2}, in both one and two dimensions.

In addition to accuracy and robustness, the tests aim to provide a thorough assessment of the NSFR method by observing the effects of variants of flux reconstruction schemes for shock-dominated problems. The tests assessing the correction parameter are performed with a polynomial degree of 3 ($p_3$). Table \ref{tab:FRschemes} list the variants of the flux reconstruction scheme investigated in this work. 
\begin{table}[ht]
\centering
\caption{Values of the correction parameter for polynomial degree 3}
\label{tab:FRschemes}
\begin{tabular}{@{}c||c@{}}
\toprule
\multicolumn{1}{l}{FR Parameter} & \multicolumn{1}{c}{c} \\ \midrule
Discontinuous Galerkin, $c_{DG}$                                    & $0$\\
Spectral-Difference, $c_{SD}$                                    & $7.44E-6$\\
Huynh's~\cite{huynh_flux_2007} $g_2$ FR, $c_{HU}$                                    & $1.32E-5$\\
 $c_{+}$                                     & $2.87E-5$\\
$c_{+\times 10}$                            & $2.87E-4$\\
\bottomrule
\end{tabular}
\end{table}

\newpage
\subsection{Low Density Convergence Test}
The first case used to verify the implementation of the PPL is the 2D Low Density case \cite{ZHANG20108918}. The problem is used to demonstrate that the scheme still achieves the expected order of accuracy. The problem is solved on a domain of size $[0,2\pi]\times[0,2\pi]$ with periodic boundary conditions and initial condition,
\begin{equation*}
    \rho_0(x,y) = 1 + 0.995\sin(x+y),\ u_0(x,y) = 1,\ v_0(x,y)=1,\ p_0(x,y)=1.
\end{equation*}
The exact solution
\begin{equation*}
    \rho(x,y,t) = 1 + 0.995\sin(x+y-2t),\ u(x,y,t) = 1,\ v(x,y,t)=1,\ p(x,y,t)=1
\end{equation*}
is used to determine the order of convergence. The final time is $t=0.1s$ and $\Delta x = \Delta y = \frac{2\pi}{N}$. The minimum expected density for this case is $\rho_{min} = 0.005$, thereby serving as a good benchmark case for positivity preservation of density. The results are shown in Table.~\ref{tab:convergence2D} for polynomial degrees $2$ and $3$ with grid size starting at $8\times8$ and successively refined to $512\times512$. The scheme achieves the expected order of accuracy in the finer grids.

\begin{table}[ht]
\centering
\caption{Low Density convergence test with initial data $1 + 0.999sin(x+y)$ and exact solution $1 + 0.999sin(x+y-2t)$, $\Delta x = \Delta y = \frac{2\pi}{N}$ and final time $t=0.1s$}
\label{tab:convergence2D}
\resizebox{0.7\textwidth}{!}{%
\begin{tabular}{@{}llllllll@{}}
\toprule
\multicolumn{8}{c}{\textbf{NSFR - $CH_{RA}$ - Roe Dissipation - $\mathbf{c_{DG}}$ - GLL - CFL0.5}}    \\ \midrule
k & $N\times N$    & dx         & $L_1$ error & Order & $L_2$ error & Order & Limiter    \\ \midrule
2 & $8\times8$     & 2.0833E-02 & 1.1166E-01  &       & 2.4880E-02  &       & \checkmark \\
  & $16\times16$   & 1.0417E-02 & 1.8730E-02  & 2.58  & 5.4392E-03  & 2.19  &            \\
  & $32\times32$   & 5.2083E-03 & 3.1915E-03  & 2.55  & 1.2218E-03  & 2.15  &            \\
  & $64\times64$   & 2.6042E-03 & 3.6526E-04  & 3.13  & 1.4299E-04  & 3.10  &            \\
  & $128\times128$ & 1.3021E-03 & 4.5685E-05  & 3.00  & 1.9041E-05  & 2.91  &            \\
  & $256\times256$ & 6.5104E-04 & 5.6583E-06  & 3.01  & 2.3462E-06  & 3.02  &            \\
  & $512\times512$ & 3.2552E-04 & 7.0088E-07  & 3.01  & 2.8993E-07  & 3.02  &            \\ \midrule
3 & $8\times8$     & 1.5625E-02 & 2.4676E-02  &       & 9.7284E-03  &       & \checkmark \\
  & $16\times16$   & 7.8125E-03 & 3.4557E-03  & 2.84  & 1.8439E-03  & 2.40  & \checkmark \\
  & $32\times32$   & 3.9062E-03 & 2.8400E-04  & 3.61  & 1.7284E-04  & 3.42  &            \\
  & $64\times64$   & 1.9531E-03 & 1.7344E-05  & 4.03  & 1.2619E-05  & 3.78  &            \\
  & $128\times128$ & 9.7656E-04 & 1.0449E-06  & 4.05  & 8.2709E-07  & 3.93  &            \\
  & $256\times256$ & 4.8828E-04 & 6.4863E-08  & 4.01  & 5.3026E-08  & 3.96  &            \\
  & $512\times512$ & 2.4414E-04 & 4.0222E-09  & 4.01  & 3.3402E-09  & 3.99  &            \\ \bottomrule
\end{tabular}%
}
\end{table}

In this investigation, the behaviour of the scheme when using GL flux nodes is also examined. To use the limiter with GLL solution nodes and GL flux nodes, it is essential to verify that the scheme still reaches the expected order of accuracy. The same initial condition and exact solution are used to determine the order of convergence. Similar to the GLL test, the results are shown in Table.~\ref{tab:convergence2D_GL} for polynomial degrees $2,3$ with grid size starting at $8\times8$ and successively refined to $512\times512$. The CFL for this test is modified to account for the sensitivities of GL nodes. The test was monitored to ensure that the limiter is turned on for coarser meshes. The scheme achieves the expected order of accuracy in the finer grids.

\begin{table}[ht]
\centering
\caption{Low Density convergence test with initial data $1 + 0.995sin(x+y)$ and exact solution $1 + 0.995sin(x+y-2t)$, $\Delta x = \Delta y = \frac{2\pi}{N}$ and final time $t=0.1s$}
\label{tab:convergence2D_GL}
\resizebox{0.7\textwidth}{!}{%
\begin{tabular}{@{}llllllll@{}}
\toprule
\multicolumn{8}{c}{\textbf{NSFR - $CH_{RA}$ - Roe Dissipation - $\mathbf{c_{DG}}$ - GL - CFL0.01}}    \\ \midrule
k & $N\times N$    & dx         & $L_1$ error & Order & $L_2$ error & Order & Limiter    \\ \midrule
2 & $8\times8$     & 2.0833E-02 & 2.0711E+00  &       & 5.3860E-01  &       & \checkmark \\
  & $16\times16$   & 1.0417E-02 & 5.5323E-02  & 5.23  & 1.5116E-02  & 5.16  & \checkmark \\
  & $32\times32$   & 5.2083E-03 & 8.2801E-03  & 2.74  & 2.7284E-03  & 2.47  &\\
  & $64\times64$   & 2.6042E-03 & 8.8661E-04  & 3.22  & 3.8200E-04  & 2.84  &\\
  & $128\times128$ & 1.3021E-03 & 9.3052E-05  & 3.25  & 4.3321E-05  & 3.14  &\\
  & $256\times256$ & 6.5104E-04 & 8.8377E-06  & 3.40  & 3.4600E-06  & 3.65  &\\
  & $512\times512$ & 3.2552E-04 & 8.8582E-07  & 3.32  & 2.7869E-07  & 3.63  &            \\ \midrule
3 & $8\times8$     & 1.5625E-02 & 2.0763E+00  &       & 6.0201E-01  &       & \checkmark \\
  & $16\times16$   & 7.8125E-03 & 2.0856E-02  & 6.64  & 7.8696E-03  & 6.26  & \checkmark \\
  & $32\times32$   & 3.9062E-03 & 1.8480E-03  & 3.50  & 9.0592E-04  & 3.12  &\\
  & $64\times64$   & 1.9531E-03 & 8.2785E-05  & 4.48  & 6.0332E-05  & 3.91  &\\
  & $128\times128$ & 9.7656E-04 & 3.7260E-06  & 4.47  & 2.9894E-06  & 4.33  &\\
  & $256\times256$ & 4.8828E-04 & 1.3273E-07  & 4.81  & 1.0178E-07  & 4.88  \\
  & $512\times512$ & 2.4414E-04 & 5.4963E-09  & 4.59  & 3.9881E-09  & 4.67  \\ \bottomrule
\end{tabular}%
}
\end{table}
\subsection{Sod Shock Tube}
\label{sec: 1DSodShockTest}
The first 1D test case used to investigate the shock capturing capabilities of the NSFR scheme is the Sod Shock Tube case. This case is a Riemann problem for the 1D Euler equations. The computational domain is $[-0.5, 0.5]$. The problem is initialized as,
\begin{equation}
    (\rho, w, p) = \begin{cases}
        (1,0,1) \qquad &\text{if} \: x < 0,\\
        (0.125, 0,0.1) \qquad &\text{if} \: x \geq 0.
    \end{cases}
\end{equation}
The final time is $t = 0.2s$. The number of cells is $N = 512$ for the tests unless stated otherwise. The exact solution of density contains a left rarefaction wave and a right shock wave with a contact discontinuity in the middle.

\subsubsection{Comparison with Strong DG}\label{sec:sod_strongDG}
\begin{figure}[h!]
\centering
\includegraphics[width=0.5\textwidth]{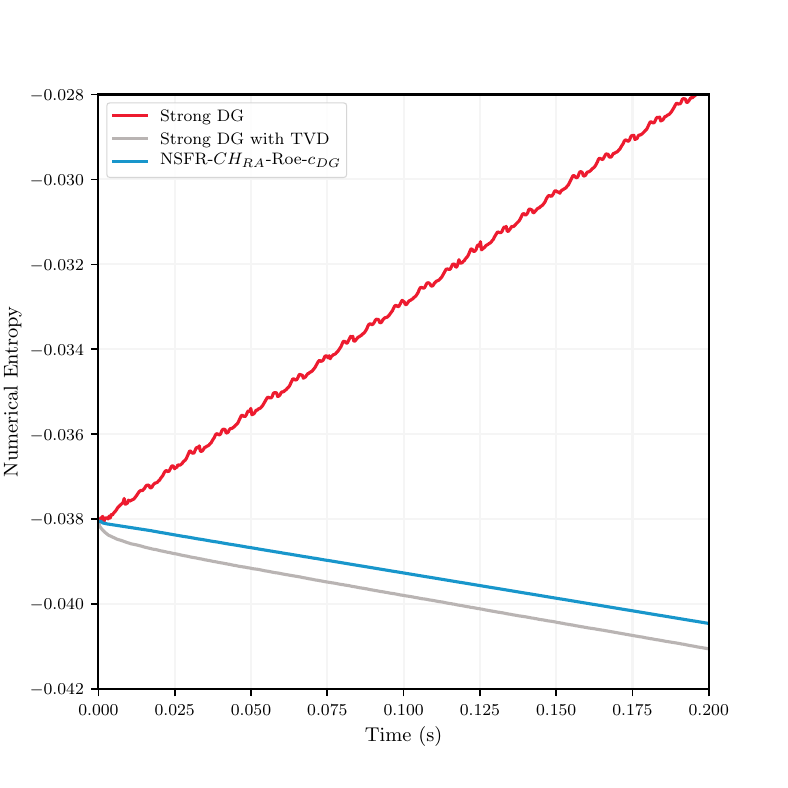}\\
\includegraphics[width=0.5\textwidth]{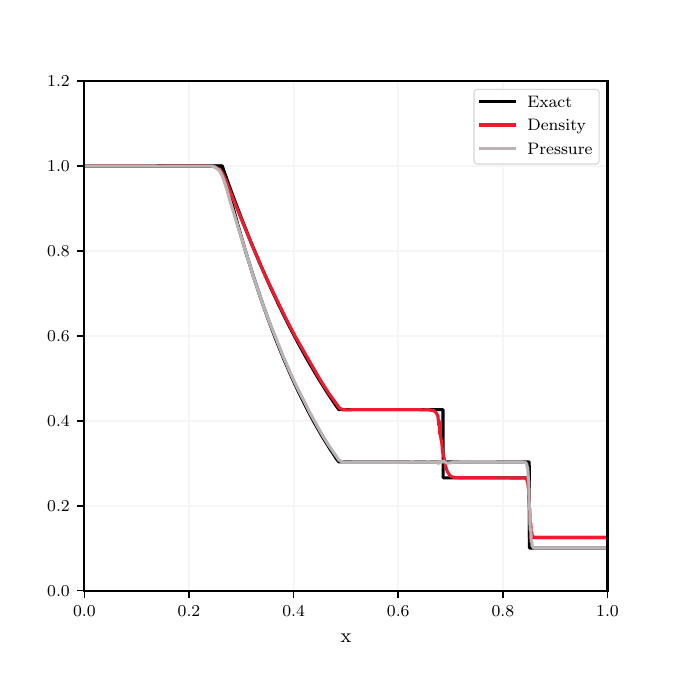}\hfill
\includegraphics[width=0.5\textwidth]{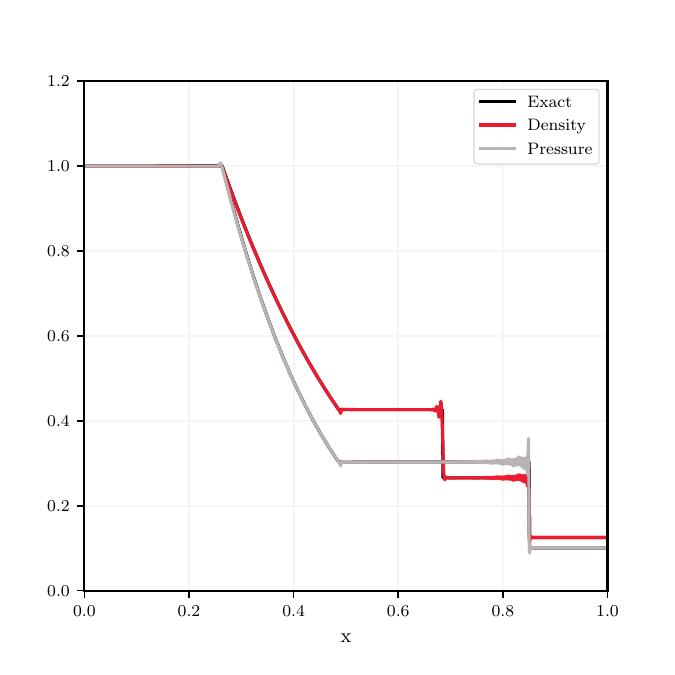}
\caption{\textit{[Sod Shock Tube]}  Entropy over time for initial three test configurations (top) and density and pressure with SSPRK3 Strong DG with TVD limiter (bottom left) and NSFR (bottom right), using $p = 3$, $CFL$ $0.5$ and the PPL at $t=0.2s$}
\label{fig:sodStrongDGComparison}
\end{figure}

The SSPRK3 Strong DG method with Roe flux and GLL nodes is used with the PPL to obtain a solution to the Sod Shock Tube test case. The $CFL$ number used to run this case is $0.5$. The case runs successfully without diverging; however, an entropy-violating shock is present in the solution. While the PPL preserves the positivity of density and pressure, the scheme itself does not prevent the presence of entropy-violating shocks. A plot of entropy over time is provided in Fig.~\ref{fig:sodStrongDGComparison}, wherein the entropy is shown to be increasing over time for this scheme.

To eliminate the entropy-violating shock, the TVD limiter introduced in \citet{chen2017entropy} is used with the PPL. This test configuration confirms that the implemented PPL can be used in conjunction with other limiting strategies. The $CFL$ number used to run the case is $0.5$. Referring back to Fig.~\ref{fig:sodStrongDGComparison}, we can see that the entropy decreases over time, and this is reflected in the density and pressure plots below. There is no entropy-violating shock present, and the density solution shows a left rarefaction wave, right shock wave and a contact discontinuity in the middle as expected. While the main features of the solution closely reflect that of the exact solution, the SSPRK3 Strong DG TVD solution contains some dissipative error at the head of the rarefaction wave and the contact discontinuity, with some dissipative error at the right shock as well. The solution does not accurately reflect the exact solution in these locations and is overly dissipative.

This test case is also run using the NSFR scheme to demonstrate the advantages of the scheme. As a result of the scheme's robustness, at a suitable $CFL$ number, the NSFR scheme can be used to solve this case without the PPL, while the SSPRK3 Strong DG method fails without the limiter when tested with a $CFL$ as low as $0.01$. The maximum $CFL$ number for which the NSFR scheme can be used without the PPL is $0.2$, which is 20 times greater than the lowest tested $CFL$ for the SSPRK3 Strong DG method. 

Using the PPL, the NSFR scheme can be used to run the test at a higher $CFL$. For this test, the $CFL$ number is set to $0.5$ so that the SSPRK3 Strong DG solutions can be compared to the NSFR solution. The solution obtained with the flux reconstruction parameter set to $c_{DG}$ is shown in Fig.~\ref{fig:sodStrongDGComparison}. The solution contains a left rarefaction wave and a right shock wave with a contact discontinuity in the middle. There is no entropy-violating shock present, and the shock transitions are sharp and are not overly dissipative. The NSFR solution does, however, contain some dispersive error in the wake of the shock. In terms of entropy, the NSFR scheme also has decreasing entropy over time as expected, but the decrease in entropy is less than that of the SSPRK3 Strong DG solution with the TVD limiter, suggesting that it is less dissipative.\\

\subsubsection{Polynomial Degree Comparison}\label{sec: sod_polydegree}
\begin{figure}[h!]
    \centering
    \includegraphics[width=.5\textwidth]{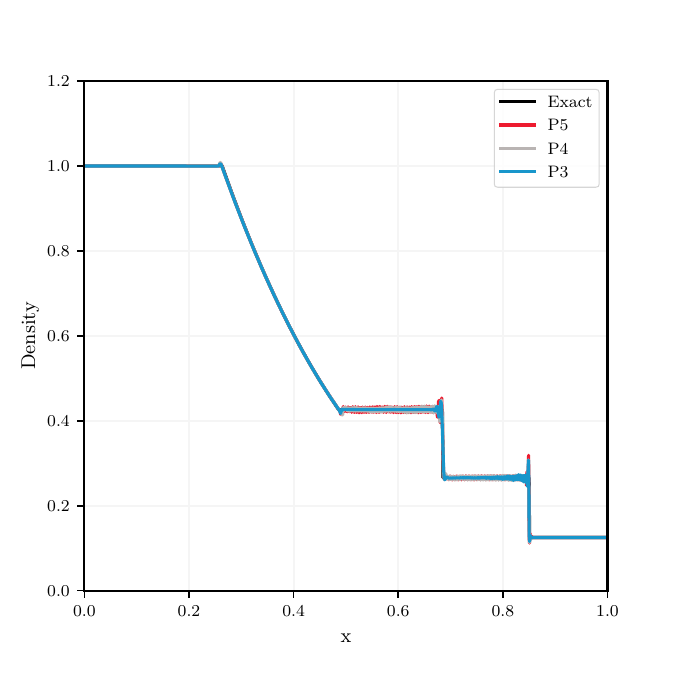}\hfill
    \includegraphics[width=.5\textwidth]{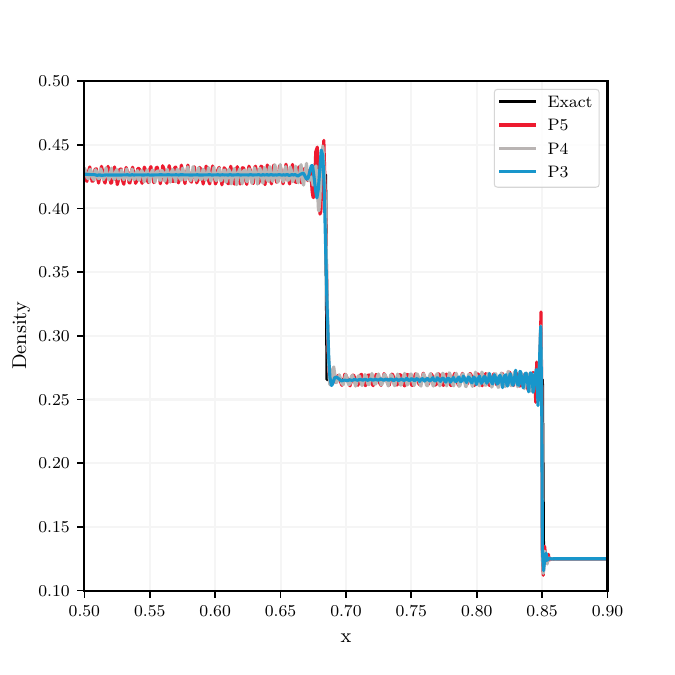}
    \caption{\textit{[Sod Shock Tube]} Density with NSFR-$CH_{RA}$, Roe dissipation, $c_{DG}$, $p = 3, 4, 5$, $CFL$ $0.5, 0.4, 0.3$ respectively with PPL at $t=0.2s$. Plots of full domain (left) and $0\leq x\leq 3.6$}
    \label{fig:sodShockNSFRPolyDegree}
\end{figure}

A comparison can also be conducted for the polynomial degree to demonstrate the ability of the scheme to remain stable without the added use of shock-capturing methods. The Sod Shock Tube case is run with polynomial degrees $p = 3, 4, 5$ with $CFL$ numbers $0.5, 0.4$ and $0.3$ respectively. The number of cells for each of the tests is $N=240,192,160$ respectively to ensure the tests have an equivalent number of degrees of freedom. The three solutions for the different polynomial degrees are compared in Fig.~\ref{fig:sodShockNSFRPolyDegree}. All three results closely follow the exact solution and have sharp shock transitions. The $p = 4$ and $p = 5$ results follow closely with the $p=3$ results but have a significant increase in oscillations in the wake of the shock as expected.

\subsubsection{CFL Condition for Two-Point Flux}\label{sec:tpf_sod}
This test is used to numerically verify that the two-point flux used for the tests in this investigation abides by a $CFL$ condition. The number of cells for this test is $N=500$. The test is run for polynomial degrees $p = 2,3,4$ with an increment of $0.0005s$ for the constant time step. The $CFL$ for a given configuration is calculated in two different approaches. The first method uses the minimum $\Delta x$ between the GLL quadrature nodes at the given degree for $\Delta x$. The second method uses the cell size for $\Delta x$. The $CFL$ is calculated as follows:
\begin{equation}\label{eq: cfl-formula}
    CFL = \frac{\Delta t}{\Delta x} \Bar{\lambda}_{max},
\end{equation}
where $\Bar{\lambda}_{max}$ denotes the maximum wavespeed across the entire domain and all time steps, calculated using the cell averaged solution vector, $\Bar{\mathbf{w}}$.
\begin{table}[]
\centering
\caption{Two-point flux positivity preservation test using the Sod Shock Tube setup with $N=500$.}
\label{tab:tpf_test_sodshock}
\resizebox{0.85\textwidth}{!}{
\begin{tabular}{rrrrrrr}
\hline
\multicolumn{7}{c}{\textbf{NSFR - $CH_{RA}$ - Roe Dissipation - $\mathbf{c_{DG}}$ - N = 500}}             \\ \hline
\multicolumn{1}{l}{k} &
  \multicolumn{1}{l}{Min. $\Delta x$} &
  \multicolumn{1}{l}{Cell Width} &
  \multicolumn{1}{l}{Time Step} &
  \multicolumn{1}{l}{Max. $\lambda$} &
  \multicolumn{1}{l}{Min. $\Delta x$ $CFL$} &
  \multicolumn{1}{l}{Cell Width $CFL$} \\ \hline
2 & 0.001   & 0.002 & 0.0021  & 2.32 & 0.95                         & 0.47                         \\
  &         &       & 0.00215 & 2.32 & 0.97                         & 0.49                         \\
  &         &       & 0.0022  & 2.32 & 0.99                         & 0.49                         \\
  &         &       & 0.00225 & 2.51 & \cellcolor[HTML]{FFA3A3}1.11 & \cellcolor[HTML]{FFA3A3}0.55 \\ \hline
3 & 0.00055 & 0.002 & 0.0014  & 2.26 & 0.86                         & 0.24                         \\
  &         &       & 0.00145 & 2.26 & 0.90                         & 0.25                         \\
  &         &       & 0.0015  & 2.27 & 0.95                         & 0.26                         \\
  &         &       & 0.00155 & 3.14 & \cellcolor[HTML]{FFA3A3}1.37 & \cellcolor[HTML]{FFA3A3}0.38 \\ \hline
4 & 0.00035 & 0.002 & 0.0008  & 2.23 & 0.71                         & 0.12                         \\
  &         &       & 0.00085 & 2.23 & 0.78                         & 0.13                         \\
  &         &       & 0.0009  & 2.24 & 0.84                         & 0.15                         \\
  &         &       & 0.00095 & 2.76 & \cellcolor[HTML]{FFA3A3}1.12 & \cellcolor[HTML]{FFA3A3}0.19 \\ \hline
\end{tabular}
}
\end{table}

For the first method of calculating the $CFL$, we observe that the test only fails when the $CFL$ exceeds $1$. This follows closely with the $CFL$ condition given in Eq.\ref{eq: cfl-condition-1}, thus numerically verifying that the sufficient condition for the $CFL$ holds for the two-point flux as well. For the second method, it is observed for $p=2$ that the test fails when the $CFL$ exceeds $0.5$, which follows closely with the observed condition given in Eq.\ref{eq: cfl-condition-2}. It is noted, however, that for $p=3$ and $4$, the $CFL$ at failure is lower than the expected condition as the change in the minimum $\Delta x$ between the different polynomial degrees is not negligible. This case numerically shows that Eq. \ref{eq: cfl-condition-1} is a sufficient condition for positivity preservation for two-point fluxes, irrespective of the polynomial degree for the current case.

\subsubsection{Two-Point Flux Comparison}

This test is also used to compare the impact of the two-point flux on the solution. In addition to the $\left(\text{CH}_{\text{RA}}\right)$ flux used throughout this investigation, the following fluxes are considered:
\begin{table}[H]
\centering
\caption{Summary of two-point fluxes}
\label{tab:2ptflux}
\begin{tabular}{@{}lrrl@{}}
\toprule
\textbf{Two-Point Flux}                 & \textbf{EC} & \multicolumn{1}{l}{\textbf{KEP}} & \textbf{PEP}                  \\ \midrule
Ismail-Roe $\left(\text{IR}\right)$\cite{ismail2009affordable}    & \checkmark   & \multicolumn{1}{l}{}             &                               \\
Kennedy-Gruber $\left(\text{KG}\right)$\cite{KENNEDY20081676} &             & \checkmark                        &                               \\
Chandrashekar $\left(\text{CH}\right)$\cite{chandrashekar2013kinetic}  & \checkmark   & \checkmark                        &                               \\
$\text{CH}_{\text{RA}}$\cite{ranocha2021preventing}    & \checkmark   & \checkmark                        & \multicolumn{1}{r}{\checkmark} \\ \bottomrule
\end{tabular}
\end{table}

The properties of the different two-point fluxes are highlighted in Table~\ref{tab:2ptflux}. The computational cost for all the fluxes is similar with the exception of the $\left(\text{KG}\right)$ flux which does not require logarithmic averaging (Eq.~\ref{eq:log_mean}) thereby incurring a lower cost than the other fluxes. Despite the $KG$ incurring a lower cost, it is not preferred for problems involving strong shocks. The issues with the $KG$ flux are highlighted in Section~\ref{sec: 2DSVSW}.

\begin{figure}[h!]
    \centering
    \includegraphics[width=.5\textwidth]{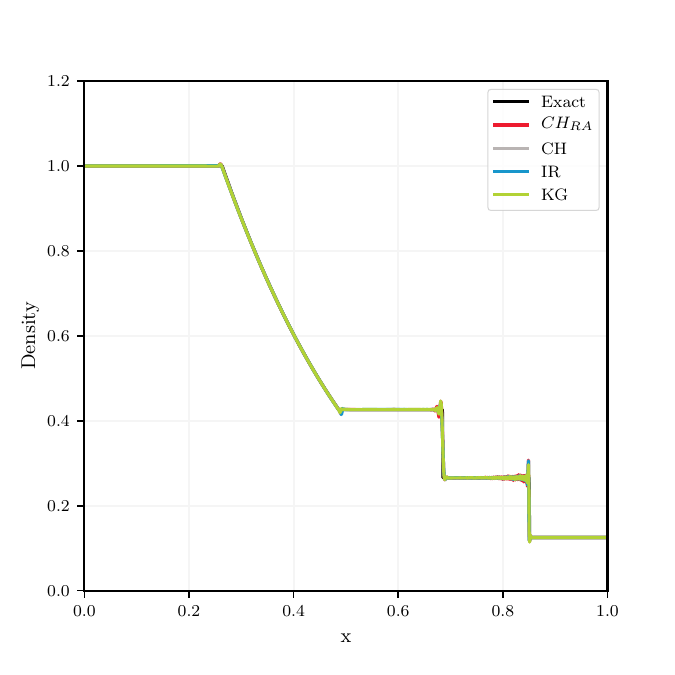}\hfill
    \includegraphics[width=.5\textwidth]{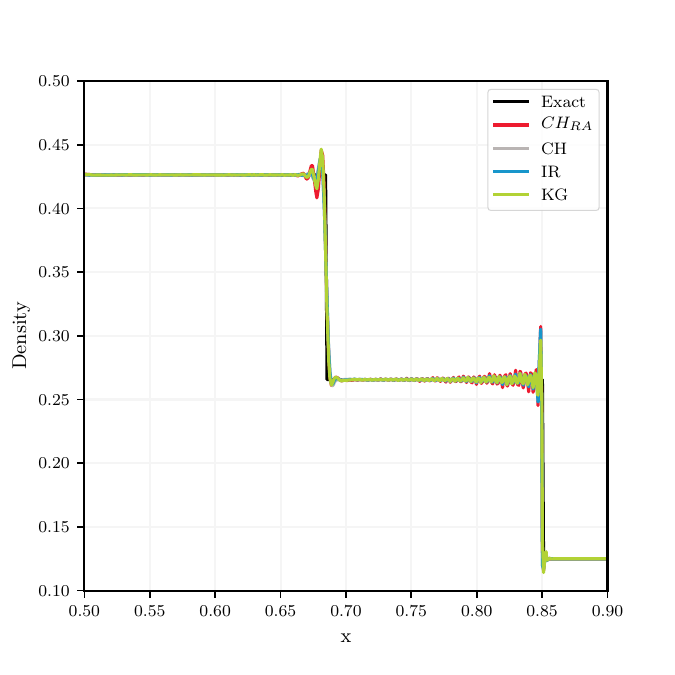}
    \caption{\textit{[Sod Shock Tube]} Density and pressure with NSFR, Roe dissipation, $c_{DG}$, $p = 3$, $CFL$ $0.5$ with PPL at $t=0.2s$ for $CH_{RA}$, CH, IR and KG}
    \label{fig:sodShockNSFR_2ptFlux}
\end{figure}

The solution for each two-point flux is shown in Fig.~\ref{fig:sodShockNSFR_2ptFlux}. The results for the different two-point fluxes follow closely with the exact solution. The shock location and magnitude are correct and the shock transitions are sharp. It is observed that the $CH_{RA}$ solution is slightly more oscillatory in the wake of the shock compared to the other 3 fluxes. Additionally, the magnitude of the overshoot at the shock is slightly greater for the $CH_{RA}$ and $IR$ fluxes compared to $CH$ and $KG$. It is also noted that while the $KG$ flux is less oscillatory in the wake of the shocks, it exhibits the greatest undershoot at the foot of the rarefaction wave out of all the fluxes.

\subsubsection{Investigation of Flux Reconstruction Schemes}
\begin{figure}[h!]
    \centering
    \includegraphics[width=.5\textwidth]{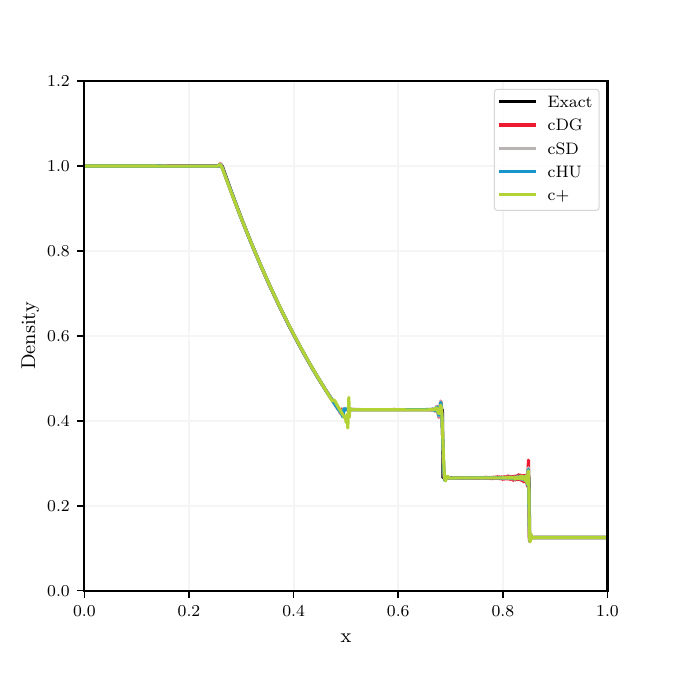}\hfill
    \includegraphics[width=.5\textwidth]{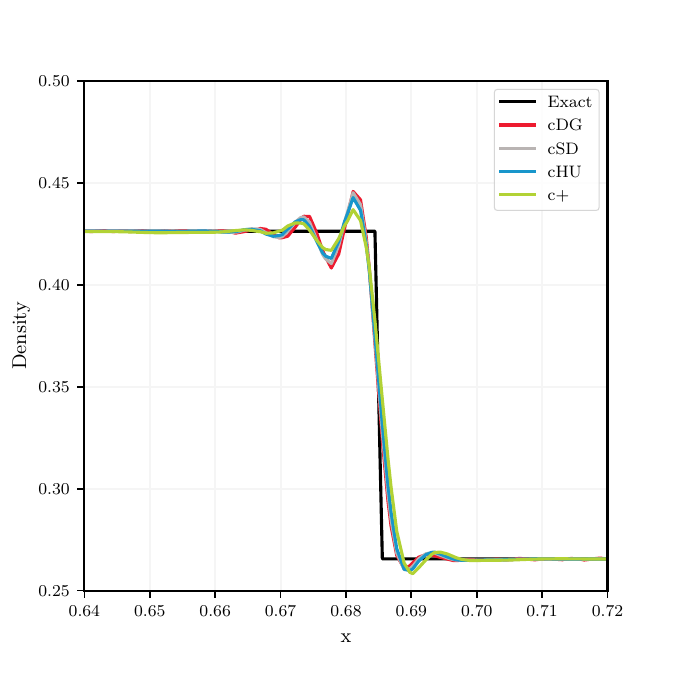}
    \caption{\textit{[Sod Shock Tube]} Density with NSFR-$CH_{RA}$, Roe dissipation, $p = 3$, $CFL$ $0.5$ respectively with PPL at $t=0.2s$ for different flux reconstruction parameters. Plots of full domain (left) and $0\leq x\leq 3.6$}
    \label{fig:sodShockFRParam}
\end{figure}

Here, we investigate variants of the flux reconstruction scheme. The Sod Shock Tube case is run with the $c_{SD}$, $c_{HU}$, and $c_{+}$ schemes with $p = 3$ and a $CFL$ of $0.5$. Using the solution obtained with the $c_{DG}$ scheme in Fig.~\ref{fig:sodStrongDGComparison}, the four different schemes are compared in Fig.~\ref{fig:sodShockFRParam}. The results for all four schemes follow closely with the exact solution, and the shocks are in the correct locations and have sharp transitions. As the value of the flux reconstruction parameter is increased from $c_{DG}$ to $c_{+}$, it is observed that the overshoot at the shock and the oscillations in the wake of the shock are diminished; which has the effect of filtering the highest mode of the solution. However, while filtering of nonphysical behaviour is an advantage, it is also observed that as the FR parameter is increased, there is a noticeable undershoot at the foot of the rarefaction wave. 

Overall, the Sod Shock Tube case demonstrates the NSFR method's ability to simulate problems involving shocks without the need for a TVB/TVD limiter and the added robustness of the method also allows for a higher $CFL$. Even when parameters such as the polynomial degree and FR parameter are varied, the method still consistently provides a solution that accurately reflects the exact solution without the need for TVB/TVD limiting. The test also sufficiently verifies the CFL condition for two-point fluxes numerically.

\subsection{Shu Osher Problem}
\label{sec: 1DShuOsherProblem}
The Shu Osher problem is a 1D Euler case with a shock propagating into an unsteady density field. This problem essentially simulates shock-turbulence interaction in 1D. The problem tests the ability of a scheme to capture a shock, the interaction of the shock with an unsteady density field and the waves that propagate downstream of the shock. The computational domain of the problem is $[-5,5]$. The problem is initialized as,
\begin{equation}
    (\rho, w, p) = \begin{cases}
        (3.857143,2.629369,10.33333) \qquad &\text{if} \: x < -4,\\
        (1+0.2sin(5x),0,1) \qquad &\text{if} \: x \geq -4.
    \end{cases}
\end{equation}
The final time for this case is $t = 1.8s$. The number of cells is $N = 128$.  As outlined in \cite{JOHNSEN20101213}, at time $1.8s$, the shock location is $x_s\approx2.39$, the contact discontinuity location is $x_c\approx0.69$, and the location of the leading acoustic wave is $x_a\approx-2.75$.  The results are compared to a reference solution obtained using the third-order SSPRK3 Strong DG method with Lax-Friedrichs flux, GLL nodes, polynomial degree of $p=0$ and $N=32768$.

\begin{figure}
    \centering
    \includegraphics[width=.5\textwidth]{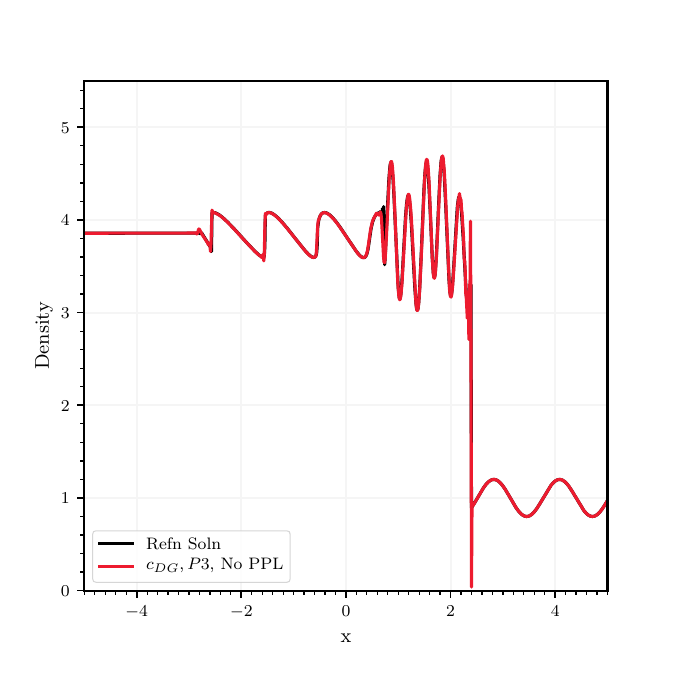}\hfill
    \includegraphics[width=.5\textwidth]{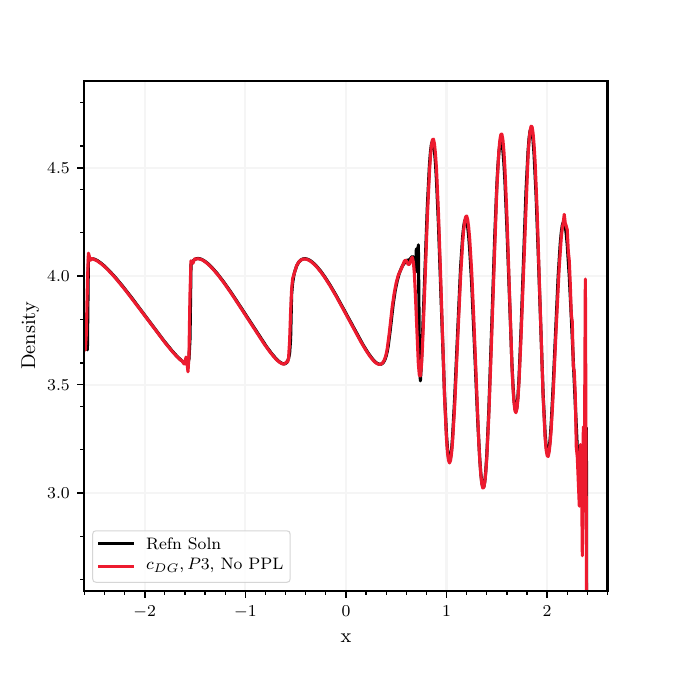}
    \caption{\textit{[Shu Osher Problem]} Density solution with NSFR-$CH_{RA}$, Roe dissipation, $c_{DG}$, no PPL at $t=1.8s$. Plots of full domain (left) and $-2.6\leq x \leq 2.6$ (right).}
    \label{fig:ShuOsherNSFR_noPPL}
\end{figure}

This test is first run using the NSFR scheme without the PPL to demonstrate its robustness. The $CFL$ number for this case is $0.01$ and the results are given in Fig.~\ref{fig:ShuOsherNSFR_noPPL}. This test cannot be run with the Strong DG scheme in the same configuration. There is a significant overshoot at the shock, but otherwise, the solution accurately reflects all the expected features, including the shock, contact discontinuity and acoustic waves.  

\subsubsection{CFL Condition for Two-Point Flux}\label{sec:tpf_shuosher}
This test is also used to numerically verify that the two-point flux used throughout this study abides by a $CFL$ condition. The number of cells for this test is $N=360$. The test is run for polynomial degrees $p = 2,3,4$ with an increment of $0.00005s$ for the constant time step. Once again, the $CFL$ for a given configuration is calculated in two different ways. The first method uses the minimum $\Delta x$ between the GLL quadrature nodes at the given degree for $\Delta x$. The second method uses the cell size for $\Delta x$. The $CFL$ is calculated using Eq.\ref{eq: cfl-formula}. The results of the numerical study are given in Table \ref{tab:tpf_test_shuosher}.

\begin{table}[]
\centering
\caption{Two-point flux positivity preservation test using the Shu Osher Problem setup with $N=360$.}
\label{tab:tpf_test_shuosher}
\resizebox{0.7\textwidth}{!}{
\begin{tabular}{rrrrrrr}
\hline
\multicolumn{7}{c}{\textbf{NSFR - $CH_{RA}$ - Roe Dissipation - $\mathbf{c_{DG}}$ - N = 360}}               \\ \hline
\multicolumn{1}{l}{k} &
  \multicolumn{1}{l}{Min. $\Delta x$} &
  \multicolumn{1}{l}{Cell Width} &
  \multicolumn{1}{l}{Time Step} &
  \multicolumn{1}{l}{Max. $\lambda$} &
  \multicolumn{1}{l}{Min. $\Delta x$ CFL} &
  \multicolumn{1}{l}{Cell Width CFL} \\ \hline
2 & 0.01389 & 0.02778 & 0.0021  & 5.45 & 0.82                         & 0.41                         \\
  &         &         & 0.00215 & 5.64 & 0.87                         & 0.44                         \\
  &         &         & 0.0022  & 5.70 & 0.90                         & 0.45                         \\
  &         &         & 0.00225 & 7.30 & \cellcolor[HTML]{FFA3A3}1.18 & \cellcolor[HTML]{FFA3A3}0.59 \\ \hline
3 & 0.00768 & 0.02778 & 0.0014  & 4.97 & 0.91                         & 0.25                         \\
  &         &         & 0.00145 & 4.99 & 0.94                         & 0.26                         \\
  &         &         & 0.0015  & 4.97 & 0.97                         & 0.27                         \\
  &         &         & 0.00155 & 5.95 & \cellcolor[HTML]{FFA3A3}1.20 & \cellcolor[HTML]{FFA3A3}0.33 \\ \hline
4 & 0.00479 & 0.02778 & 0.0008  & 4.90 & 0.82                         & 0.14                         \\
  &         &         & 0.00085 & 4.90 & 0.87                         & 0.15                         \\
  &         &         & 0.0009  & 4.93 & 0.93                         & 0.16                         \\
  &         &         & 0.00095 & 4.90 & \cellcolor[HTML]{FFA3A3}0.97 & \cellcolor[HTML]{FFA3A3}0.17 \\ \hline
\end{tabular}
}
\end{table}

Similar to the results obtained using Sod Shock Tube, this test also demonstrates that the two-point flux abides by the condition given in Eq.\ref{eq: cfl-condition-1} for the first method. The $p=4$ case fails for a $CFL$ slightly less than $1$, but it is sufficiently close enough to the bound that it is considered to satisfy the condition. The results also show that the two-point flux abides by Eq.\ref{eq: cfl-condition-2} for a polynomial degree of 2. This test case further establishes Eq. \ref{eq: cfl-condition-1} as a sufficient condition for positivity preservation of pressure for two-point fluxes irrespective of polynomial degree.

\subsubsection{Investigation of Flux Reconstruction Schemes}
\begin{figure}
    \centering
    \includegraphics[width=.5\textwidth]{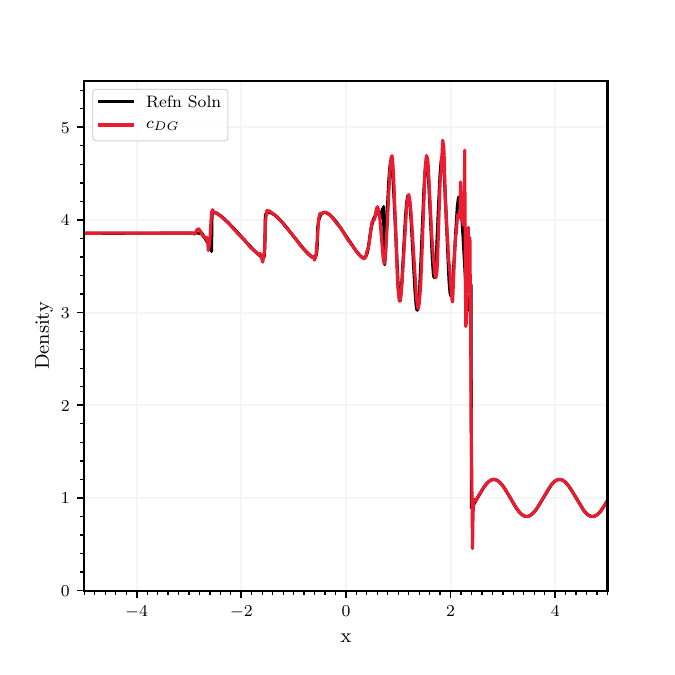}\hfill
    \includegraphics[width=.5\textwidth]{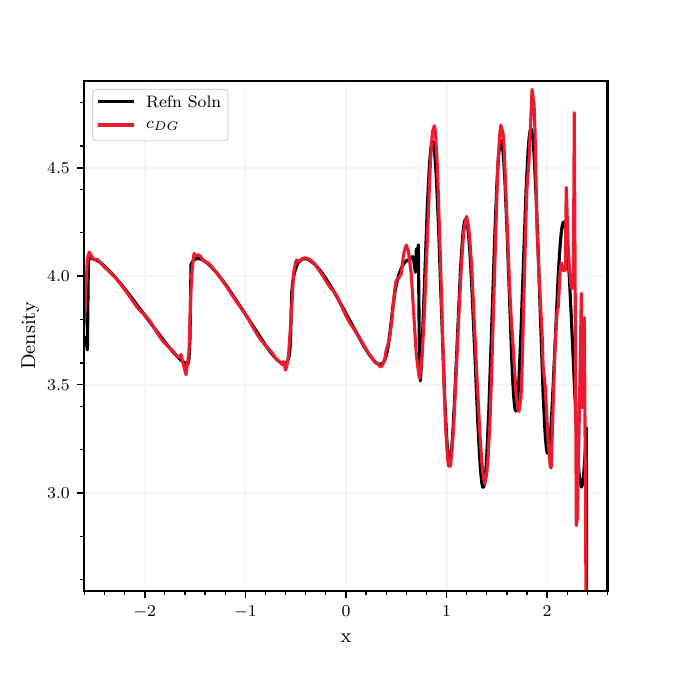}
    \caption{\textit{[Shu Osher Problem]} Density solution with NSFR-$CH_{RA}$, Roe dissipation, $c_{DG}$, PPL at $t=1.8s$. Plots of full domain (left) and $-2.6\leq x \leq 2.6$ (right).}
    \label{fig:ShuOsherNSFRcDG}
\end{figure}

This test is also run using the $c_{DG}$ and $c_+$ schemes. The $CFL$ number used for both configurations is 0.5. This increased $CFL$ is a good indicator of the added robustness the PPL provides.
The results obtained with the $c_{DG}$ scheme are shown in Fig.~\ref{fig:ShuOsherNSFRcDG}.The $c_{DG}$ results capture all the locations of the shock, contact discontinuity and acoustic waves accurately, similar to the results in Fig.~\ref{fig:ShuOsherNSFR_noPPL} despite having a $CFL$ that is 50 times greater. Comparing the limited solution to the reference solution, however, it is evident that there are overshoots in the vicinity of the shock, and there are oscillations present immediately upstream of the shock.

\begin{figure}[ht]
    \centering
    \includegraphics[width=.5\textwidth]{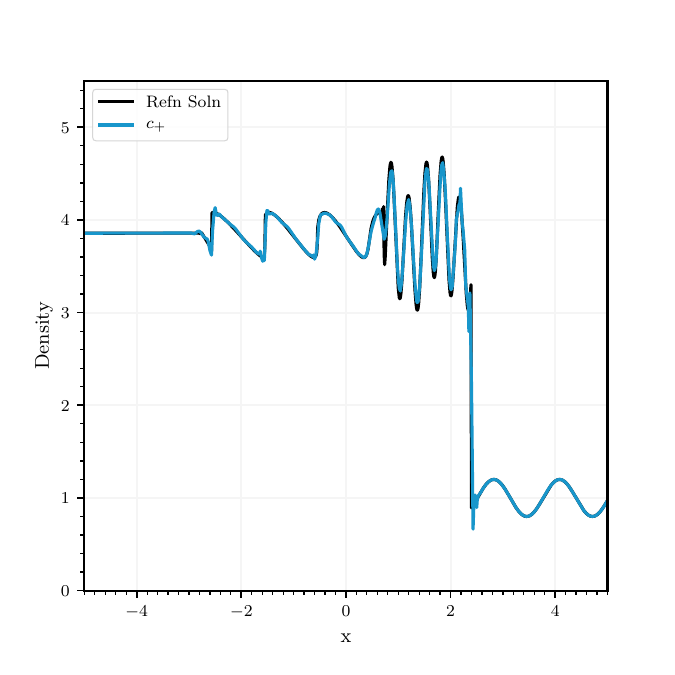}\hfill
    \includegraphics[width=.5\textwidth]{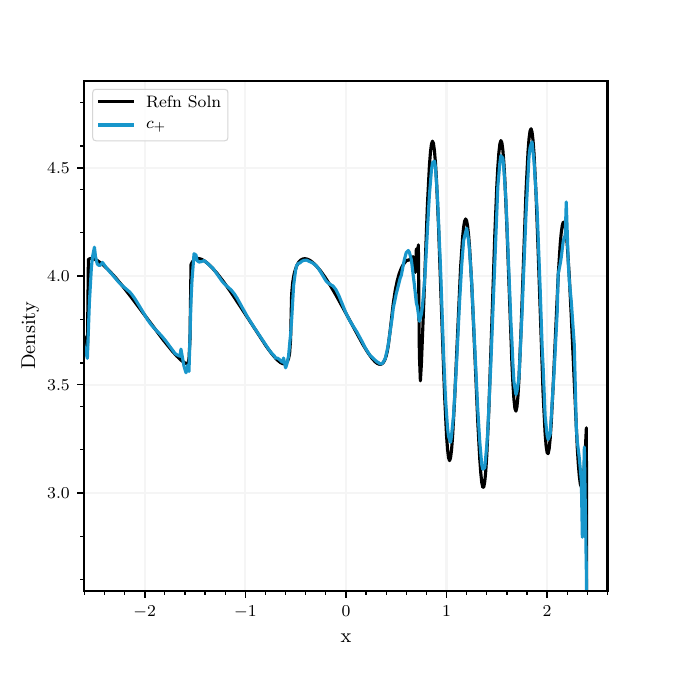}
    \caption{\textit{[Shu Osher Problem]} Density solution with NSFR-$CH_{RA}$, Roe dissipation, $c_{+}$, PPL at $t=1.8s$. Plots of full domain (left) and $-2.6\leq x \leq 2.6$ (right).}
    \label{fig:ShuOsherNSFRc+}
\end{figure}

The results obtained using the $c_{+}$ scheme are shown in Fig.~\ref{fig:ShuOsherNSFRc+}. The plot shows that the shock, contact discontinuity, and leading acoustic wave are all accurately captured and the shock does not exhibit an overshoot as in the $c_{DG}$ scheme and the oscillations in the wake of the shock have been dampened. However, comparing the $c_+$ solution to the reference solution, it is apparent that there is room for improvement as the dissipation added to dampen the oscillations and overshoot also diminish the magnitude of the waves upstream of the shock. 

The plot of the results between $-2.6\leq x\leq 2.6$ shows the mitigation of the overshoot and oscillations in the $c_{+}$ results, but it also highlights the diminished magnitude of the waves downstream of the shock which is better captured by $c_{DG}$. Using $c_{+}$ to filter out the highest mode produces fewer oscillatory results but leads to diminished peaks. An adaptive application of the flux reconstruction parameter wherein FR, $c_{+}$, is used in shock locations while discontinuous Galerkin, $c_{DG}$, is used elsewhere could lead to improved results that achieve the desired magnitudes while reducing overshoots and oscillations in the presence of shocks and would be part of a future study.

\subsubsection{Impact of Flux Nodes}
\begin{figure}[]
    \centering
    \includegraphics[width=.42\textwidth]{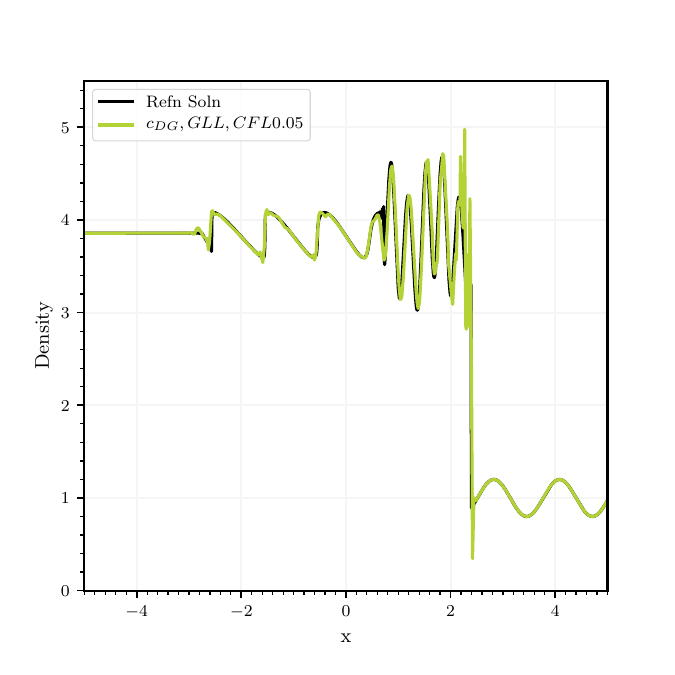}
    \includegraphics[width=.42\textwidth]{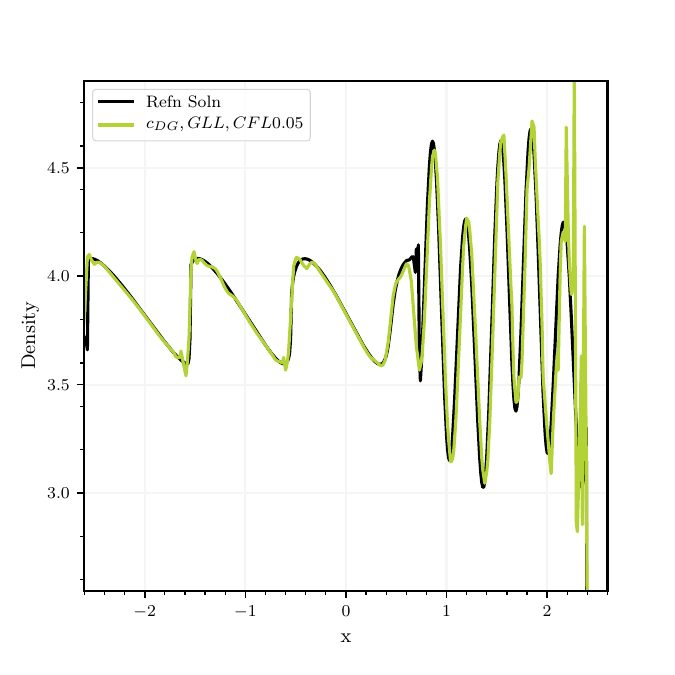}\\
    \includegraphics[width=.42\textwidth]{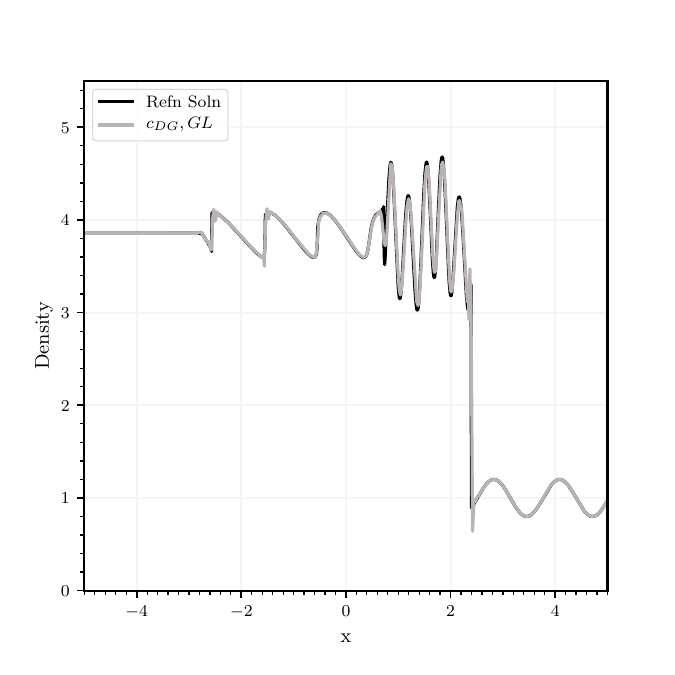}
    \includegraphics[width=.42\textwidth]{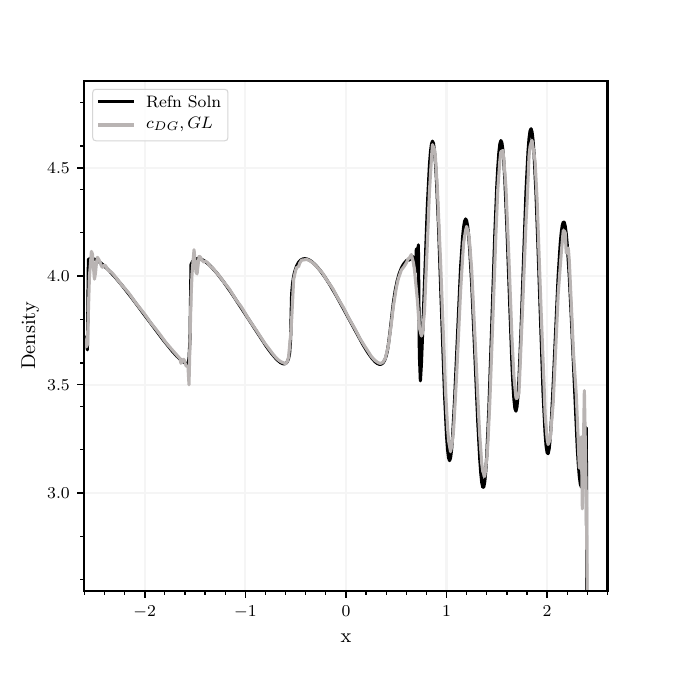}
    \caption{\textit{[Shu Osher Problem]} Density solution with NSFR-$CH_{RA}$, Roe dissipation, $c_{DG}$, PPL at $t=1.8s$. Plots of GLL nodes (top) and GL nodes (bottom). Plots of full domain (left) and $-2.6\leq x \leq 2.6$ (right).}
    \label{fig:ShuOsherNSFRGL}
\end{figure}

This test case is used to demonstrate the impact of the choice of flux nodes. The test is run once again using the same parameters as the $c_{DG}$ scheme, with the exception of the flux nodes, which are set to GL. As a result of GL quadrature being used, the $CFL$ is also subsequently lowered. As discussed in \citet{chan2018discretely}, for quadratures such as GL, the solution can spike when evaluating the entropy-projected conservative variables at surface points, and to mitigate this effect, the $CFL$ is set to 0.05. The test is also rerun with the lowered $CFL$ for GLL flux nodes to facilitate a better comparison of the results. As shown in Fig.~\ref{fig:ShuOsherNSFRGL}, the results for the GL quadrature are more accurate than the GLL results. There is less overshoot in the wake of the shock, and there is less oscillatory behaviour, but it comes at the cost of the $CFL$ being lowered. The test can be run with a significantly higher $CFL$ using GLL flux nodes, and the issues of overshoot and oscillations can be addressed through the flux reconstruction parameter as shown by the $c_{+}$ scheme, thereby making GLL flux nodes preferable over GL flux nodes.

\subsubsection{Polynomial Degree Comparison}
\begin{figure}[ht]
    \centering
    \includegraphics[width=.5\textwidth]{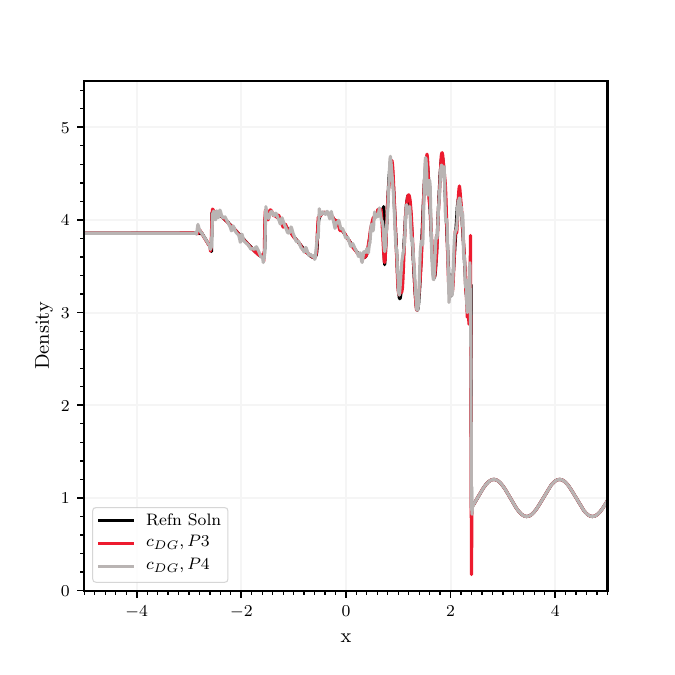}\hfill
    \includegraphics[width=.5\textwidth]{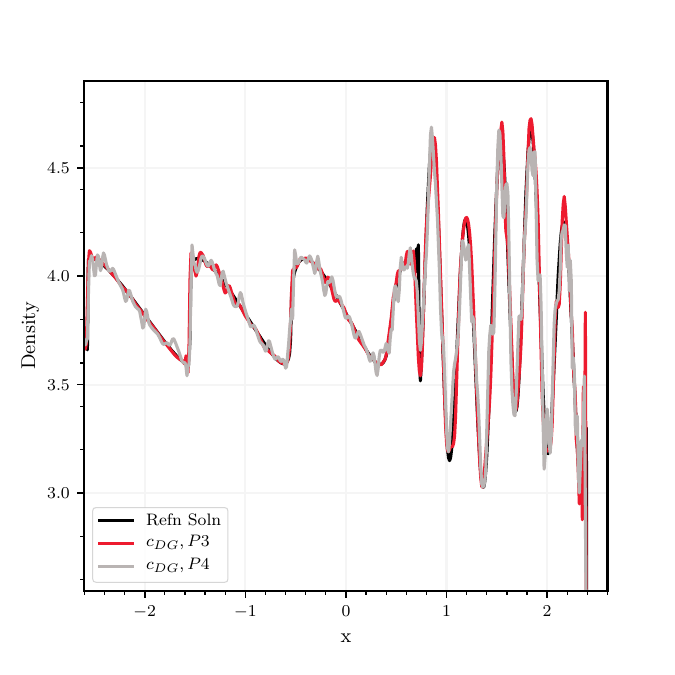}
    \caption{\textit{[Shu Osher Problem]} Density solution with NSFR-$CH_{RA}$, Roe dissipation, $c_{DG}$, PPL at $t=1.8s$ for $p=3,4$. Plots of full domain (left) and $-2.6\leq x \leq 2.6$ (right).}
    \label{fig:ShuOsherNSFR_polydegree_34}
\end{figure}

A polynomial degree comparison is also presented for the Shu Osher test case. The test is run for $p=3$ and $p=4$ at a $CFL$ of $0.5$. The results for polynomial degrees $3$ and $4$ are shown in Fig.~\ref{fig:ShuOsherNSFR_polydegree_34}. Both solutions closely follow the reference solution but the higher polynomial degree solution contains more oscillations in the waves downstream of the shock similar to the Sod Shock Tube polynomial degree comparison. However, unlike the Sod Shock Tube results, the higher polynomial degrees display less overshoot at the shock than the $p=3$ results. 

\subsection{Leblanc Shock Tube}
\label{sec: 1DLeblancShockTest}
The modified version of the Leblanc Shock Tube case seen in \citet{ZHANG20108918} is used to verify the 1D implementation of the PPL. Due to the extreme nature of the shock, this case is a good test of the robustness of the NSFR scheme in conjunction with the PPL. This case is a Riemann problem of 1D Euler equations. The computational domain is $[-10.0, 10.0]$. The problem is initialized as,
\begin{equation}
    (\rho, w, p) = \begin{cases}
        (2,0,10^{9}) \qquad &\text{if} \: x < 0,\\
        (0.001,0,1) \qquad &\text{if} \: x \geq 0.
    \end{cases}
\end{equation}
The final time for this case is $t = 1\times 10^{-4}$. The number of cells is $N = 512$.\\

To further highlight the advantages of the NSFR method, it is noted that to run this case using the SSPRK3 Strong DG method, both the PPL and TVD limiters are required. The resulting solution is similar to those seen in Sec. \ref{sec:sod_strongDG}, where the use of the TVD limiter leads to dissipative error in the vicinity of the shock and contact discontinuity, which does not accurately reflect that which is seen in literature. 

\begin{figure}
    \centering
    \includegraphics[width=.5\textwidth]{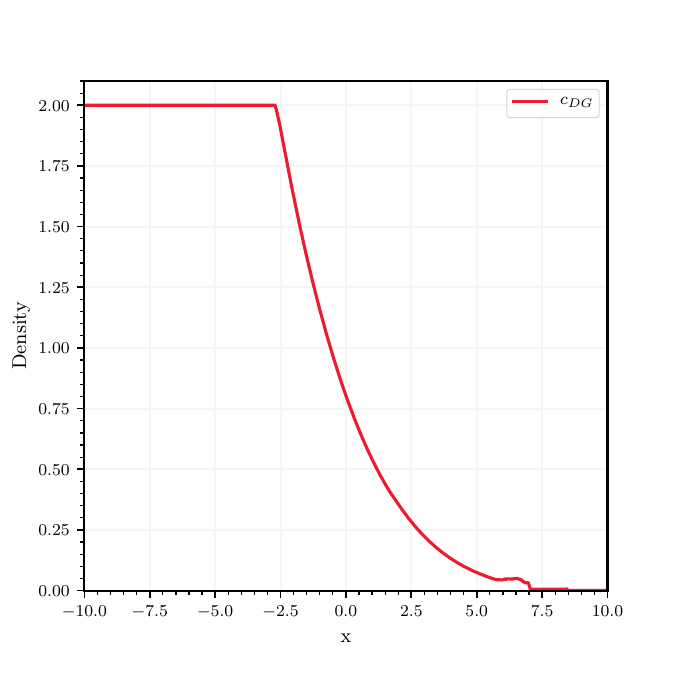}\hfill
    \includegraphics[width=.5\textwidth]{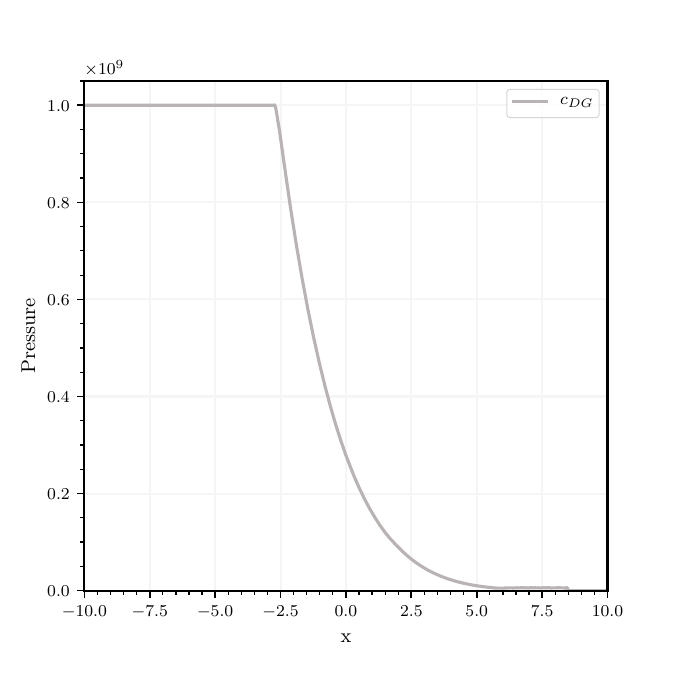}
    \caption{\textit{[Leblanc Shock Tube]} Density (left) and pressure (right) with NSFR-$CH_{RA}$, Roe dissipation, $c_{DG}$, $p = 3$, $CFL$ $0.06$ and PPL at $t=0.0001s$}
    \label{fig:LeblancShockNSFRcDG}
\end{figure}

The $CFL$ number used to run this case is $0.06$ to best compare the results with the $c_{DG}$ scheme. The $c_{DG}$ scheme is used to run the Modified 1D Leblanc Shock Tube case. The maximum $CFL$ number for this case using $c_{DG}$ is 0.06. The solution is shown in Fig.~\ref{fig:LeblancShockNSFRcDG}. The result captures the location and magnitude of the shocks very well. The solution follows closely with results seen in \citet{ZHANG20108918} and \citet{xu2022third}. This test case shows that the PPL can be used in conjunction with NSFR schemes to ensure that the positivity preservation property is satisfied.\\

\subsubsection{Investigation of Flux Reconstruction Schemes}
\begin{figure}
    \centering
    \includegraphics[width=.5\textwidth]{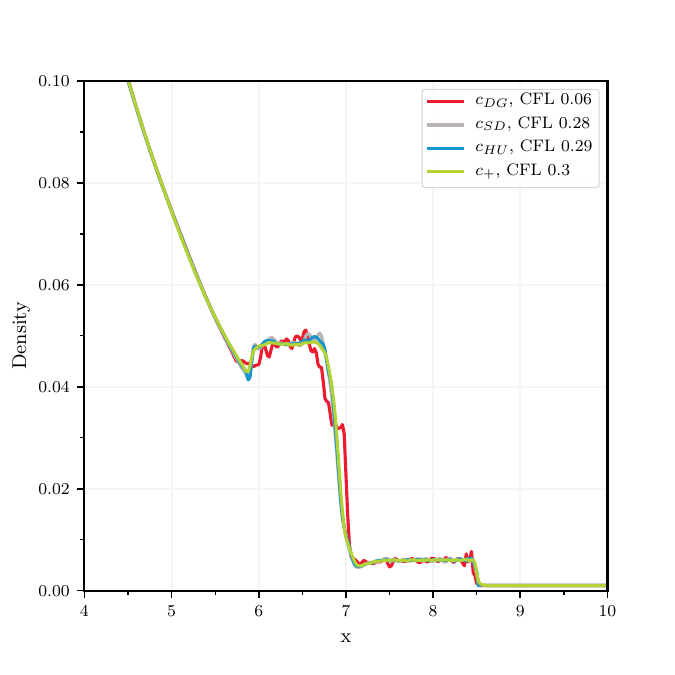}\hfill
    \includegraphics[width=.5\textwidth]{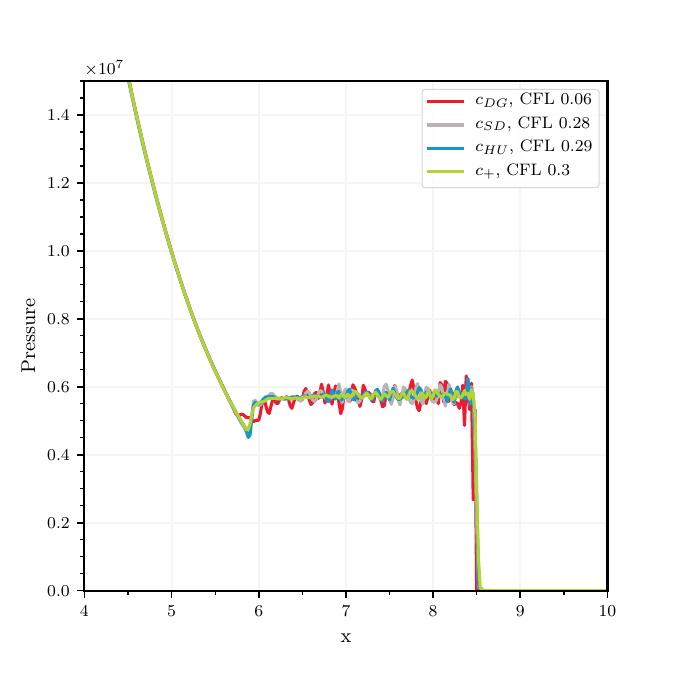}
    \caption{\textit{[Leblanc Shock Tube]} Density (left) and pressure (right) with NSFR-$CH_{RA}$, Roe dissipation, $p = 3$, using the PPL, comparing $c_{DG}$, $c_{SD}$, $c_{HU}$, $c_{+}$ for $4 \leq x \leq 10$ at $t=0.2s$}
    \label{fig:LeblancShockNSFRComparison}
\end{figure}

To demonstrate the impact of variants of flux reconstruction schemes, the case is also run using $c_{SD}$, $c_{HU}$, and $c_{+}$. The maximum $CFL$ number for each of the flux reconstruction parameters is $0.28,0.29,0.3$ respectively. The three schemes have similar density and pressure profiles, so to compare the impact of the flux reconstruction parameter, the solution in the wake of the shock ($5 \leq x \leq 8.75$) is presented for all the schemes in Fig.~\ref{fig:LeblancShockNSFRComparison}. When comparing the different solutions, it can be seen that the oscillations in the wake of the shock are dampened for both density and pressure when $c>0$ despite being run with a $CFL$ that is nearly five times greater. As the value of $c$ increases, the oscillations become increasingly dampened, while preserving the magnitude and location of the shock. While the increase in the value of $c$ can dampen the oscillations, it also results in a slight undershoot at the foot of the rarefaction wave similar to the Sod Shock Tube results.\\ 

\begin{figure}
    \centering
    \includegraphics[width=.5\textwidth]{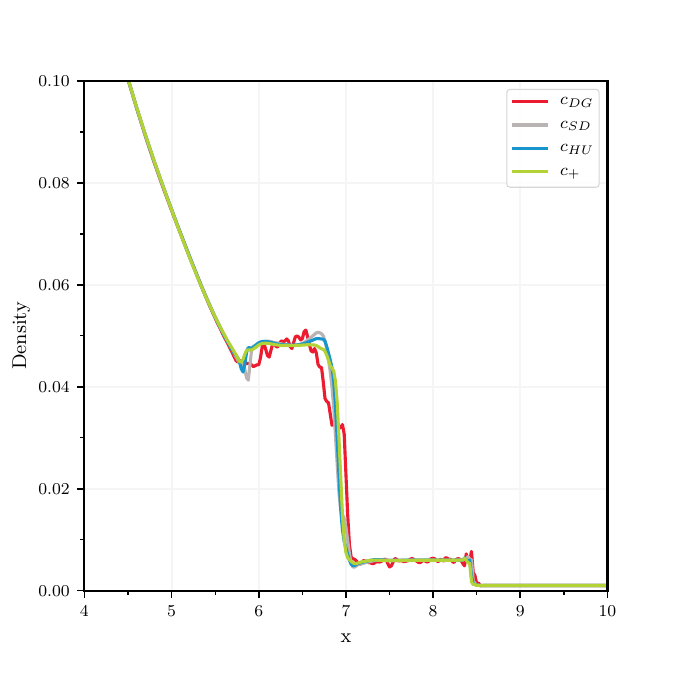}\hfill
    \includegraphics[width=.5\textwidth]{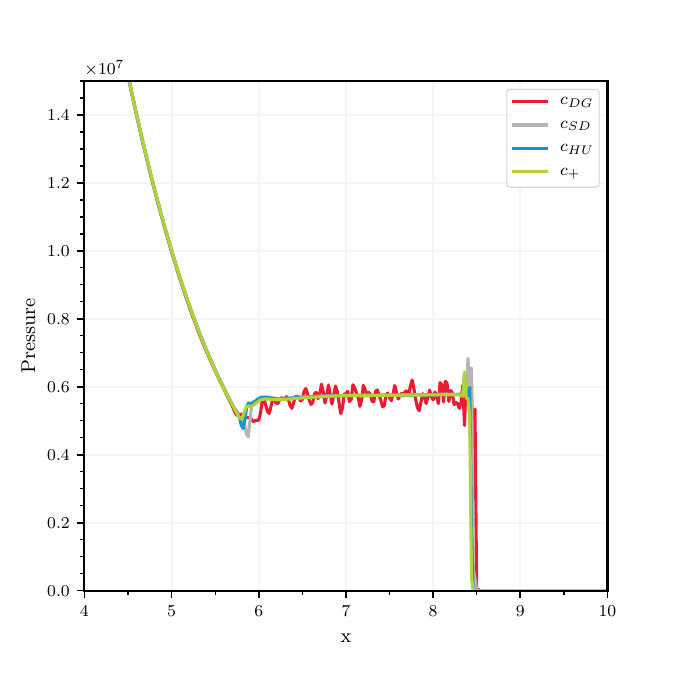}
    \caption{\textit{[Leblanc Shock Tube]} Density (left) and pressure (right) with NSFR-$CH_{RA}$, Roe dissipation, $p = 3$, using the PPL, comparing $c_{DG}$, $c_{SD}$, $c_{HU}$, $c_{+}$ for $4 \leq x \leq 10$ with $CFL=0.06$ at $t=0.0001s$}
    \label{fig:LeblancShockNSFRComparison_cfl0.06}
\end{figure}

The previous test configuration shows the $CFL$ advantage associated with an increase in the value of $c$, but the results also exhibit better shock capturing capabilities as the oscillations in the vicinity of the shock are mitigated. To further show the added advantage of shock capturing, the tests for the different flux reconstruction parameters are run again with a $CFL$ of $0.06$ for all values of $c$. The results are shown in Fig.~\ref{fig:LeblancShockNSFRComparison_cfl0.06}. While we still observe overshoot at the shock and undershoot at the foot of the rarefaction wave for the FR parameters that are greater than zero, we also see that the oscillations previously seen in the wake of the shock are almost entirely eliminated at a lower $CFL$. This highlights the effectiveness of flux reconstruction schemes for shock capturing and demonstrates the potential of an adaptive approach, where the FR parameter is tuned to maximize timestep efficiency while suppressing oscillations.
\break
\subsection{Strong Vortex Shock Wave Interaction}\label{sec: 2DSVSW}
The implementation of the PPL in two dimensions is verified using the Strong Vortex Shock Wave (SVSW) interaction case \cite{rault2003shock} which is a benchmark case from the 5th international workshop on high-order CFD methods for the American Institute of Aeronautics and Astronautics conference led by \citet{HOW5}. This case can be run using the SSPRK3 Strong DG method with the limiter; however, the results contain severe oscillations and are not consistent with published results. To run this case with the SSPRK3 Strong DG method, both the PPL and a TVD limiter would be required. The NSFR scheme, however, can be used to run this test even without the application of the PPL. This allows us to compare the solutions obtained with and without the PPL to further confirm that the limiter does not affect the accuracy of the underlying scheme.

The computational domain is $[0,2]\times[0,1]$ with a grid of $N_x\times N_y = 300\times150$. The initial condition of the problem involves a stationary shock ($M_s = 1.5$) at $x=0.5$ and a strong vortex ($M_v = 0.9$) centered at $(x_c,y_c) = (0.25,0.5)$. The left boundary is a supersonic inlet, the right boundary is a subsonic outlet and the bottom and top boundaries are walls. The final time for this simulation is $t=0.7s$. 

\subsubsection{Verification of Limiter Accuracy}
\begin{figure}[h]
    \centering
    \includegraphics[height=0.225\textheight]{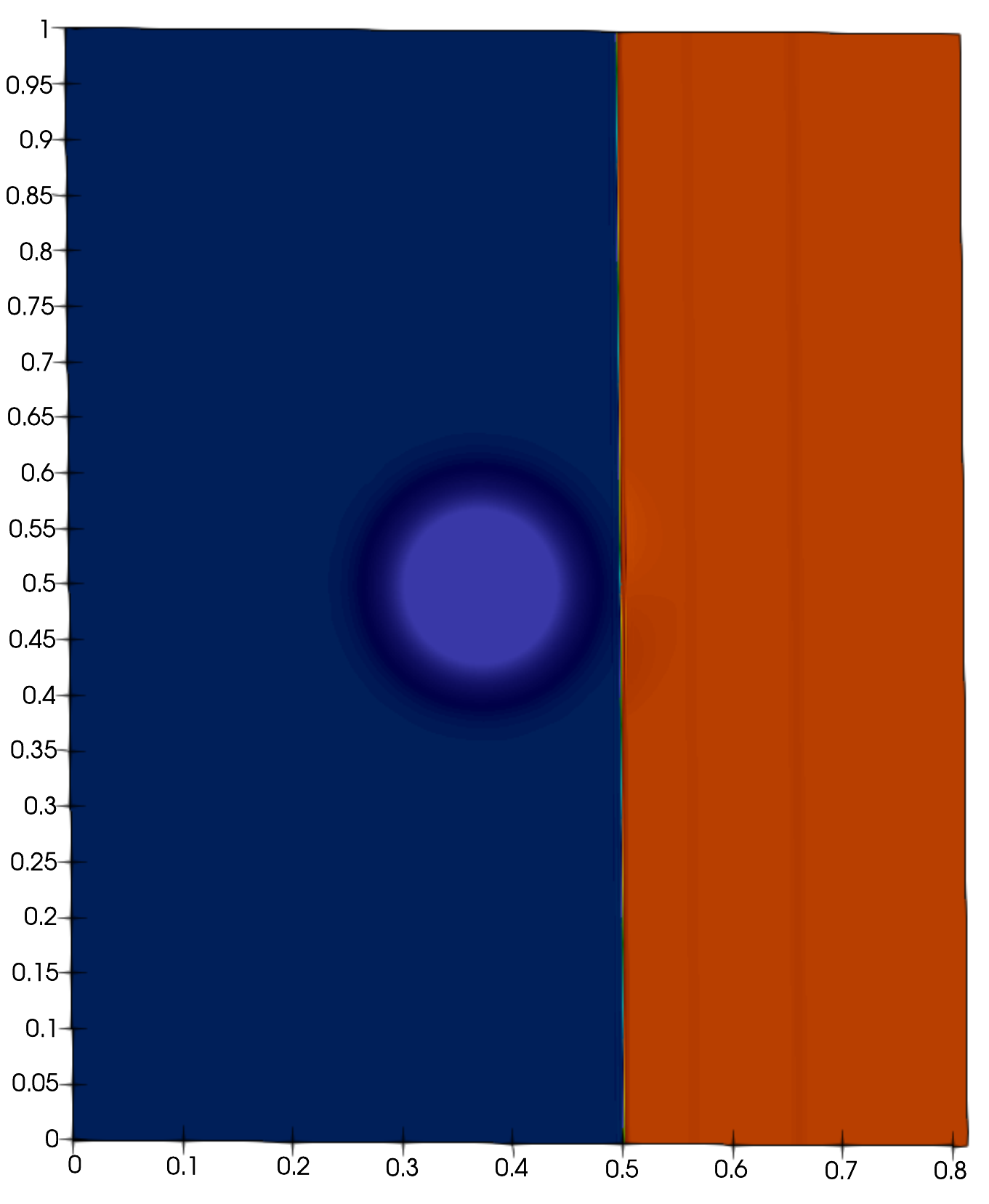}\includegraphics[height=0.225\textheight]{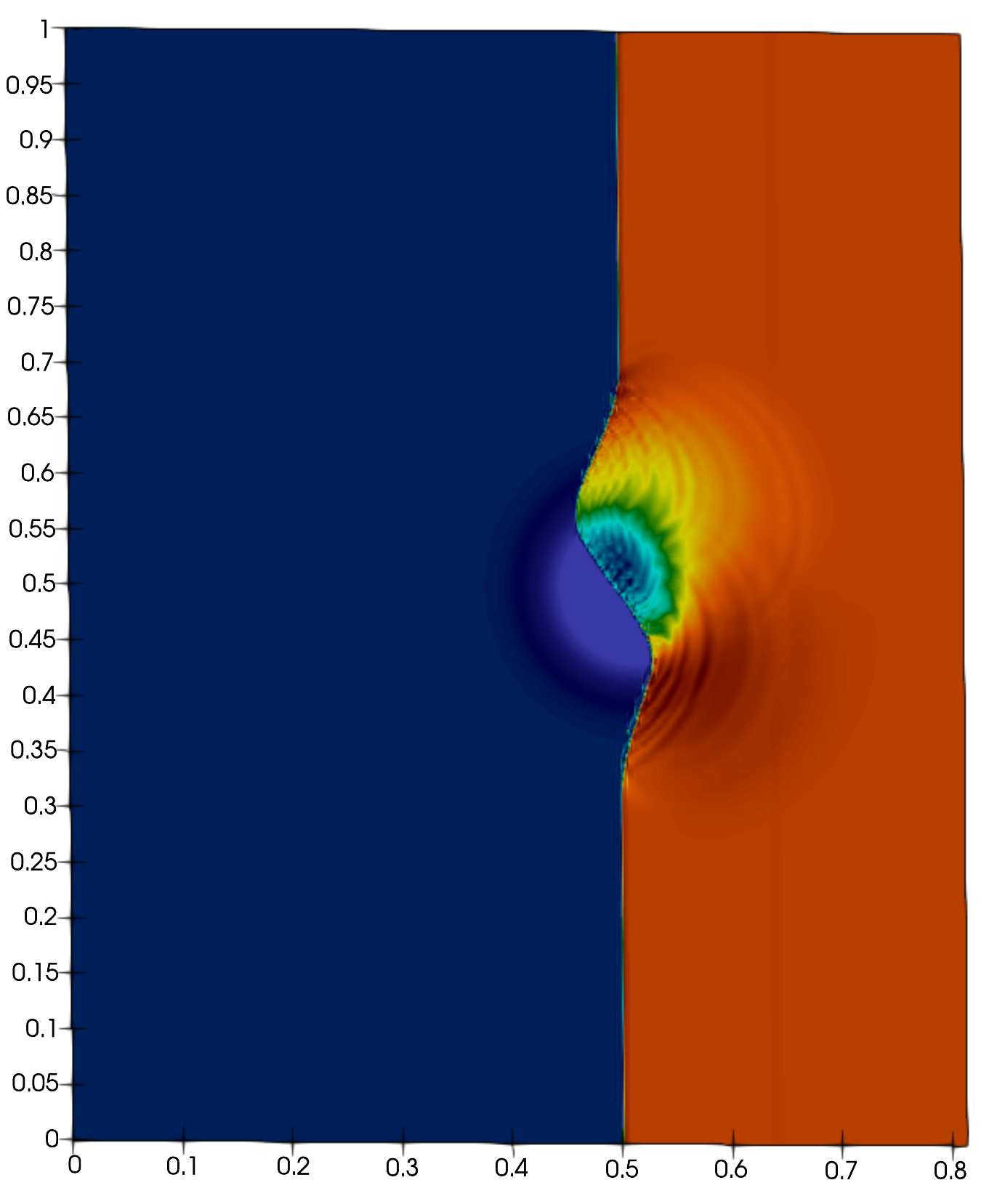}\includegraphics[height=0.225\textheight]{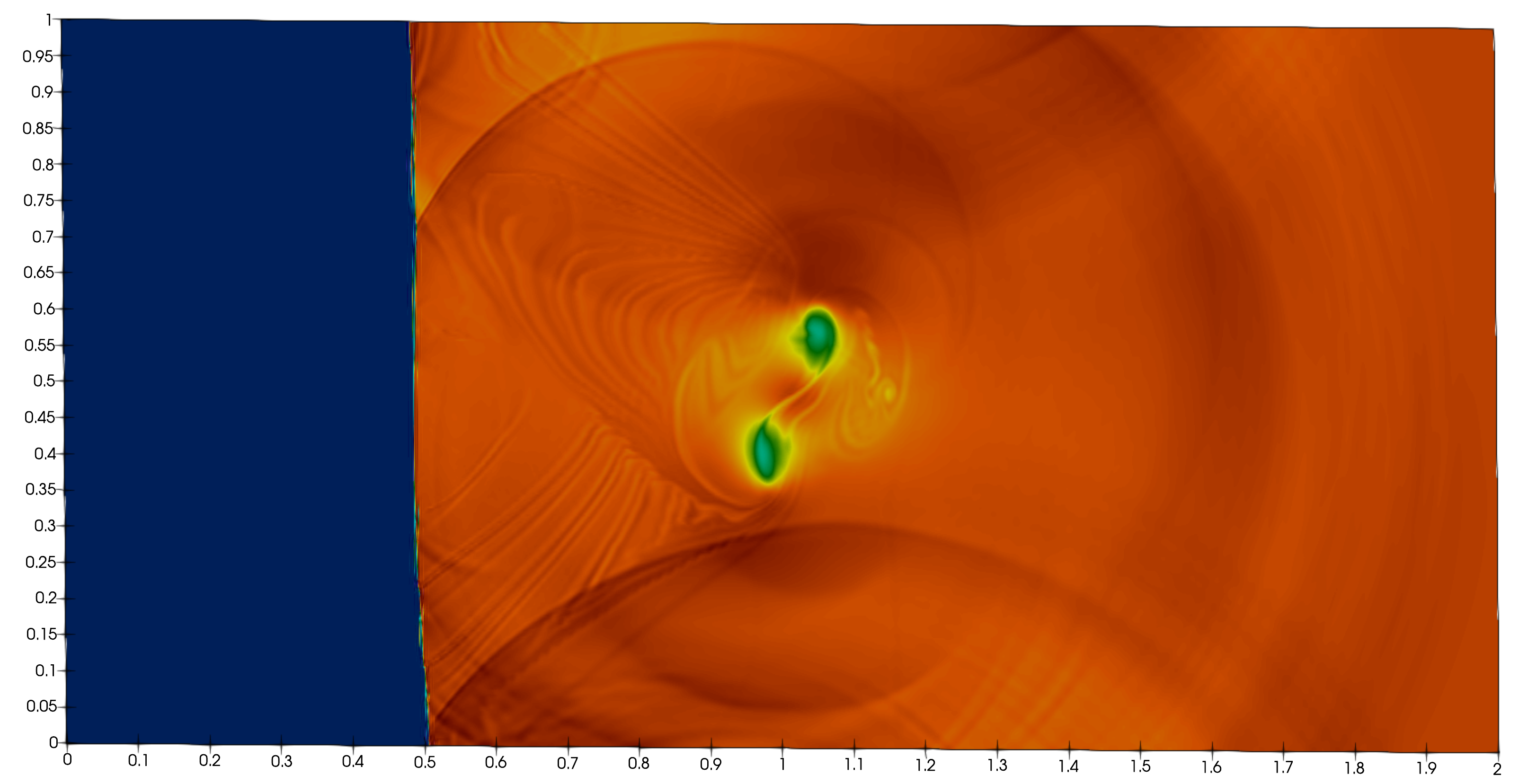}
    \caption{\textit{[SVSW]} Density solution with polynomial degree 3, NSFR-$CH_{RA}$, Roe dissipation, $c_{DG}$ at $t=0.07s$ (Left), $t=0.15s$ (Middle) and $t=0.7s$ (Right)}
    \label{fig:SVSW_noLim}
\end{figure}

The tests are run with the maximum possible $CFL$ before failure. Without the limiter, the test can be run with the $c_{DG}$ scheme for a $CFL$ of $0.27$. With the limiter, the test can be run with a $CFL$ of $0.58$ which is a little over two times greater. The result obtained using the $c_{DG}$ scheme without the PPL is shown in Fig.~\ref{fig:SVSW_noLim}. As expected, the strong vortex is split into two separate vortices due to the compression effects when passing through the shock. The solution contains cylindrical acoustic waves downstream as expected and the shock also cuts off the sound waves centered on the vortex core, resulting in alternating expansion and compression regions.

\begin{figure}
    \centering
    \includegraphics[width=.45\textwidth]{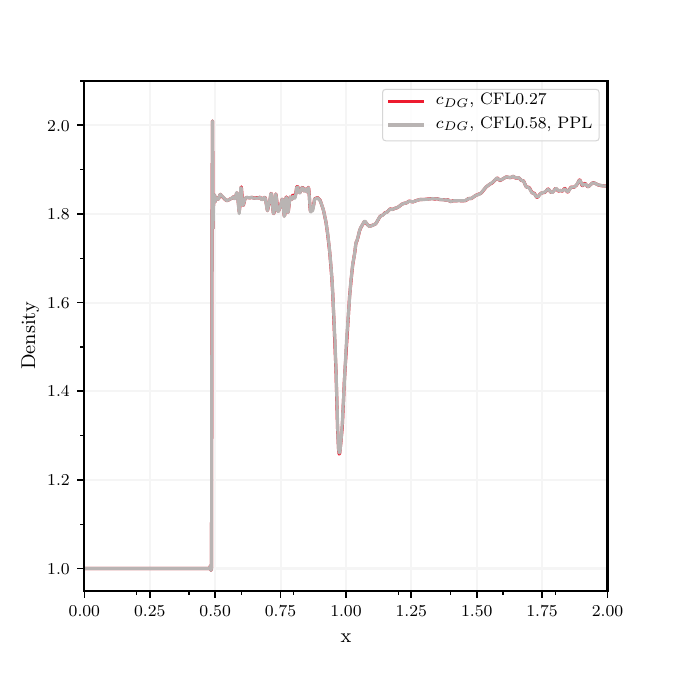}
    \includegraphics[width=.45\textwidth]{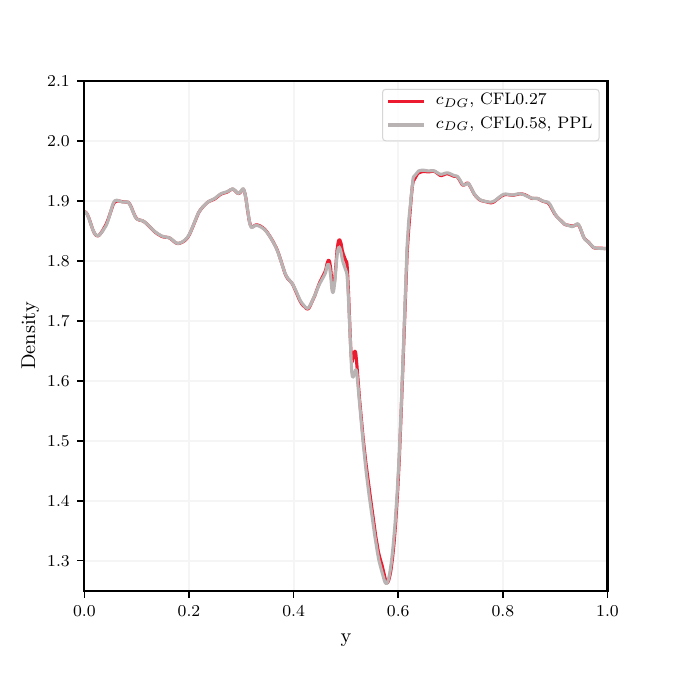}
    \caption{\textit{[SVSW]} Density solution with NSFR-$CH_{RA}$, Roe dissipation, $c_{DG}$, with and without PPL at $t=0.7s$. Plot along $y = 0.4$ (left) and $x = 1.05$ (right).}
    \label{fig:SVSWLines}
\end{figure}

The result obtained with the limiter is visually identical to the results in Fig.~\ref{fig:SVSW_noLim}. The solution contains the expected features, including the two vortex cores, cylindrical acoustic waves, and cutoff sound waves. The solution closely reflects the solution obtained without the limiter, despite being run with a CFL that is a little over two times greater. To better demonstrate that the limiter does not affect the underlying scheme, density line plots are shown. The first density plot is along the line $y = 0.4$ and the second plot is along the line $x = 1.05$. The first plot shows the stationary shock and the first vortex center, while the second plot shows the second vortex center. The plots are shown in Fig.~\ref{fig:SVSWLines}. In both plots, the solution with the limiter follows closely with the unlimited solution and reflects all the features accurately. The only noticeable difference is that the magnitude of the limited solution density is smaller at certain locations. As the limited solution reflects the unlimited solution closely despite having a much larger CFL, it is evident that the limiter does not affect the underlying scheme and its accuracy.

\subsubsection{Comparison with Strong DG}
\begin{figure}[h]
    \centering
    \includegraphics[height=0.25\textheight]{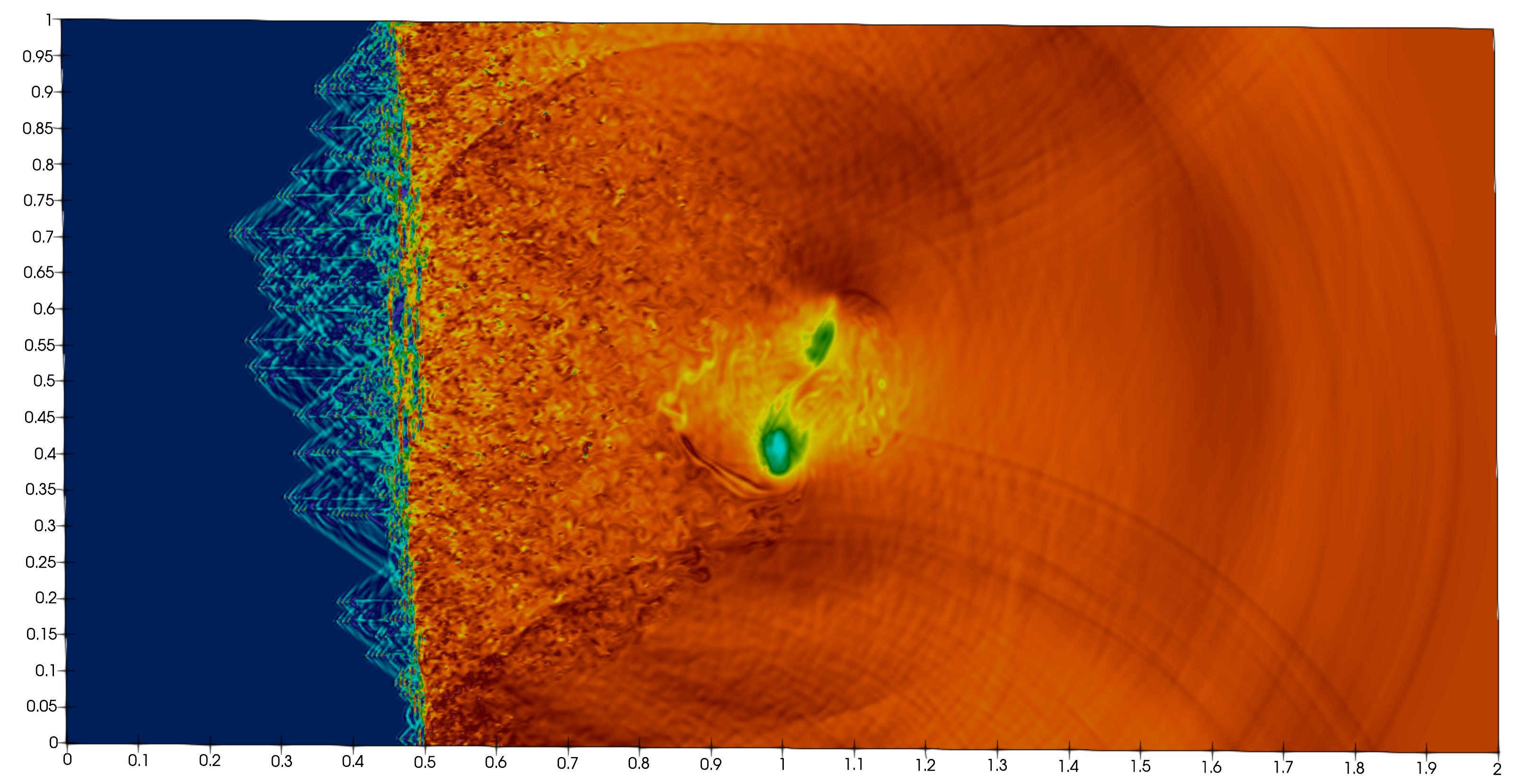}
    \includegraphics[height=0.25\textheight]{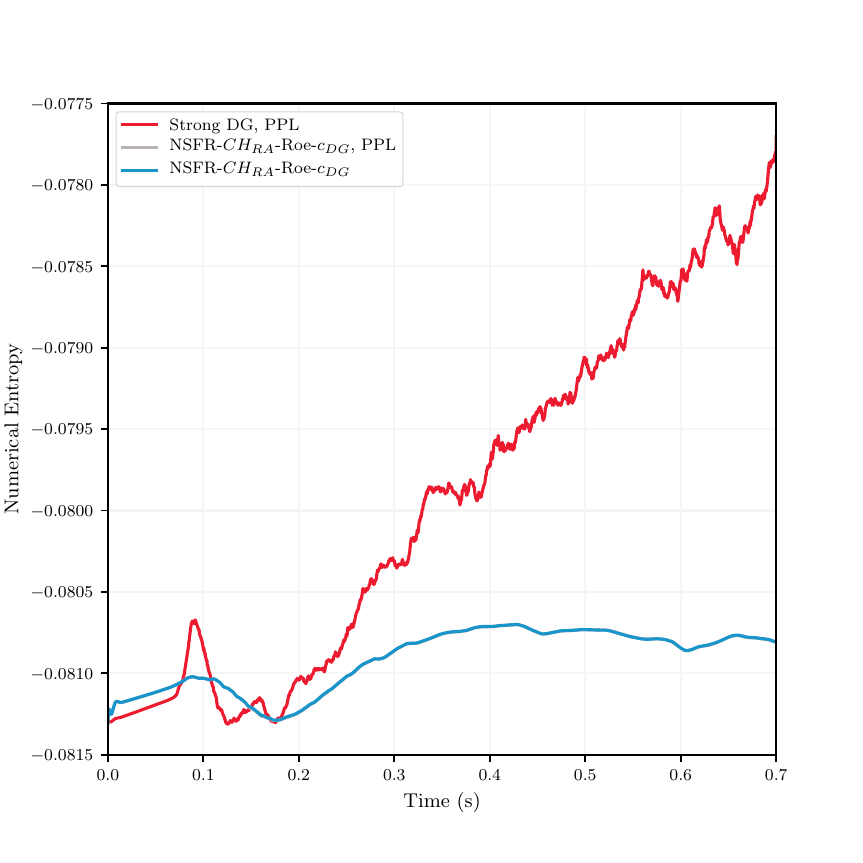}
    \caption{\textit{[SVSW]} Top: Density solution with the SSPRK3 Strong DG method at $t = 0.7s$. Bottom: Plot of entropy over time for the NSFR-$CH_{RA}$, Roe dissipation, $c_{DG}$ method with and without the PPL and SSPRK3 Strong DG method with the PPL}
    \label{fig:svsw_strongdgentropy}
\end{figure}

To highlight the difference between the SSPRK3 Strong DG method and the NSFR method,  the Strong DG results obtained with the Lax-Friedrichs Flux and PPL with a CFL of 0.2 are shown in Fig.~\ref{fig:svsw_strongdgentropy} along with the plot of the entropy over time. The entropy over time for the $c_{DG}$ results with and without the limiter are nearly identical, whereas the SSPRK3 Strong DG solution has increasing entropy over time. This lack of entropy stability is reflected in the Strong DG results, which contain severe oscillations.

\subsubsection{CFL Condition for Two-Point Flux}\label{sec:tpf_svsw}
This test is also used to numerically verify that the two-point flux abides by the $CFL$ conditions, Eq.\ref{eq: cfl-condition-1} and \ref{eq: cfl-condition-2}, in two dimensions. The number of cells for this test is $N_x\times N_y =202\times101$. The test is run for polynomial degrees $p = 2,3,4$ with an increment of $0.00005s$ for the constant time step. Once again, the $CFL$ for a given configuration is calculated in two different ways. The first method uses the minimum $\Delta x$ between the GLL quadrature nodes at the given degree for $\Delta x$. The second method uses the cell size for $\Delta x$. The $CFL$ is calculated using Eq.\ref{eq: cfl-formula}. The results of the numerical study are given in Table \ref{tab:tpf_test_svsw}.

\begin{table}[]
\centering
\caption{Two-point flux positivity preservation test using the Strong Vortex Shock Wave setup with $N_x\times N_y = 202\times101$.}
\label{tab:tpf_test_svsw}
\resizebox{0.7\textwidth}{!}{
\begin{tabular}{rrrrrrr}
\hline
\multicolumn{7}{r}{\textbf{NSFR - $CH_{RA}$ - Roe Dissipation - $\mathbf{c_{DG}}$ - $N_x \times N_y = 202\times101$}} \\ \hline
k & Min.    & Cell Width & Time Step & Max. $\lambda$ & Min. dx CFL                  & Cell Width CFL               \\ \hline
2 & 0.00495 & 0.00990099 & 0.00095   & 3.90           & 0.75                         & 0.37                         \\
  &         &            & 0.001     & 3.90           & 0.79                         & 0.39                         \\
  &         &            & 0.00105   & 3.90           & 0.83                         & 0.41                         \\
  &         &            & 0.0011    & 4.23           & \cellcolor[HTML]{FFA3A3}0.94 & \cellcolor[HTML]{FFA3A3}0.47 \\ \hline
3 & 0.00274 & 0.00990099 & 0.00045   & 3.90           & 0.64                         & 0.18                         \\
  &         &            & 0.0005    & 3.90           & 0.71                         & 0.20                         \\
  &         &            & 0.00055   & 3.90           & 0.78                         & 0.22                         \\
  &         &            & 0.0006    & 9.37           & \cellcolor[HTML]{FFA3A3}2.05 & \cellcolor[HTML]{FFA3A3}0.57 \\ \hline
4 & 0.00171 & 0.00990099 & 0.00025   & 3.90           & 0.57                         & 0.10                         \\
  &         &            & 0.0003    & 3.90           & 0.68                         & 0.12                         \\
  &         &            & 0.00035   & 3.90           & 0.80                         & 0.14                         \\
  &         &            & 0.0004    & 6.75           & \cellcolor[HTML]{FFA3A3}1.58 & \cellcolor[HTML]{FFA3A3}0.27 \\ \hline
\end{tabular}
}
\end{table}

From the previous results, it is evident that the 1D condition can be relaxed up to a $CFL$ of 1. For higher dimensions, this condition still applies as demonstrated by this case wherein the test only fails if the $CFL$ is close to or greater than $1$. For the second method, the $p=2$ configuration fails when the $CFL$ is in the range of $0.5$ as given by Eq.\ref{eq: cfl-condition-2} and the other polynomial degrees fail much earlier as observed in the 1D cases.
This test case numerically validates the condition that positivity is preserved in $2D$ if the $CFL$ is $\leq 1$ as given in Eq. \ref{eq: cfl-condition-1}, irrespective of polynomial degree.

\subsubsection{Impact of Flux Nodes}
\begin{figure}[h]
    \centering
    \includegraphics[width=0.9\textwidth]{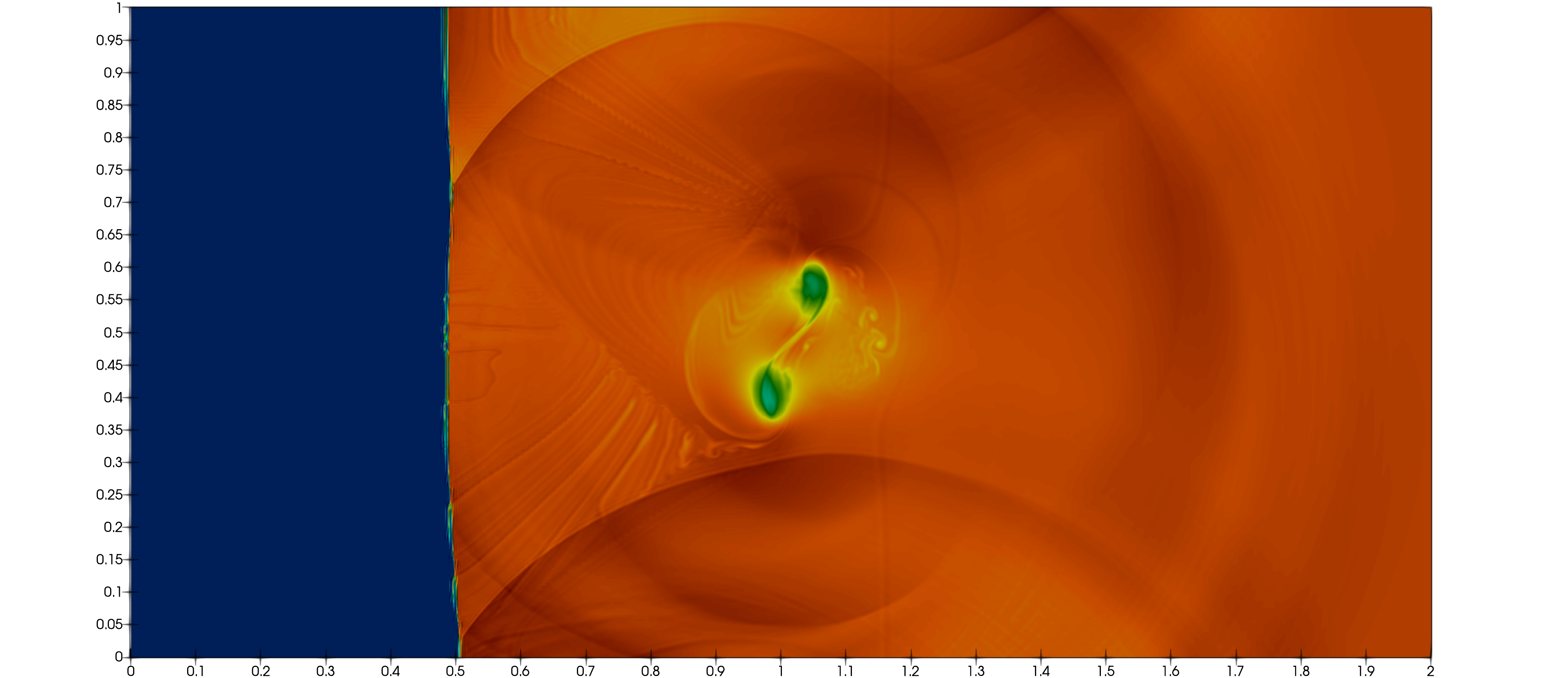}\\
    \caption{\textit{[SVSW]} Density solution with NSFR-$CH_{RA}$, Roe dissipation, $c_{DG}$, PPL using GL flux nodes at $t=0.7s$}
    \label{fig:svsw_GL}
\end{figure}

The application of the PPL for the NSFR scheme with GL flux nodes is also investigated using this test case. For this test, the CFL is lowered to $0.07$ due to the sensitivity of the entropy-projected conservative variables at surface points as discussed for the Shu Osher case in Section \ref{sec: 1DShuOsherProblem}. Additionally, the grid size is changed to $302\times151$ to ensure that the stationary shock is not at the boundaries of the cell. The result obtained with the GL flux nodes is shown in Fig.~\ref{fig:svsw_GL}, and it is visually identical to the results in Fig.~\ref{fig:SVSW_noLim}. The solution contains the expected features, including the two vortex cores, cylindrical acoustic waves, and cutoff sound waves. This case is a good verification that the limiter works as intended when GL flux nodes are used. 

\begin{figure}
    \centering
    \includegraphics[width=.45\textwidth]{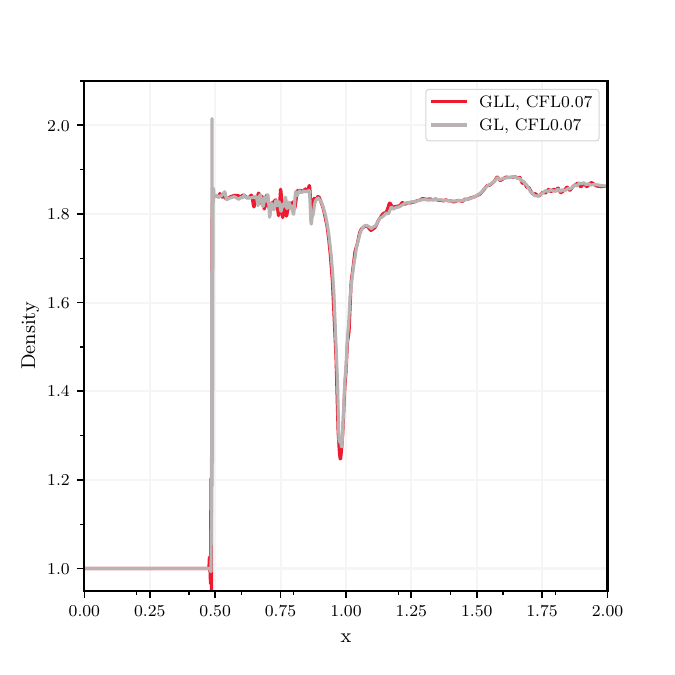}
    \includegraphics[width=.45\textwidth]{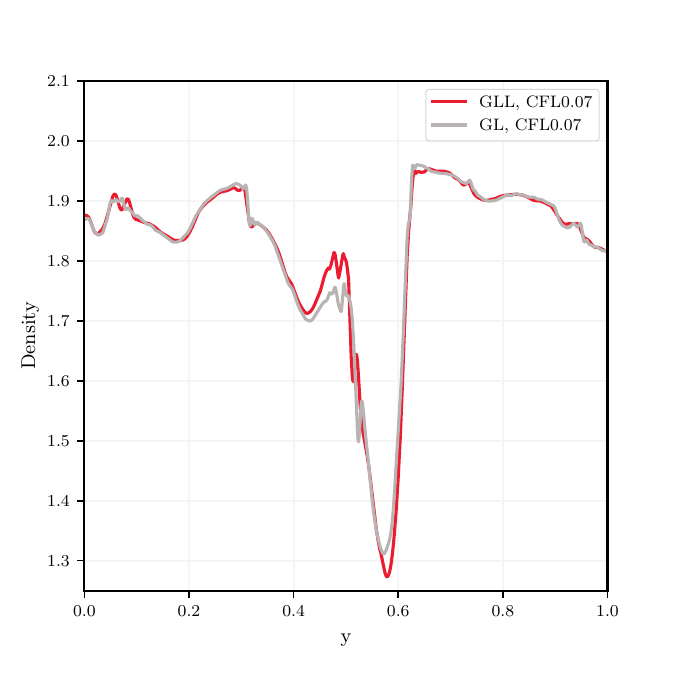}
    \caption{\textit{[SVSW]} Density solution with NSFR-$CH_{RA}$, Roe dissipation, $c_{DG}$, with GLL and GL flux nodes at $t=1.8s$. Plot along $y = 0.4$ (left) and $x = 1.05$ (right).}
    \label{fig:SVSWLines_GL}
\end{figure}

The density plots along the lines $y=0.4$ and $x=1.05$ are also provided for the GL results. To best compare the collocated and uncollocated schemes, the GLL test is rerun for a CFL of $0.07$. The plots are shown in Fig.~\ref{fig:SVSWLines_GL}. Both plots accurately reflect the solution features, but differences can be seen in magnitude for the vortex in the plot along $x=1.05$. The GLL solution also exhibits slightly more oscillatory behaviour. The GL solution has fewer oscillations as a result of the greater integration strength of GL nodes. Overall, this case demonstrates that the implementation of the PPL behaves as expected for GL flux nodes.

\subsubsection{Two-Point Flux Comparison}
\begin{figure}[h]
    \centering
    \includegraphics[width=0.35\textwidth]{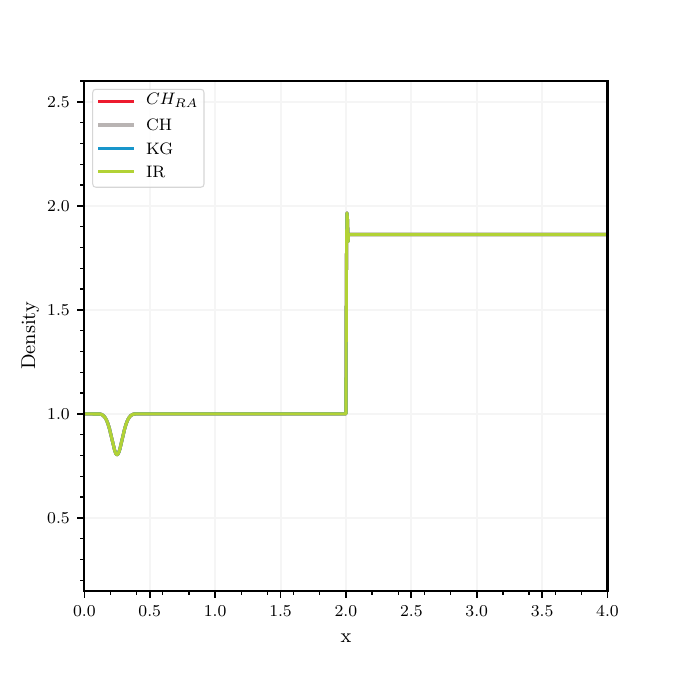}\includegraphics[width=0.35\textwidth]{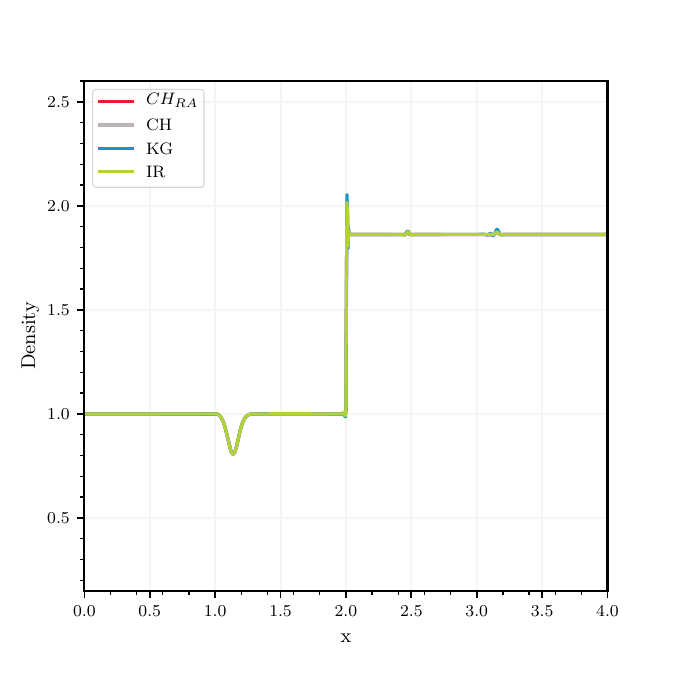}\includegraphics[width=0.35\textwidth]{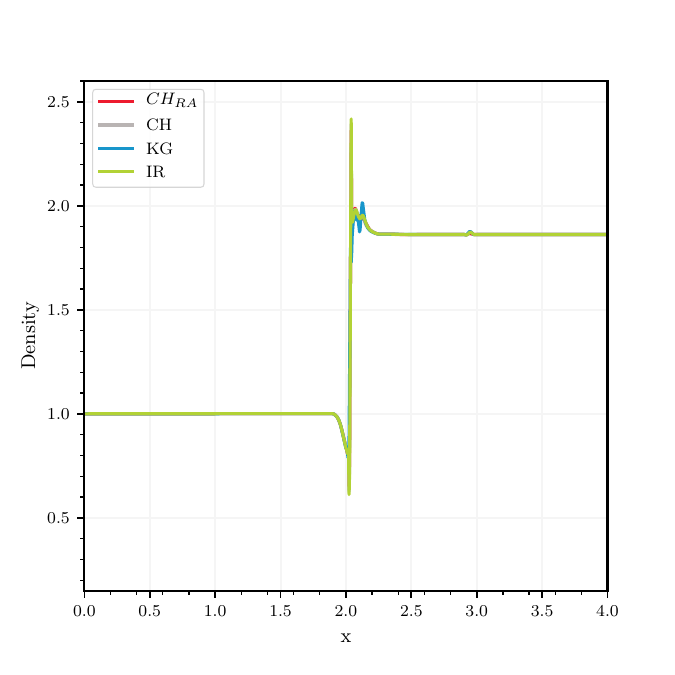}\\
   \includegraphics[width=0.35\textwidth]{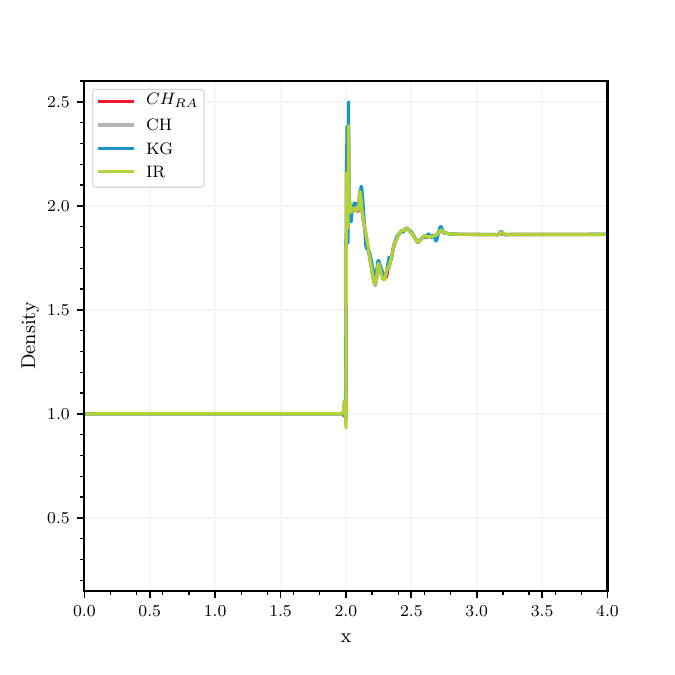}\includegraphics[width=0.35\textwidth]{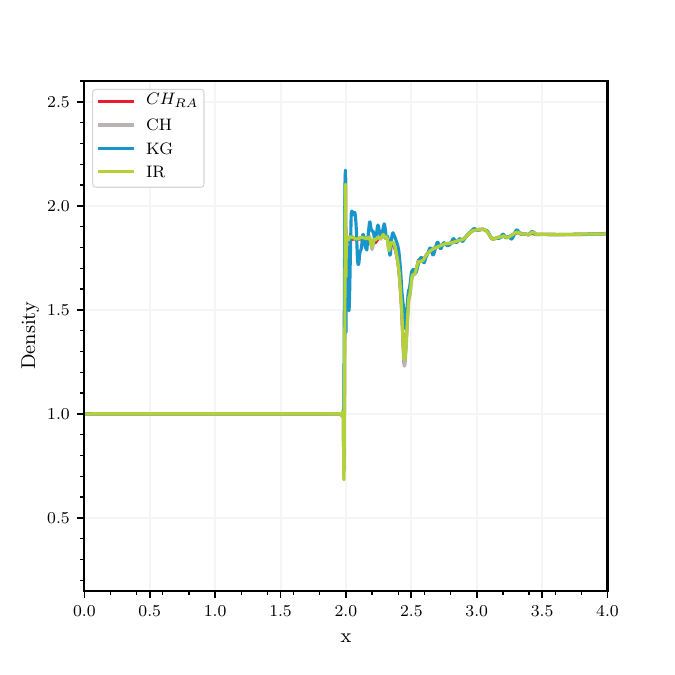}\includegraphics[width=0.35\textwidth]{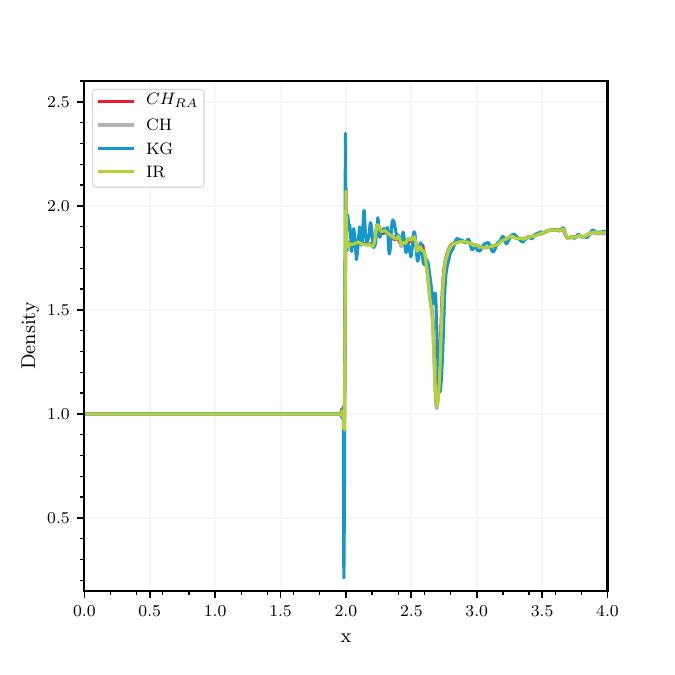}
    \caption{\textit{[Extended SVSW]} Plot of density along $y=0.425$ with NSFR, Roe dissipation, $c_{DG}$ using different entropy conserving fluxes. Top row: $t=0.00s$ (Left), $t=0.50s$ (Middle) and $t=1.00s$ (Right). Bottom row: $t=1.25s$ (Left), $t=1.50s$ (Middle) and $t=1.75s$ (Right)}
    \label{fig:SVSW_ecFlux}
\end{figure}

This test is also used to investigate the behaviour of different two-point fluxes for strong shock problems. To best demonstrate that the behaviour of the fluxes differs in the presence of shocks, the SVSW test is modified slightly. The domain is expanded to $[0,4]\times[0,1]$, the stationary shock is placed at $x=2.0$, and the final time is set to $t=1.75s$. This allows for the vortex to be tracked for a longer period of time before it interacts with the shock and will enable better comparison of the four fluxes. As the domain has been expanded, the grid size is also changed to $200\times50$ to ensure that $\Delta x = \Delta y$.  A plot of density along $y=0.425$ is shown for 6 different time steps in Fig.~\ref{fig:SVSW_ecFlux}. 

\begin{figure}[h]
    \centering
    \includegraphics[width=0.6\textwidth]{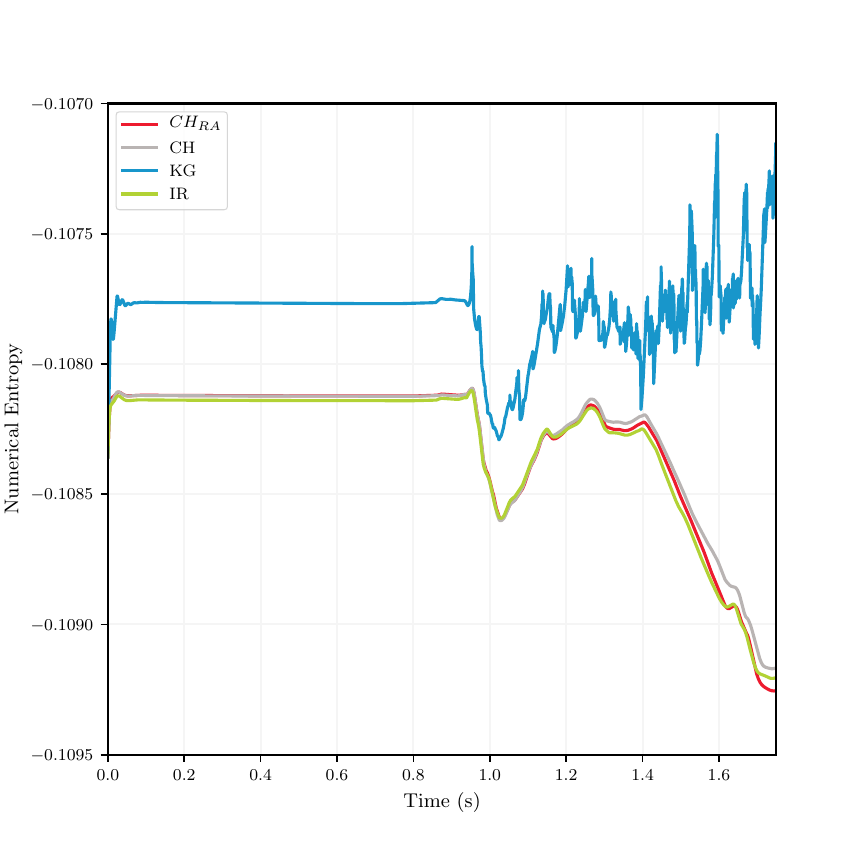}
    \caption{\textit{[Extended SVSW]} Plot of entropy over time with NSFR, Roe dissipation, $c_{DG}$, PPL for the different entropy conserving fluxes.}
    \label{fig:SVSW_ecFlux_entropy}
\end{figure}

The plots of the density before the vortex interacts with the shock show that in the region upstream of the shock, the different fluxes produce identical results. Downstream of the shock, oscillations are observed of slightly varying magnitudes. The plots in the bottom row show the density along $y=0.425$ after the vortex interacts with the shock. The results for the IR, CH and $CH_{RA}$ fluxes are relatively similar. Comparing the KG flux to the three aforementioned fluxes, it is obvious that the results obtained with the KG flux are more oscillatory than the rest. The behaviour of KG is similar to that of the SSPRK3 Strong DG method. To verify that this behaviour is related to entropy stability, a plot of entropy over time for all four fluxes is presented in Fig.~\ref{fig:SVSW_ecFlux_entropy}. The plot for the entropy over time demonstrates that the added dissipation of the KG flux is insufficient to ensure entropy stability, as expected. This leads to increasing entropy over time, which presents itself as severe oscillations in the solution similar to those seen in the Strong DG results.   

\break
\subsection{Shock Diffraction}
\label{sec: 2DShockDiffraction}
The shock diffraction case is a 2D problem of a right-moving shock that passes a backward-facing corner. The Mach 5.09 shock diffracts at a $90\text{\textdegree}$ edge. The solution is expected to have several key features. These features include the incident shock, the diffracted shock, the reflected expansion wave and the vortex roll-up as highlighted in \cite{hillier1991computation}. This case is a good test of the positivity-preserving strategy as it is easy to get nonphysical values below the $90\text{\textdegree}$ corner and to the right of the corner.
The final time used is $T=2.3$. The computational domain is $[0,1]\times[6,11]\cup[1,13]\times[0,11]$. The initial condition is
\begin{equation}
    (\rho, u, v, p) = \begin{cases}
        (\rho_{\text{left}}, u_{\text{left}}, v_{\text{left}}, p_{\text{left}}) \qquad &\text{if} \: x < 0.5\\
        (1.4,0,0,1) \qquad &\text{if} \: x > 0.5\\
    \end{cases}\\
\end{equation}
Where $\rho_{\text{left}} = 7.041132906907898$, $u_{\text{left}} = 4.07794695481336$, $v_{\text{left}} = 0$, $p_{\text{left}} = 30.05945$. The conditions for the left boundary is inflow for ${x=0,6\leq y\leq 11}$ and wall boundary for ${x=1,0\leq y\leq6}$. The bottom boundary is a wall for ${0\leq x\leq1,y=6}$ and outflow for (${y=0,1\leq x\leq13}$). The right boundary, ${x=13,0\leq y\leq11}$, is outflow and the top boundary, ${y=11,0\leq x\leq13}$, is a wall boundary. All the schemes are run with a $CFL$ of $0.45$.

\begin{figure}
    \centering
    \includegraphics[width=0.5\linewidth]{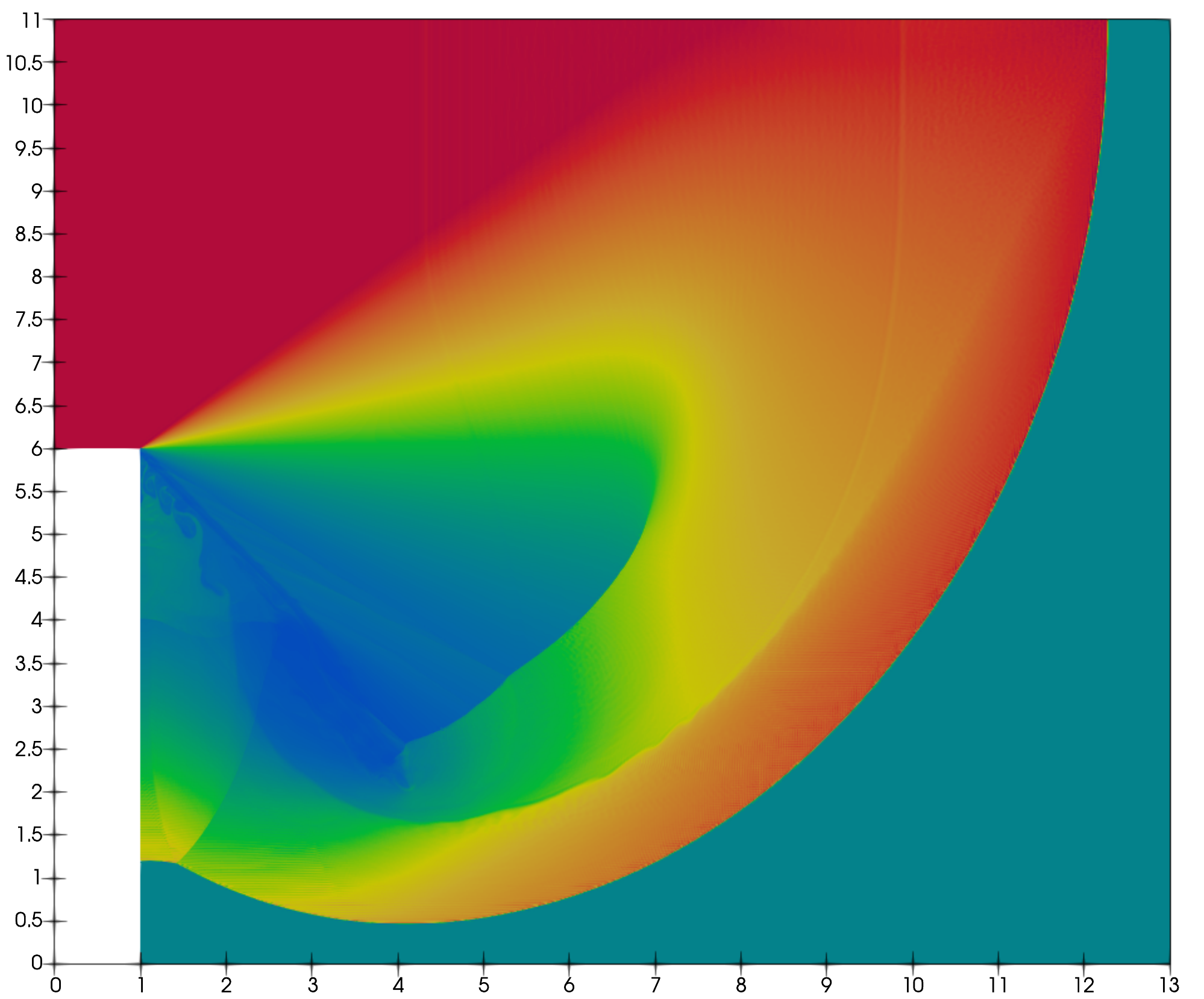}
    \caption{\textit{[Shock Diffraction]} Density solution with NSFR-$CH_{RA}$, $c_{DG}$, Roe dissipation, $p=3$, PPL at $t=2.3s$}
    \label{fig:shockDiffraction_cdg}
\end{figure}

The test fails due to nonphysical values at the given CFL when using the PPL with only the modifications implemented by \citet{WANG2012653}. However, with the modifications introduced in this paper, we are able to obtain results that are comparable to the solution in \citet{ZHANG20108918}, which uses a TVB limiter in addition to the original PPL. The results with the modified limiter are shown in Fig.~\ref{fig:shockDiffraction_cdg}.

\subsubsection{Investigation of Flux Reconstruction Schemes}
\begin{figure}
    \centering
    \includegraphics[width = 0.47\textwidth]{shockDiff_cDG.pdf}
    \includegraphics[width = 0.47\textwidth]{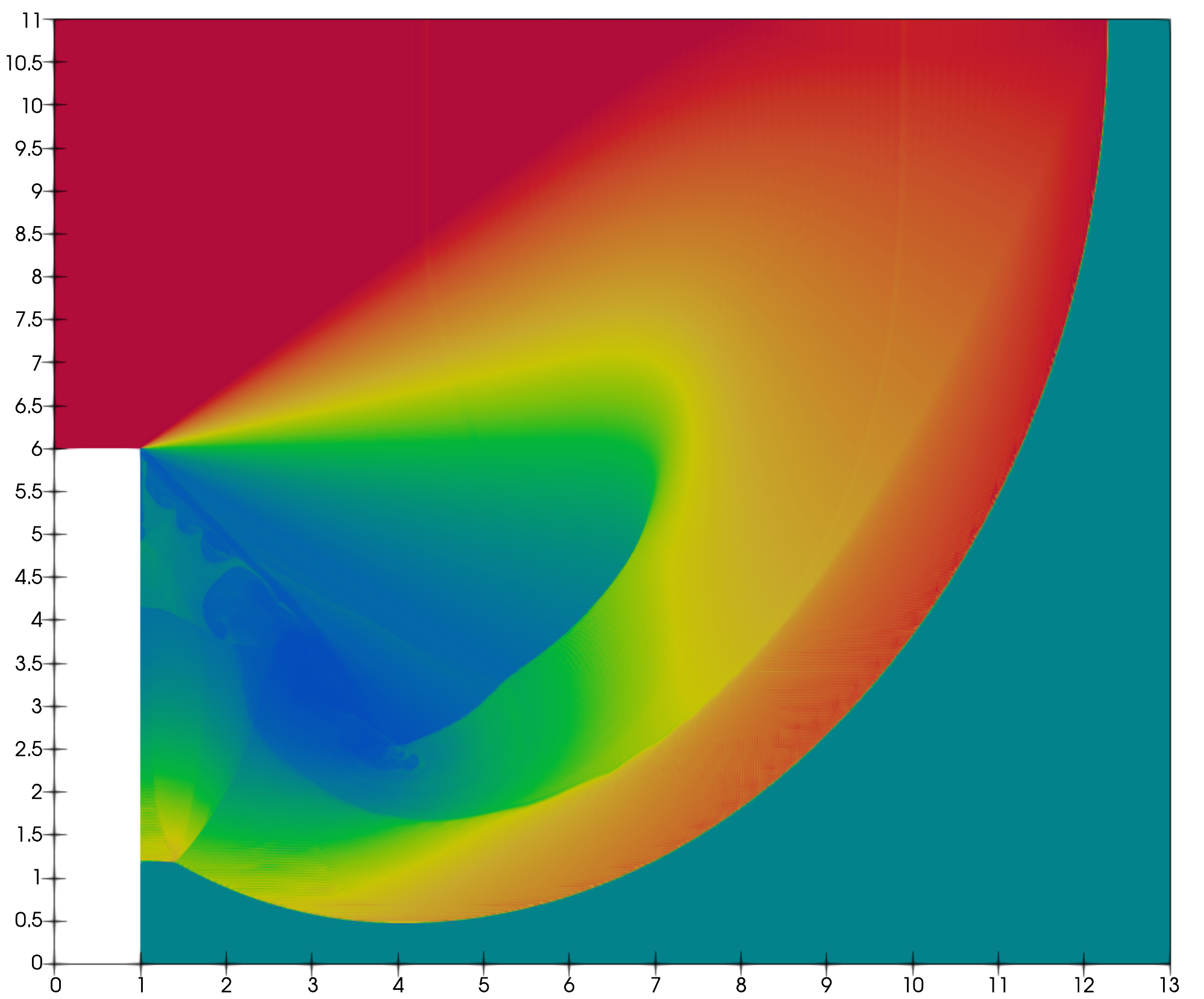}\\
    \includegraphics[width = 0.47\textwidth]{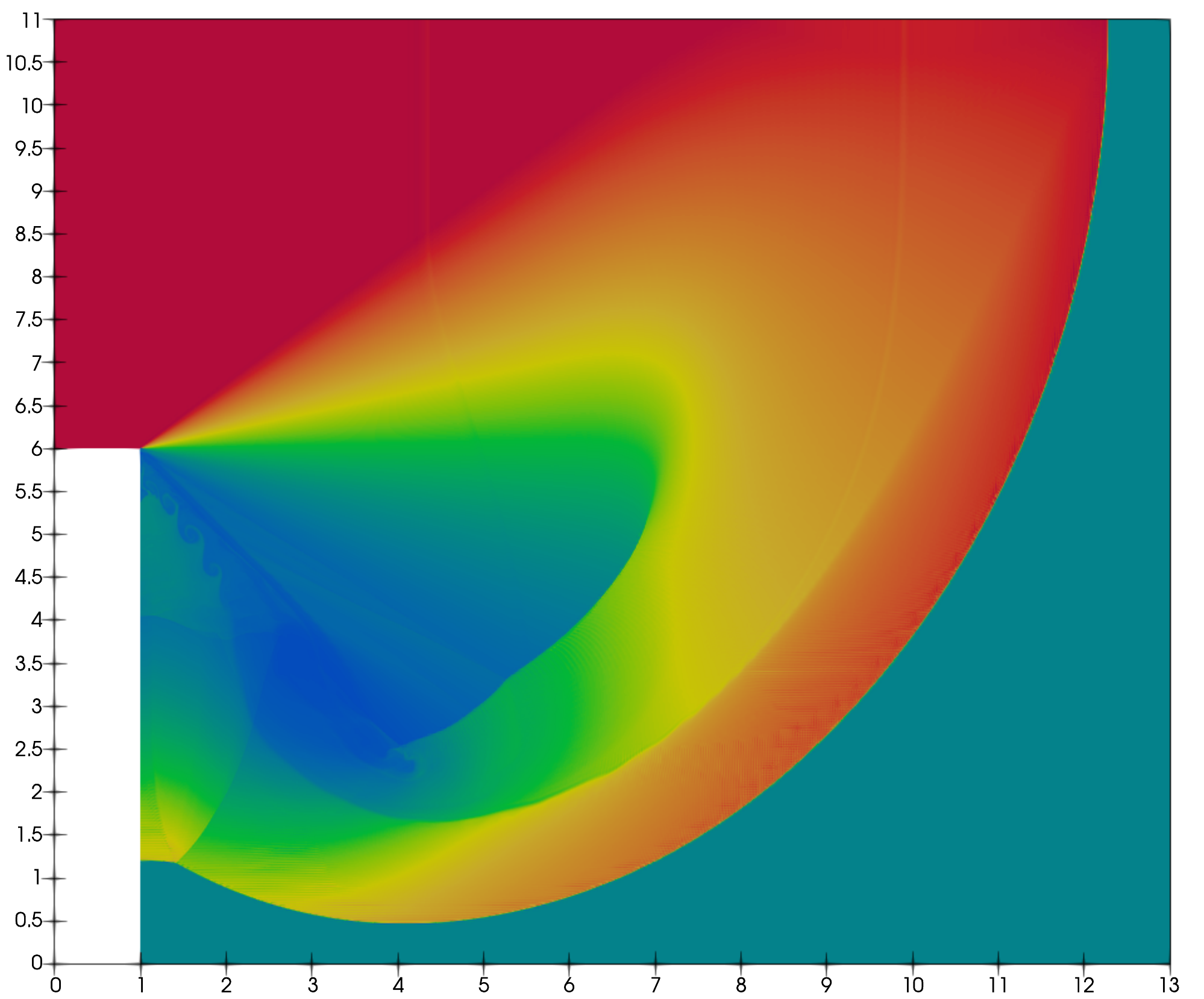}
    \includegraphics[width = 0.47\textwidth]{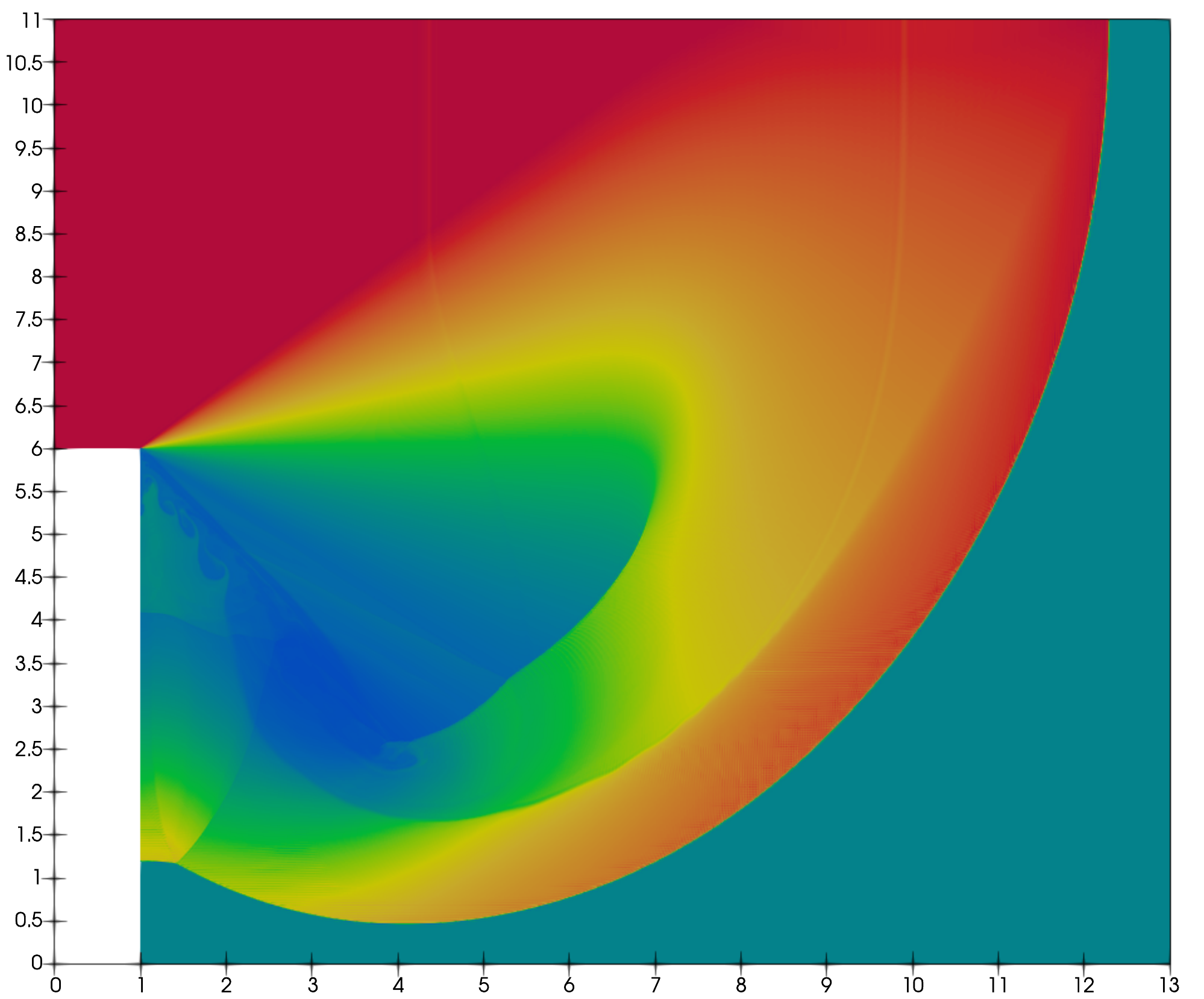}
    \caption{\textit{[Shock Diffraction]} Density solution with NSFR-$CH_{RA}$, Roe dissipation, $p=3$, PPL at $t=2.3s$, using different correction parameters - $c_{DG}$ (top left), $c_{SD}$ (top right), $c_{HU}$ (bottom left), $c_{+}$ (bottom right).}
    \label{fig:ShockDiffraction_cParamFinalDensity}
\end{figure}

This test is also used to compare the impact of various flux reconstruction schemes. Four different schemes are run including $c_{DG}$, $c_{SD}$, $c_{HU}$ , and $c_{+}$. It is run with a grid of $416\times352$.  The final result for each configuration is shown in Fig.~\ref{fig:ShockDiffraction_cParamFinalDensity}. All four results reflect the expected features, including the vortex roll-up and Kelvin-Helmholtz instabilities (KHI). It should be noted that the $c_{DG}$ solution has more noise across the solution and a greater overshoot at the shock compared to the other three schemes. 

\begin{figure}
    \centering
    \includegraphics[width=0.47\textwidth]{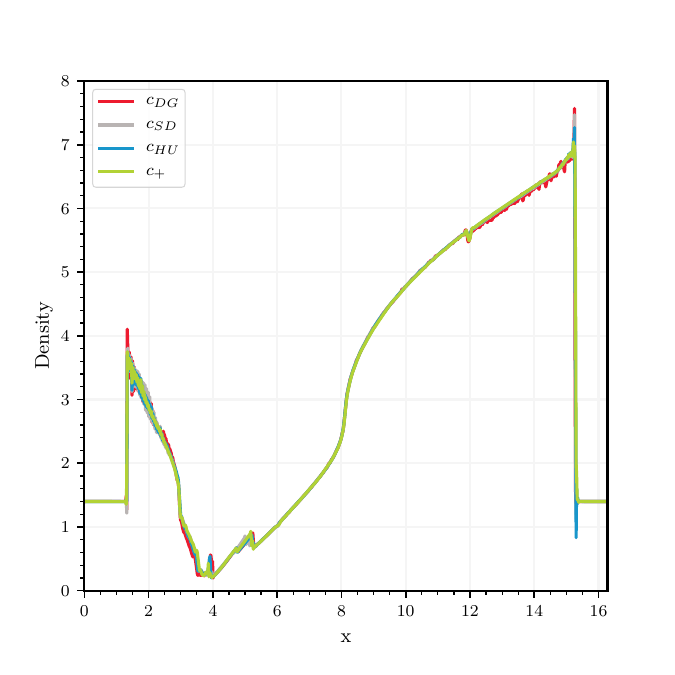}
    \includegraphics[width=0.47\textwidth]{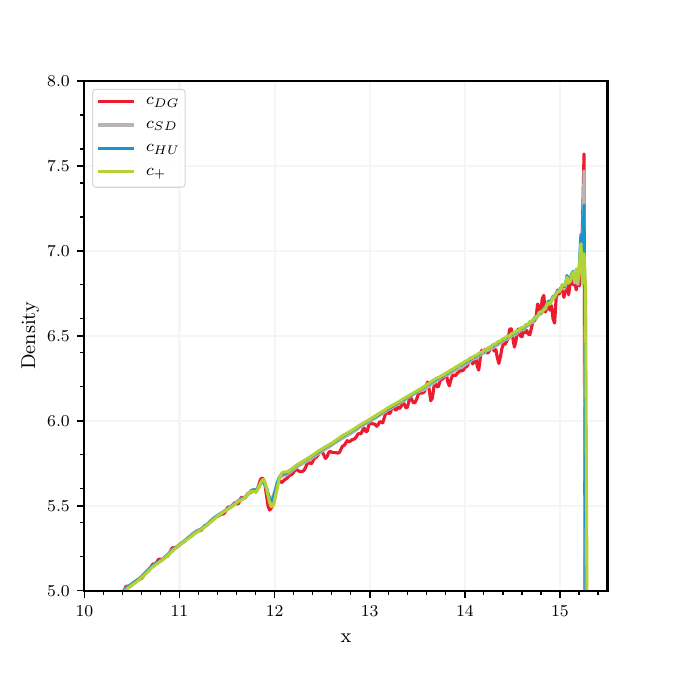}
    \caption{\textit{[Shock Diffraction]} Plot of density across $(1,0)\rightarrow(13,11)$ with NSFR, Roe dissipation, PPL at $t=2.3s$ using different flux reconstruction correction parameters (left) with a closer look at the wake of the shock (right)}
    \label{fig:ShockDiffraction_diffFluxes_plotAcross}
\end{figure}

To demonstrate the noise present in the $c_{DG}$ solution, a plot of density along the line $(1,0)\rightarrow(13,11)$ is provided for all the schemes in Fig.~\ref{fig:ShockDiffraction_diffFluxes_plotAcross}. The plot is also zoomed in to the wake of the shock to better show the difference in the solutions. In the areas of the plot that are upstream of the shock, comparing the schemes is difficult due to the presence of the KHI. In the plot focusing on the wake of the shock, we can see that there is a significant amount of oscillations in the $c_{DG}$ solution compared to the other solutions. In addition to the oscillations, the $c_{DG}$ solution also has the greatest overshoot at the shock. The overshoot observed at the shock is mitigated as the correction parameter increases in value.

\begin{figure}
    \centering
    \includegraphics[width=0.47\textwidth]{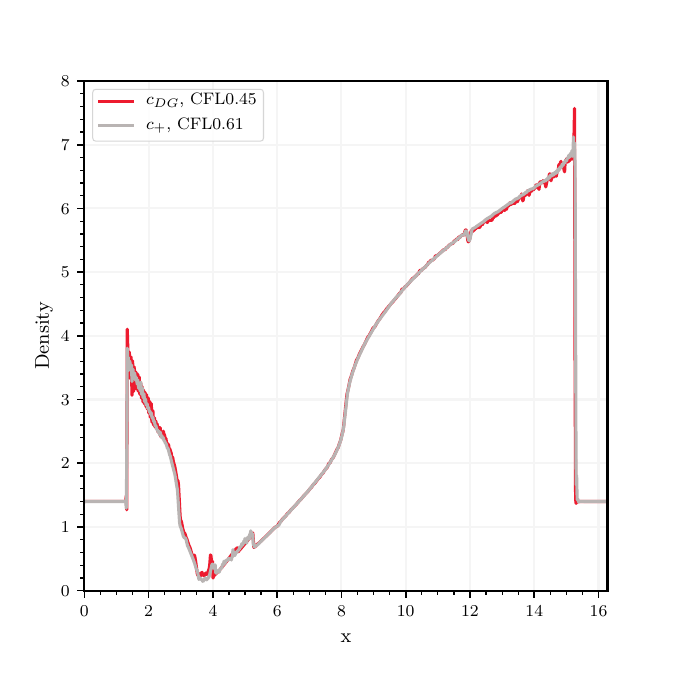}
    \includegraphics[width=0.47\textwidth]{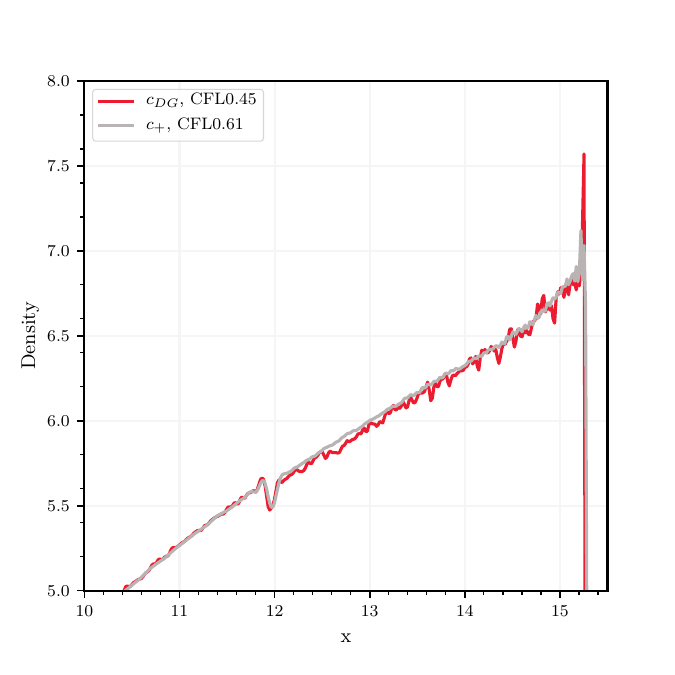}
    \caption{\textit{[Shock Diffraction]} Plot of density across $(1,0)\rightarrow(13,11)$ with NSFR, Roe dissipation, PPL at $t=2.3s$ (left) with a closer look at the wake of the shock (right)}
    \label{fig:ShockDiffraction_c+_maxCFL}
\end{figure}

The test configurations for this case so far have demonstrated the effects of increasing the $c$ value with the $CFL$ kept consistent across all runs. Next, we change the $CFL$ for the $c_+$ configuration to show that the oscillations are mitigated even for a larger timestep. This test is run for the maximum CFL before failure for $c_+$, which is $0.61$. The density plot along the line $(1,0)\rightarrow(13,11)$ is shown in Fig.~\ref{fig:ShockDiffraction_c+_maxCFL}. In addition to the higher $CFL$ case, the $c_{DG}$ results from the previous run are also shown for comparison. The plots show that the increase in $c$ value continues to mitigate oscillations and overshoot, providing a better solution than $c_{DG}$, despite the significant increase in the time step size. The added robustness and the mitigation of oscillations by the $c_+$ scheme make it a very attractive implementation for strong shock problems that involve complex features.

\break
\subsection{Double Mach Reflection}
The Double Mach Reflection (DMR) problem introduced by Woodward and Colella involves sending a Mach 10 shock into a reflecting wall. The original physical problem has the shock moving horizontally into a wall inclined at a 30\textdegree angle. However, the problem is generally solved with a shock moving diagonally into a horizontal wall on the rectangular domain. The rectangular domain is $[0,4]\times[0,1]$ with the initial condition,
\begin{equation}
    (\rho, u, v, p) = \begin{cases}
        (8,\frac{33\sqrt{3}}{8},-\frac{33}{8},116.5) \qquad &\text{if} \: y > \sqrt{3}(x-\frac{1}{6})\\
        (1.4,0,0,1) \qquad &\text{if} \: y < \sqrt{3}(x-\frac{1}{6})
    \end{cases}
\end{equation}
The left boundary is inflow, and the right boundary is outflow. The bottom boundary corresponding to $0\leq x<\frac{1}{6}$ has the post-shock condition. The bottom boundary corresponding to $\frac{1}{6}<x\leq4$ is a wall boundary. In most implementations of the DMR case, the top boundary is set to a function that follows the motion of the shock. To minimize implementation effort, the domain is extended in the $y$-direction such that the top boundary doesn't affect the domain of the test case. Thus, the computational domain is $[0,4]\times[0,3]$. The final time is $t=0.2s$. As outlined in \citet{vevek2019alternative}, the solution at the final time is a self-similar shock structure with two triple points and two slip lines. The primary triple point is where the incident shock, reflected shock, and Mach stem meet, while the secondary triple point is where the reflected shock, secondary Mach stem, and secondary reflected shock meet. The solution also contains two slip lines. The reflected shock interacts with the primary slip line and generates KHI along the slip line. A jet is also present near the wall. The complex features are best captured using a high-order scheme; thus, the case serves as a good indicator of the resolution of a scheme.

Running this test in the configurations given in the following subsections using the limiter in \cite{WANG2012653}, results in the test failing due to nonphysical values. However, with the new modifications, we are able to obtain results that are comparable to solutions seen in literature. The timestep efficiency introduced by this modification is evident in this case and the Shock Diffraction case.

\subsubsection{Polynomial Degree Comparison}
\begin{figure}
    \centering
    \includegraphics[width=0.85\textwidth,trim={0cm 3cm 0 3.5cm},clip]{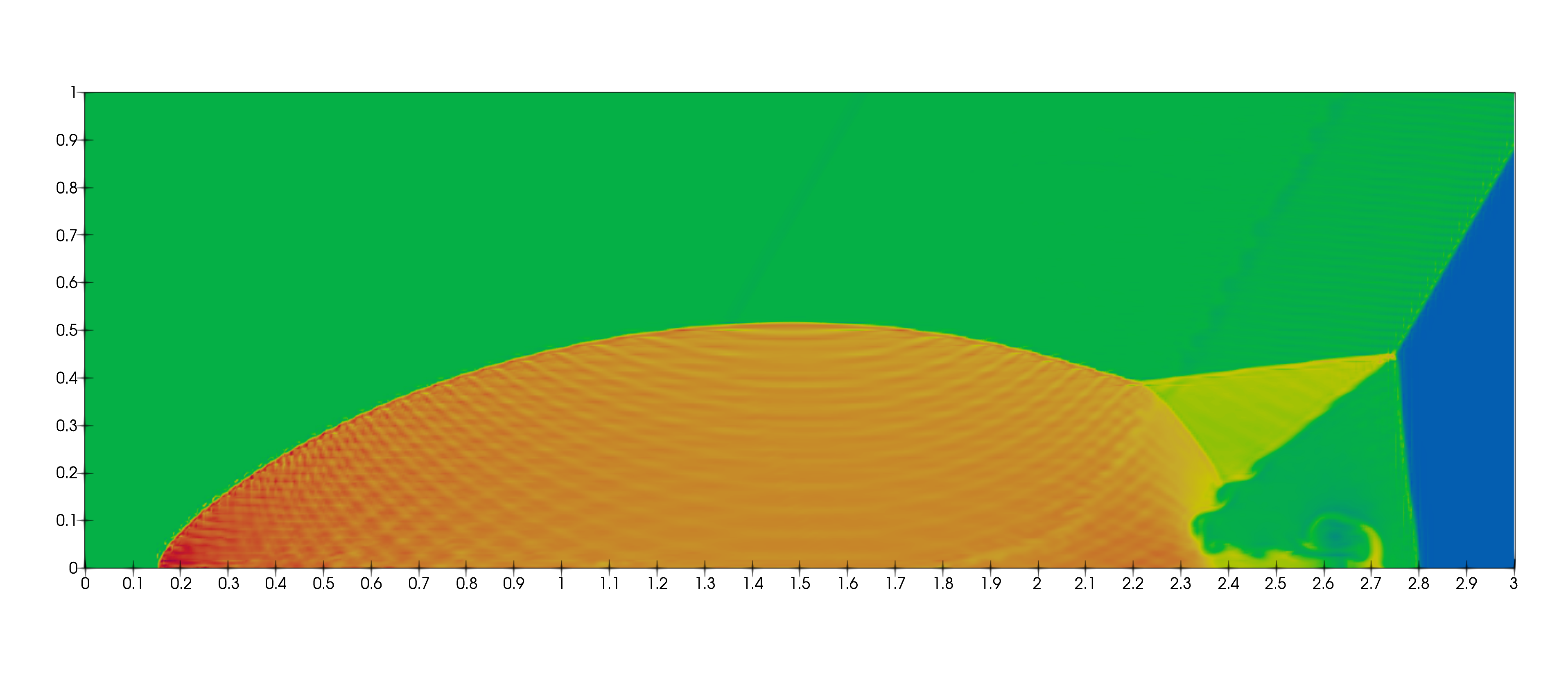}\\
    \includegraphics[width=0.85\textwidth,trim={0cm 3cm 0 3.5cm},clip]{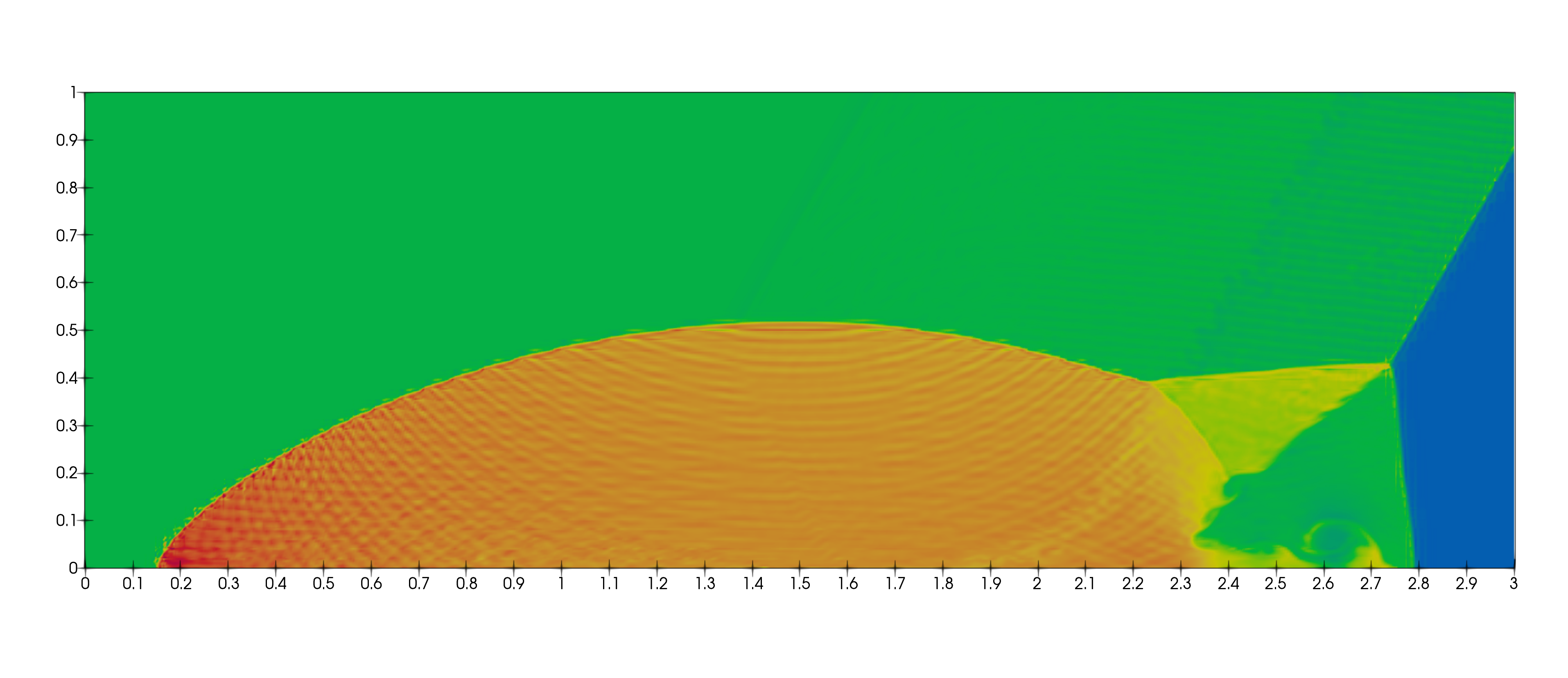}\\
    \includegraphics[width=0.85\textwidth,trim={0cm 3cm 0 3.5cm},clip]{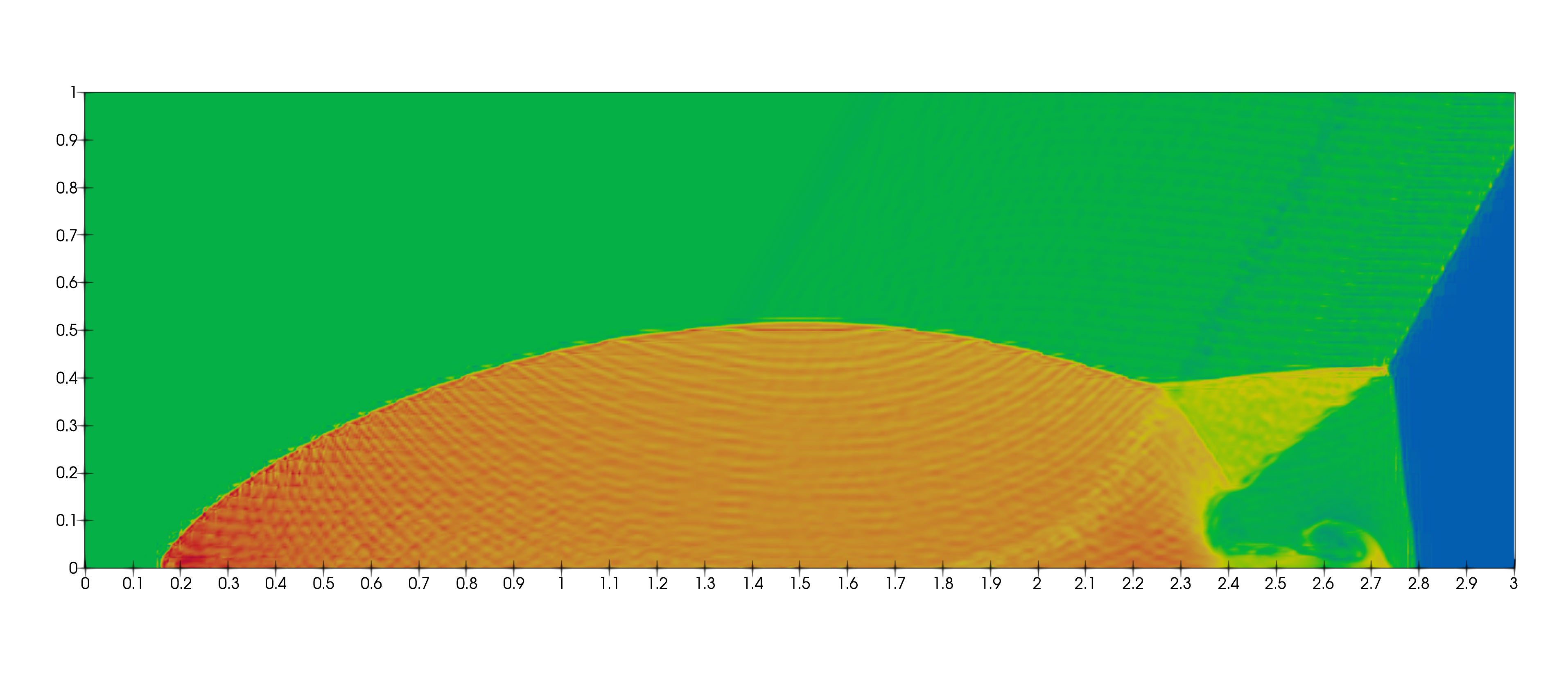}\\
    \caption{\textit{[DMR]} Density solution with NSFR, Roe dissipation, $c_{DG}$, PPL at $t=0.2s$, using different polynomial degrees - $3$ (top), $4$ (middle), $5$ (bottom)}
    \label{fig:DMRPolyDegree}
\end{figure}

The DMR case is first run with increasing polynomial orders while maintaining an equivalent number of degrees of freedom to demonstrate the ability to use higher order polynomials for cases involving strong shocks. The polynomial degrees used in this case are $p = 3, 4, 5$, and a $CFL$ number of $0.15$ is chosen for all degrees. These results are obtained using relatively coarse grids with $240\times60$, $192\times48$, and $160\times40$ for the test space for each polynomial degree respectively. The results are shown in Fig.~\ref{fig:DMRPolyDegree}.

Looking at the $p=3$ results first, we see all the expected solution features are present despite the coarseness of the grid. The solution contains two triple points, one located at $x\approx2.75$ and the other at $x\approx2.2$. There are also two slip lines present. The first slip line is prominent; where, the interaction of the reflected shock and the slip line is evident from the KHI that are generated. The secondary slip line is often difficult to discern in DMR solutions, and the presence of the slip line is a good indicator of the resolution of the scheme. In the $p=3$ solution, the secondary slip line is present and can be identified at $x\approx1.85$ at the bottom boundary or at $x\approx 2.6$ at the top boundary.

As the polynomial degree is increased, we expect that there will be an increase in oscillatory behaviour as demonstrated for the Sod Shock Tube case in Sec.\ref{sec: sod_polydegree}. Despite the presence of the spurious oscillations, it is expected that the solutions still contain all the expected features. In both the $p=4$ and $p=5$ solutions, it is evident that all the expected solution features are present and follow closely with the features described in the $p=3$ results. The secondary slip line is also clearly present in the solution despite the increase in oscillatory behaviour seen in the area under the secondary reflected shock. This test case demonstrates that higher polynomial degrees can be used for cases involving strong shocks.

\subsubsection{Investigation of Flux Reconstruction Schemes}
\begin{figure}[ht]
    \centering
    \includegraphics[width=0.85\textwidth,trim={0cm 1.5cm 0cm 6cm},clip]{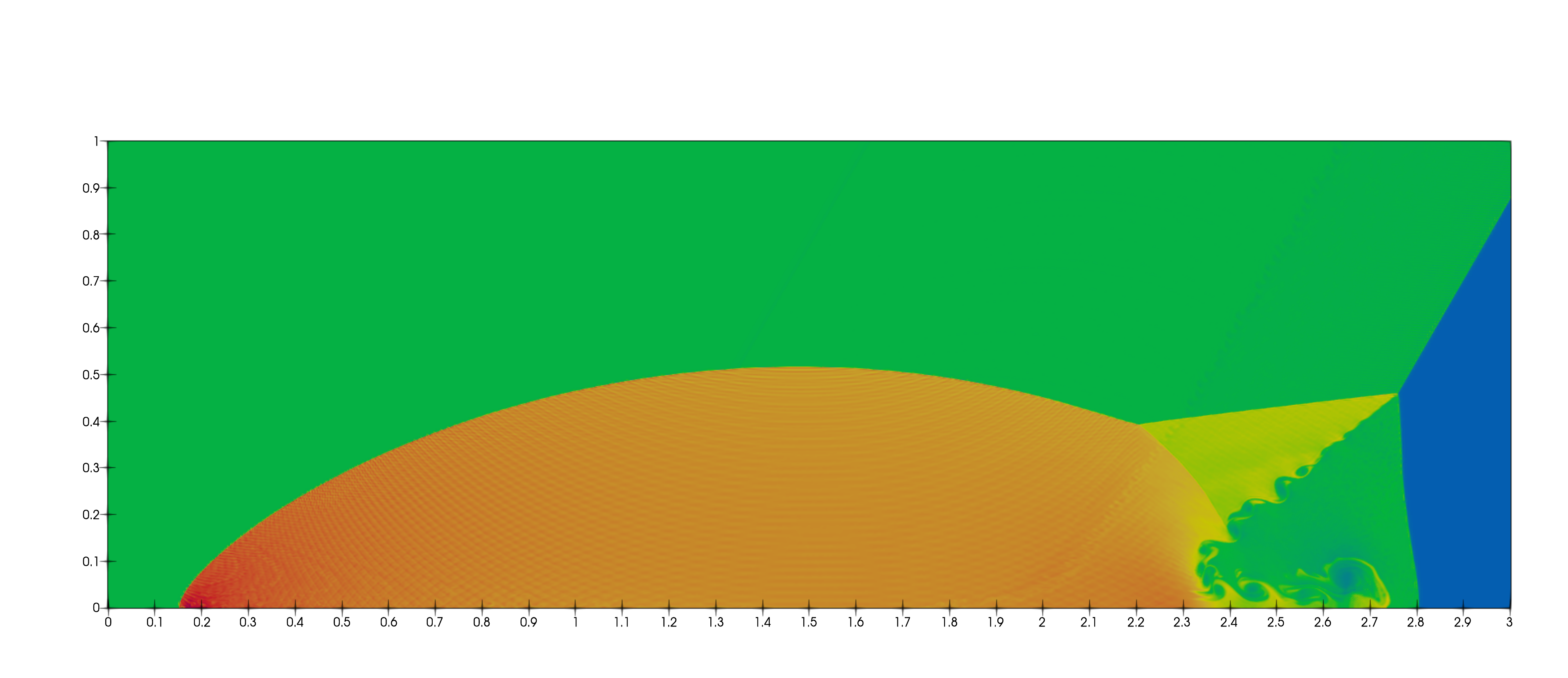}\\
    \includegraphics[width=0.85\textwidth,trim={0cm 1.5cm 0 6cm},clip]{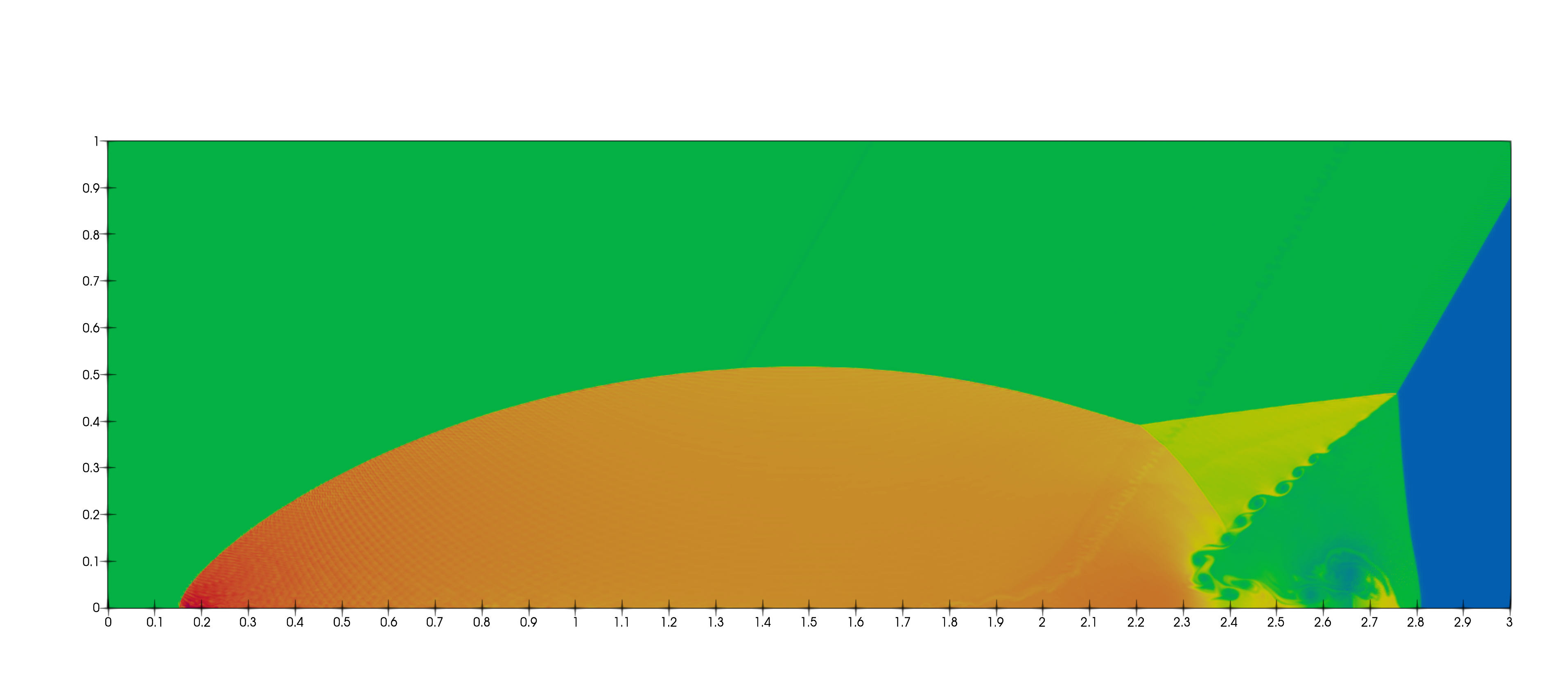}\\
    \caption{\textit{[DMR]} Density solution with NSFR, Roe dissipation, PPL at $t=0.2s$, using different flux reconstruction parameters - $c_{DG}$ (top), $c_{+}$ (bottom)}
    \label{fig:DMRFRParam}
\end{figure}

In addition to the study of polynomial degrees, the DMR case is also run with different flux reconstruction schemes to compare the effects of the $c$ value on the complex features. These tests are run with a finer grid of $720\times180$ for the test space ($720\times540$ for the full computational domain). The $CFL$ for $c_{DG}$ is 0.5 and for $c_+$ it is 0.58, which is the maximum $CFL$ for which the case can be run for the respective schemes. The results are shown in Fig.~\ref{fig:DMRFRParam}.

In both sets of results, there are two triple points and two slip lines present. The instabilities along the primary slip line, as well as the jet near the wall, are well defined. When comparing the primary slip line for the $c_{DG}$ and $c_{+}$ solutions, it is evident that the $c_{DG}$ slip line is much more oscillatory due to the noise in the solution. As a result, the KHI for the $c_{DG}$ solution is not as distinct as the KHI in the $c_{+}$ solution, which is unaffected by noise and has a distinct slip line without significant oscillations. The difference in noise between the two sets of results leads to another notable distinction between the two solutions. In the $c_{DG}$ solution, the secondary slip line is difficult to discern due to the noise in its vicinity, whereas the slip line is distinct and clearly visible in the $c_{+}$ solution due to the noise being filtered out. While both schemes are able to capture the main features of the solution accurately, the $c_{+}$ solution presents a clear advantage due to the distinct features and dampened noise.

\break
\subsection{High Mach Astrophysical Jet}
The High Mach Astrophysical Jet test case is an extreme benchmark to demonstrate the robustness of entropy-stable high-order schemes. This test is used to verify the robustness of the NSFR scheme with the PPL in two dimensions and demonstrates the advantages of the flux reconstruction scheme. This case was originally proposed by \cite{ha2005numerical} and it simulates an astrophysical jet with $Ma \approx 80$. Another version of this case is considered in \cite{ZHANG20108918} with an astrophysical jet of $Ma \approx 2000$. The two versions are outlined in this section, along with the results obtained using the NSFR scheme. Due to the high Mach numbers, this test is often run with either a TVB limiter as seen in \cite{ZHANG20108918} or with the imposition of TVD-like bounds as seen in \cite{hennemann2021provably}. Without the use of these strategies, running either version of this case using the SSPRK3 Strong DG method or the $c_{DG}$ scheme would require an extremely small $CFL$ number. A $CFL$ number as low as $0.0005$ was tested for both, but still failed due to nonphysical values. This case demonstrates the added robustness of the flux reconstruction, $c_+$ scheme, and the advantages of higher flux reconstruction correction parameters in extreme cases such as this.

\subsubsection{Mach 80 Astrophysical Jet}
The computational domain for the Mach 80 astrophysical jet is $[0,2]\times[-0.5,0.5]$ with the initial condition
\begin{equation}
    (\rho, u, v, p) = (0.5,0,0,0.4127).
\end{equation}
The top, right and bottom boundaries are outflow boundaries. The left boundary is set to $(\rho, u, v, p) = (0.5,30,0,0.4127)$ for $y\in[-0.05,0.05]$ and $(\rho, u, v, p) = (0.5,0,0,0.4127)$ otherwise. The final time is $t = 0.07$, and $N_x \times N_y = 320 \times 160$ with $\gamma = \frac{5}{3}$.

\begin{figure}[ht]
    \centering
    \includegraphics[width=\textwidth,trim={0cm 2cm 0cm 0cm},clip]{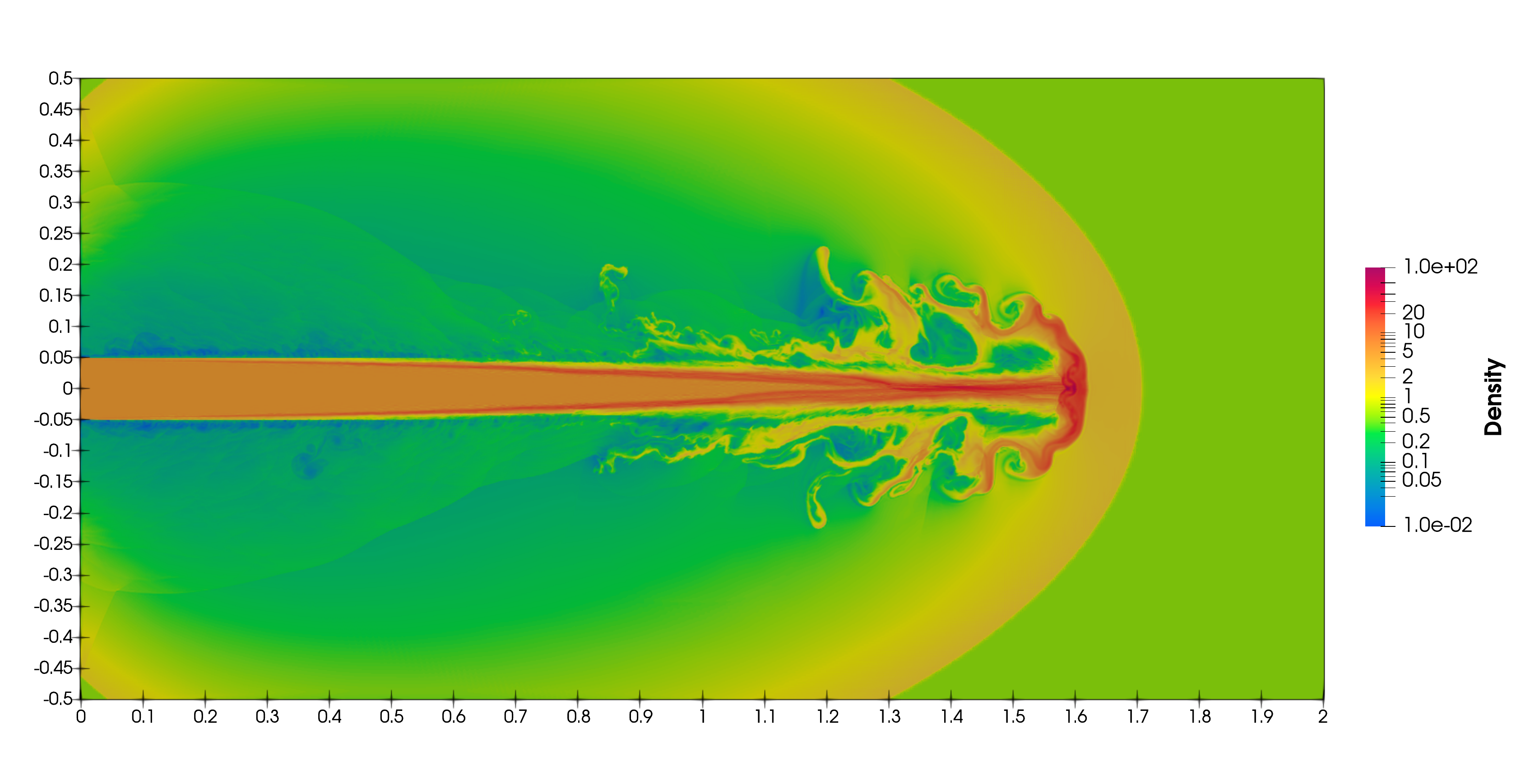}
    \caption{\textit{[Astrophysical Jet Mach 80]} Density solution with NSFR-$CH_{RA}$, Roe dissipation, $c_{+}$, PPL at $t=0.07s$.}
    \label{fig:AstroMach80}
\end{figure}

 Using the added robustness of the $c_+$ scheme, this case was run for a $CFL$ number of $0.01$. The results are shown in Fig.~\ref{fig:AstroMach80}. Even without the use of TVD limiting, the positivity of density and pressure is preserved. The results contain the expected solution features, including a bow shock, terminal Mach disk in the tip of the jet, and KHI.

\break
\subsubsection{Mach 2000 Astrophysical Jet}
The computational domain for the Mach 2000 astrophysical jet is $[-0.5,0.5]^2$ with the initial condition
\begin{equation}
    (\rho, u, v, p) = (0.5,0,0,0.4127)
\end{equation}
The top, right, and bottom boundaries are outflow boundaries. The left boundary is set to $(\rho, u, v, p) = (0.5,800,0,0.4127)$ for $y\in[-0.05,0.05]$ and $(\rho, u, v, p) = (0.5,0,0,0.4127)$ otherwise. The hypersonic inflow condition on the left boundary corresponds to a Mach number of $Ma=2156.91$. The final time is $t = 0.001$, and $N_x \times N_y = 300 \times 300$ with $\gamma = \frac{5}{3}$.

\begin{figure}
    \centering
    \includegraphics[width=0.75\textwidth,trim={0cm 2cm 0cm 0cm},clip]{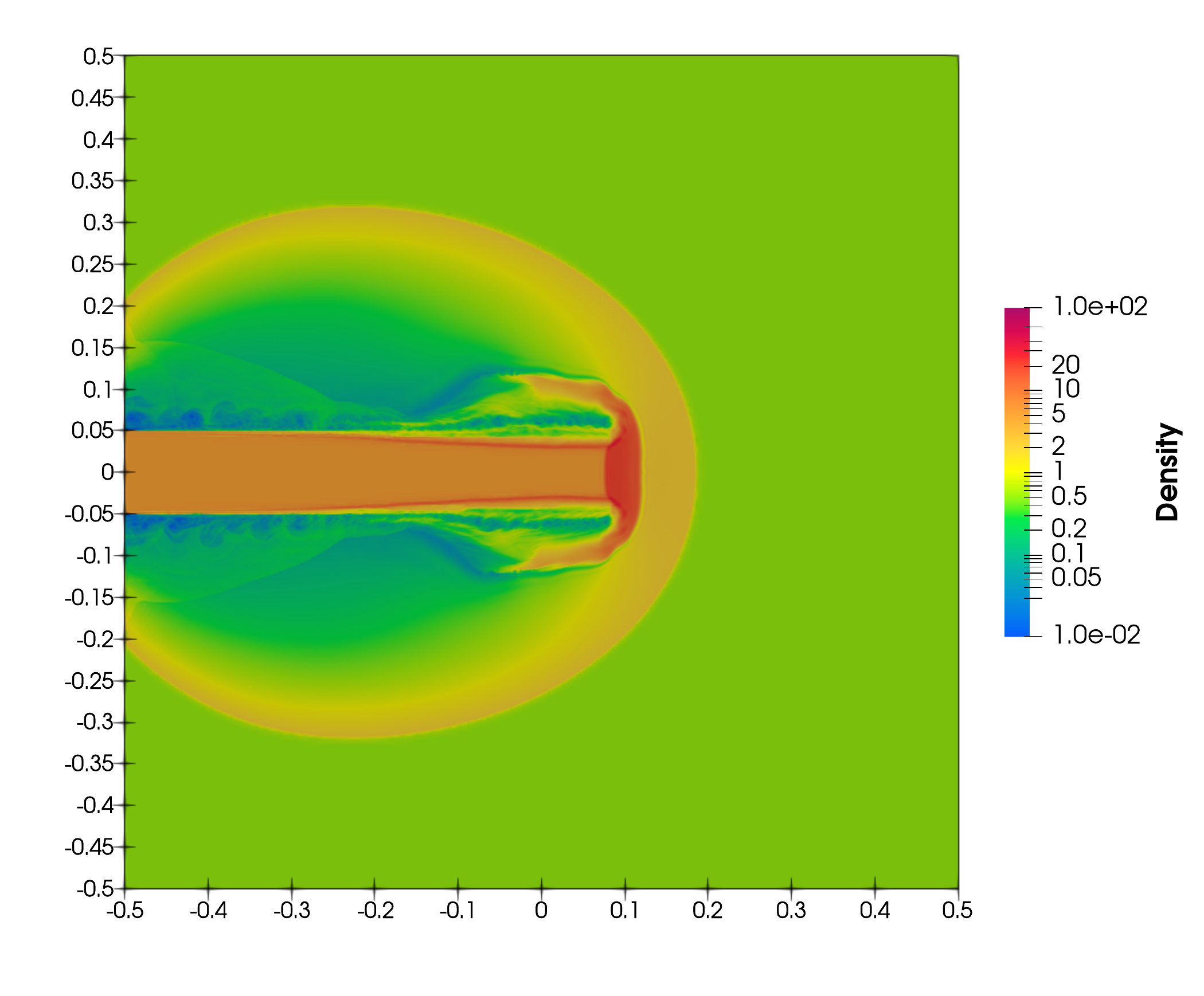}
    \caption{\textit{[Astrophysical Jet Mach 2000]} Density solution with NSFR-$CH_{RA}$, Roe dissipation, $c_{+}\times10$, PPL at $t=0.001s$.}
    \label{fig:AstroMach2000}
\end{figure}

Similar to the Mach 80 Astrophysical Jet case, due to the lack of TVD limiting, running this case with the SSPRK3 Strong DG method and $c_{DG}$ method is not feasible. Using the $c_+$ flux reconstruction parameter also requires an extremely small $CFL$ number to run successfully. To allow for a higher $CFL$ number, the flux reconstruction parameter is increased to ten times the value of $c_+$. With $c = 3.76E-2$ and a $CFL$ number of $0.0005$, the case is successfully run to completion while preserving the positivity of density and pressure. The results are shown in Fig.~\ref{fig:AstroMach2000}. The features of the solution are well captured, including the bow shock, terminal Mach disk and the unstable rollups. In previous tests, such as the DMR case, variants of flux reconstruction schemes were proven to have an advantage due to the dampening of noise, which produces a cleaner solution. In the Mach 80 and Mach 2000 cases, it is shown that the flux reconstruction scheme with higher values of the parameter, $c$, also provides the added advantage of robustness. 

\section{Conclusion}\label{sec: conclusion}
In this paper, the NSFR schemes \cite{CicchinoNonlinearlyStableFluxReconstruction2021, cicchino2022provably}, which have been successfully applied to unsteady compressible flows in arbitrary curvilinear coordinates, have been extended to shock-dominated problems using a modified implementation of the positivity-preserving limiter originally developed by \citet{ZHANG20108918}. Building on the improvements by \citet{WANG2012653}, the limiter is further modified to detect minimum values of density and pressure at the solution nodes, addressing the limitation of the original approach, which ensured positivity only at the quadrature nodes but not at the solution nodes. The enhancements to the positivity-preserving limiter were required to simulate all two-dimensional strong shock-dominated cases. The newly modified limiter guarantees positivity for the entire duration of the test at the same grid size and time step. With the robustness enhanced, this limiter is used to investigate the shock-capturing capabilities of the NSFR scheme. In the 1D Sod Shock Tube and 2D Strong Vortex Shock Wave (SVSW) Interaction cases, the advantages of the NSFR scheme over standard DG methods is established as the scheme can be used to run these cases without the limiter while the DG method requires the use of the PPL and a TVD limiter. In the 1D Leblanc Shock Tube case, the advantages of the NSFR method are clearly demonstrated, specifically the results for the $c_+$ FR scheme, which presents an essentially oscillation-free solution despite the extreme nature of this benchmark problem. The 2D Low Density case highlights the high-order accuracy of the scheme and verifies that the modified limiter preserves high-order accuracy. The shock-dominated cases also present a thorough analysis of the impact of two-point fluxes, quadrature nodes and variants of flux reconstruction schemes through the choice of the $c$ parameter. For the two-point flux investigation, the SVSW showcases the insufficient dissipation of the KG flux, which leads to increasing entropy. Additionally, the 1D Sod Shock Tube, 1D Shu Osher and 2D SVSW case are also used to numerically verify a CFL condition to preserve positivity for two-point fluxes. In the case of quadrature nodes, the GL quadrature is preferred as it has a stronger quadrature strength; however, due to the sensitivities of the entropy projections at the face, the GLL quadrature presents a significant timestep advantage over GL nodes as showcased in the 1D Shu Osher and 2D SVSW cases. For the FR parameter, all the test cases demonstrate that the timestep efficiency is maximized while suppressing oscillations as the value of the FR parameter is increased. The advantage of flux reconstruction is best highlighted in the astrophysical Mach jet case. With an FR scheme with $c$ being ten times greater than that of  $c_+$, the flow can be resolved without the need for TVD limiting. 

It has been demonstrated that FR schemes maximize the time step and increase robustness at the cost of reduced accuracy in smooth regions of the flow. This is due to the filtering of the highest mode as FR schemes with larger values of $c$ are considered.  The DG method is certainly more accurate, and we wish to preserve this accuracy in smooth regions while also mitigating oscillations in the vicinity of shocks, maximizing the time step, and ensuring robustness. Thus, the ability to scale the FR parameter based on local discontinuities represents a key advantage of the NSFR scheme, making an adaptive implementation of this parameter a promising direction for future research. In addition, an investigation of the NSFR scheme's performance for viscous shock problems and shock-turbulence interaction cases is part of ongoing research to further develop the approach to be suitable for a wide range of flows.

\section*{Acknowledgments}
\indent Sai Shruthi Srinivasan thanks the Vadasz Family Foundation and McGill Engineering Doctoral Award (MEDA). The authors would also like to thank Alexander Cicchino for helpful discussions throughout the course of this work. This research was enabled in part by support provided by Calcul Quebec and the Digital Research Alliance of Canada.

\newpage
\bibliographystyle{model1-num-names}
\bibliography{refs}

\begin{thebibliography}{58}
\expandafter\ifx\csname natexlab\endcsname\relax\def\natexlab#1{#1}\fi
\providecommand{\url}[1]{\texttt{#1}}
\providecommand{\href}[2]{#2}
\providecommand{\path}[1]{#1}
\providecommand{\DOIprefix}{doi:}
\providecommand{\ArXivprefix}{arXiv:}
\providecommand{\URLprefix}{URL: }
\providecommand{\Pubmedprefix}{pmid:}
\providecommand{\doi}[1]{\href{http://dx.doi.org/#1}{\path{#1}}}
\providecommand{\Pubmed}[1]{\href{pmid:#1}{\path{#1}}}
\providecommand{\bibinfo}[2]{#2}
\ifx\xfnm\relax \def\xfnm[#1]{\unskip,\space#1}\fi
\bibitem[{Zhang and Shu(2010)}]{ZHANG20103091}
\bibinfo{author}{X.~Zhang}, \bibinfo{author}{C.-W. Shu},
\newblock \bibinfo{title}{On maximum-principle-satisfying high order schemes
  for scalar conservation laws},
\newblock \bibinfo{journal}{Journal of Computational Physics}
  \bibinfo{volume}{229} (\bibinfo{year}{2010}) \bibinfo{pages}{3091--3120}.
\bibitem[{Reed and Hill(1973)}]{reed1973triangular}
\bibinfo{author}{W.~H. Reed}, \bibinfo{author}{T.~R. Hill},
  \bibinfo{title}{Triangular mesh methods for the neutron transport equation},
  \bibinfo{type}{Technical Report}, Los Alamos Scientific Lab., N. Mex.(USA),
  \bibinfo{year}{1973}.
\bibitem[{Hesthaven and Warburton(2007)}]{NDG}
\bibinfo{author}{J.~S. Hesthaven}, \bibinfo{author}{T.~Warburton},
  \bibinfo{title}{Nodal discontinuous Galerkin methods: algorithms, analysis,
  and applications}, \bibinfo{publisher}{Springer Science \& Business Media},
  \bibinfo{year}{2007}.
\bibitem[{Gassner(2013)}]{gassner2013skew}
\bibinfo{author}{G.~J. Gassner},
\newblock \bibinfo{title}{A skew-symmetric discontinuous {Galerkin} spectral
  element discretization and its relation to {SBP-SAT} finite difference
  methods},
\newblock \bibinfo{journal}{SIAM Journal on Scientific Computing}
  \bibinfo{volume}{35} (\bibinfo{year}{2013}) \bibinfo{pages}{A1233--A1253}.
\bibitem[{{Del Rey Fern\'andez} et~al.(2014){Del Rey Fern\'andez}, Boom, and
  Zingg}]{fernandez2014generalized}
\bibinfo{author}{D.~C. {Del Rey Fern\'andez}}, \bibinfo{author}{P.~D. Boom},
  \bibinfo{author}{D.~W. Zingg},
\newblock \bibinfo{title}{A generalized framework for nodal first derivative
  summation-by-parts operators},
\newblock \bibinfo{journal}{Journal of Computational Physics}
  \bibinfo{volume}{266} (\bibinfo{year}{2014}) \bibinfo{pages}{214--239}.
\bibitem[{{Del Rey Fern\'andez} et~al.(2020){Del Rey Fern\'andez}, Carpenter,
  Dalcin, Zampini, and Parsani}]{fernandez2019entropy}
\bibinfo{author}{D.~C. {Del Rey Fern\'andez}}, \bibinfo{author}{M.~H.
  Carpenter}, \bibinfo{author}{L.~Dalcin}, \bibinfo{author}{S.~Zampini},
  \bibinfo{author}{M.~Parsani},
\newblock \bibinfo{title}{Entropy stable h/p-nonconforming discretization with
  the summation-by-parts property for the compressible {Euler} and
  {Navier}--{Stokes} equations},
\newblock \bibinfo{journal}{SN Partial Differential Equations and Applications}
  \bibinfo{volume}{1} (\bibinfo{year}{2020}) \bibinfo{pages}{1--54}.
\bibitem[{Crean et~al.(2018)Crean, Hicken, {Del Rey Fern\'andez}, Zingg, and
  Carpenter}]{crean2018entropy}
\bibinfo{author}{J.~Crean}, \bibinfo{author}{J.~E. Hicken},
  \bibinfo{author}{D.~C. {Del Rey Fern\'andez}}, \bibinfo{author}{D.~W. Zingg},
  \bibinfo{author}{M.~H. Carpenter},
\newblock \bibinfo{title}{Entropy-stable summation-by-parts discretization of
  the {Euler} equations on general curved elements},
\newblock \bibinfo{journal}{Journal of Computational Physics}
  \bibinfo{volume}{356} (\bibinfo{year}{2018}) \bibinfo{pages}{410--438}.
\bibitem[{Chan(2018)}]{chan2018discretely}
\bibinfo{author}{J.~Chan},
\newblock \bibinfo{title}{On discretely entropy conservative and entropy stable
  discontinuous {Galerkin} methods},
\newblock \bibinfo{journal}{Journal of Computational Physics}
  \bibinfo{volume}{362} (\bibinfo{year}{2018}) \bibinfo{pages}{346--374}.
\bibitem[{Ranocha et~al.(2016)Ranocha, {\"O}ffner, and
  Sonar}]{ranocha2016summation}
\bibinfo{author}{H.~Ranocha}, \bibinfo{author}{P.~{\"O}ffner},
  \bibinfo{author}{T.~Sonar},
\newblock \bibinfo{title}{Summation-by-parts operators for correction procedure
  via reconstruction},
\newblock \bibinfo{journal}{Journal of Computational Physics}
  \bibinfo{volume}{311} (\bibinfo{year}{2016}) \bibinfo{pages}{299--328}.
\bibitem[{Fisher et~al.(2013)Fisher, Carpenter, Nordstr{\"o}m, Yamaleev, and
  Swanson}]{fisher2013discretely}
\bibinfo{author}{T.~C. Fisher}, \bibinfo{author}{M.~H. Carpenter},
  \bibinfo{author}{J.~Nordstr{\"o}m}, \bibinfo{author}{N.~K. Yamaleev},
  \bibinfo{author}{C.~Swanson},
\newblock \bibinfo{title}{Discretely conservative finite-difference
  formulations for nonlinear conservation laws in split form: {Theory} and
  boundary conditions},
\newblock \bibinfo{journal}{Journal of Computational Physics}
  \bibinfo{volume}{234} (\bibinfo{year}{2013}) \bibinfo{pages}{353--375}.
\bibitem[{Montoya and Zingg(2021)}]{montoya2021unifying}
\bibinfo{author}{T.~Montoya}, \bibinfo{author}{D.~W. Zingg},
\newblock \bibinfo{title}{A unifying algebraic framework for discontinuous
  {Galerkin} and flux reconstruction methods based on the summation-by-parts
  property},
\newblock \bibinfo{journal}{arXiv preprint arXiv:2101.10478v1}
  (\bibinfo{year}{2021}).
\bibitem[{Tadmor(1984)}]{tadmor1984skew}
\bibinfo{author}{E.~Tadmor},
\newblock \bibinfo{title}{Skew-self adjoint form for systems of conservation
  laws},
\newblock \bibinfo{journal}{Journal of Mathematical Analysis and Applications}
  \bibinfo{volume}{103} (\bibinfo{year}{1984}) \bibinfo{pages}{428--442}.
\bibitem[{Huynh(2007)}]{huynh_flux_2007}
\bibinfo{author}{H.~T. Huynh},
\newblock \bibinfo{title}{A {flux} {reconstruction} {approach} to
  {high}-{order} {schemes} {including} {discontinuous} {Galerkin} {methods}},
\newblock \bibinfo{publisher}{American Institute of Aeronautics and
  Astronautics}, \bibinfo{year}{2007}. \DOIprefix\doi{10.2514/6.2007-4079}.
\bibitem[{Vincent et~al.(2011)Vincent, Castonguay, and
  Jameson}]{vincent_new_2011}
\bibinfo{author}{P.~E. Vincent}, \bibinfo{author}{P.~Castonguay},
  \bibinfo{author}{A.~Jameson},
\newblock \bibinfo{title}{A {new} {class} of {high}-{order} {energy} {stable}
  {flux} {reconstruction} {schemes}},
\newblock \bibinfo{journal}{Journal of Scientific Computing}
  \bibinfo{volume}{47} (\bibinfo{year}{2011}) \bibinfo{pages}{50--72}.
\bibitem[{Jameson(2010)}]{jameson_proof_2010}
\bibinfo{author}{A.~Jameson},
\newblock \bibinfo{title}{A {proof} of the {stability} of the {spectral}
  {difference} {method} for {all} {orders} of {accuracy}},
\newblock \bibinfo{journal}{Journal of Scientific Computing}
  \bibinfo{volume}{45} (\bibinfo{year}{2010}) \bibinfo{pages}{348--358}.
\bibitem[{Wang and Gao(2009)}]{wang2009unifying}
\bibinfo{author}{Z.~J. Wang}, \bibinfo{author}{H.~Gao},
\newblock \bibinfo{title}{A unifying lifting collocation penalty formulation
  including the discontinuous {Galerkin}, spectral volume/difference methods
  for conservation laws on mixed grids},
\newblock \bibinfo{journal}{Journal of Computational Physics}
  \bibinfo{volume}{228} (\bibinfo{year}{2009}) \bibinfo{pages}{8161--8186}.
\bibitem[{Zwanenburg and Nadarajah(2016)}]{zwanenburg_equivalence_2016}
\bibinfo{author}{P.~Zwanenburg}, \bibinfo{author}{S.~Nadarajah},
\newblock \bibinfo{title}{Equivalence between the {energy} {stable} {flux}
  {reconstruction} and {filtered} {discontinuous} {Galerkin} {schemes}},
\newblock \bibinfo{journal}{Journal of Computational Physics}
  \bibinfo{volume}{306} (\bibinfo{year}{2016}) \bibinfo{pages}{343--369}.
\bibitem[{Allaneau and Jameson(2011)}]{allaneau_connections_2011}
\bibinfo{author}{Y.~Allaneau}, \bibinfo{author}{A.~Jameson},
\newblock \bibinfo{title}{Connections between the filtered discontinuous
  {Galerkin} method and the flux reconstruction approach to high order
  discretizations},
\newblock \bibinfo{journal}{Computer Methods in Applied Mechanics and
  Engineering} \bibinfo{volume}{200} (\bibinfo{year}{2011})
  \bibinfo{pages}{3628--3636}.
\bibitem[{{\AA}lund and Nordstr{\"o}m(2019)}]{aalund2019encapsulated}
\bibinfo{author}{O.~{\AA}lund}, \bibinfo{author}{J.~Nordstr{\"o}m},
\newblock \bibinfo{title}{Encapsulated high order difference operators on
  curvilinear non-conforming grids},
\newblock \bibinfo{journal}{Journal of Computational Physics}
  \bibinfo{volume}{385} (\bibinfo{year}{2019}) \bibinfo{pages}{209--224}.
\bibitem[{Liu et~al.(2006)Liu, Vinokur, and Wang}]{liu_spectral_2006}
\bibinfo{author}{Y.~Liu}, \bibinfo{author}{M.~Vinokur}, \bibinfo{author}{Z.~J.
  Wang},
\newblock \bibinfo{title}{Spectral difference method for unstructured grids
  {I}: {Basic} formulation},
\newblock \bibinfo{journal}{Journal of Computational Physics}
  \bibinfo{volume}{216} (\bibinfo{year}{2006}) \bibinfo{pages}{780--801}.
\bibitem[{Abe et~al.(2018)Abe, Morinaka, Haga, Nonomura, Shibata, and
  Miyaji}]{abe2018stable}
\bibinfo{author}{Y.~Abe}, \bibinfo{author}{I.~Morinaka},
  \bibinfo{author}{T.~Haga}, \bibinfo{author}{T.~Nonomura},
  \bibinfo{author}{H.~Shibata}, \bibinfo{author}{K.~Miyaji},
\newblock \bibinfo{title}{Stable, non-dissipative, and conservative
  flux-reconstruction schemes in split forms},
\newblock \bibinfo{journal}{Journal of Computational Physics}
  \bibinfo{volume}{353} (\bibinfo{year}{2018}) \bibinfo{pages}{193--227}.
\bibitem[{Cicchino et~al.(2022{\natexlab{a}})Cicchino, Nadarajah, and {Del Rey
  Fernández}}]{CicchinoNonlinearlyStableFluxReconstruction2021}
\bibinfo{author}{A.~Cicchino}, \bibinfo{author}{S.~Nadarajah},
  \bibinfo{author}{D.~C. {Del Rey Fernández}},
\newblock \bibinfo{title}{Nonlinearly stable flux reconstruction high-order
  methods in split form},
\newblock \bibinfo{journal}{Journal of Computational Physics}
  (\bibinfo{year}{2022}{\natexlab{a}}) \bibinfo{pages}{111094}.
\bibitem[{Cicchino et~al.(2022{\natexlab{b}})Cicchino, Del Rey~Fern{\'a}ndez,
  Nadarajah, Chan, and Carpenter}]{cicchino2022provably}
\bibinfo{author}{A.~Cicchino}, \bibinfo{author}{D.~C. Del Rey~Fern{\'a}ndez},
  \bibinfo{author}{S.~Nadarajah}, \bibinfo{author}{J.~Chan},
  \bibinfo{author}{M.~H. Carpenter},
\newblock \bibinfo{title}{Provably stable flux reconstruction high-order
  methods on curvilinear elements},
\newblock \bibinfo{journal}{Journal of Computational Physics}
  \bibinfo{volume}{463} (\bibinfo{year}{2022}{\natexlab{b}})
  \bibinfo{pages}{111259}.
\bibitem[{Cicchino and Nadarajah(2025)}]{cicchino2025discretely}
\bibinfo{author}{A.~Cicchino}, \bibinfo{author}{S.~Nadarajah},
\newblock \bibinfo{title}{Discretely nonlinearly stable weight-adjusted flux
  reconstruction high-order method for compressible flows on curvilinear
  grids},
\newblock \bibinfo{journal}{Journal of Computational Physics}
  \bibinfo{volume}{521} (\bibinfo{year}{2025}) \bibinfo{pages}{113532}.
\bibitem[{Dafermos(2016)}]{CMDAFERMOS}
\bibinfo{author}{C.~Dafermos}, \bibinfo{title}{Hyperbolic Conservation Laws in
  Continuum Physics}, \bibinfo{year}{2016}.
  \DOIprefix\doi{10.1007/978-3-662-49451-6}.
\bibitem[{Harten(1997)}]{harten1997highresolution}
\bibinfo{author}{A.~Harten},
\newblock \bibinfo{title}{High resolution schemes for hyperbolic conservation
  laws},
\newblock \bibinfo{journal}{Journal of computational physics}
  \bibinfo{volume}{135} (\bibinfo{year}{1997}) \bibinfo{pages}{260--278}.
\bibitem[{Richtmyer et~al.(1950)}]{richtmyer1950method}
\bibinfo{author}{R.~Richtmyer}, et~al.,
\newblock \bibinfo{title}{{A method for the numerical calculation of
  hydrodynamic shocks}},
\newblock \bibinfo{journal}{Journal of Applied Physics} \bibinfo{volume}{21}
  (\bibinfo{year}{1950}) \bibinfo{pages}{232--237}.
\bibitem[{Wintermeyer et~al.(2017)Wintermeyer, Winters, Gassner, and
  Kopriva}]{wintermeyer2017entropy}
\bibinfo{author}{N.~Wintermeyer}, \bibinfo{author}{A.~R. Winters},
  \bibinfo{author}{G.~J. Gassner}, \bibinfo{author}{D.~A. Kopriva},
\newblock \bibinfo{title}{An entropy stable nodal discontinuous {Galerkin}
  method for the two dimensional shallow water equations on unstructured
  curvilinear meshes with discontinuous bathymetry},
\newblock \bibinfo{journal}{Journal of Computational Physics}
  \bibinfo{volume}{340} (\bibinfo{year}{2017}) \bibinfo{pages}{200--242}.
\bibitem[{Wintermeyer et~al.(2018)Wintermeyer, Winters, Gassner, and
  Warburton}]{wintermeyer2018entropy}
\bibinfo{author}{N.~Wintermeyer}, \bibinfo{author}{A.~R. Winters},
  \bibinfo{author}{G.~J. Gassner}, \bibinfo{author}{T.~Warburton},
\newblock \bibinfo{title}{{An entropy stable discontinuous Galerkin method for
  the shallow water equations on curvilinear meshes with wet/dry fronts
  accelerated by GPUs}},
\newblock \bibinfo{journal}{Journal of Computational Physics}
  \bibinfo{volume}{375} (\bibinfo{year}{2018}) \bibinfo{pages}{447--480}.
\bibitem[{Dumbser et~al.(2014)Dumbser, Zanotti, Loub{\`e}re, and
  Diot}]{dumbser2014posteriori}
\bibinfo{author}{M.~Dumbser}, \bibinfo{author}{O.~Zanotti},
  \bibinfo{author}{R.~Loub{\`e}re}, \bibinfo{author}{S.~Diot},
\newblock \bibinfo{title}{{A posteriori subcell limiting of the discontinuous
  Galerkin finite element method for hyperbolic conservation laws}},
\newblock \bibinfo{journal}{Journal of Computational Physics}
  \bibinfo{volume}{278} (\bibinfo{year}{2014}) \bibinfo{pages}{47--75}.
\bibitem[{Hennemann et~al.(2021)Hennemann, Rueda-Ram{\'\i}rez, Hindenlang, and
  Gassner}]{hennemann2021provably}
\bibinfo{author}{S.~Hennemann}, \bibinfo{author}{A.~M. Rueda-Ram{\'\i}rez},
  \bibinfo{author}{F.~J. Hindenlang}, \bibinfo{author}{G.~J. Gassner},
\newblock \bibinfo{title}{A provably entropy stable subcell shock capturing
  approach for high order split form {DG} for the compressible {Euler}
  equations},
\newblock \bibinfo{journal}{Journal of Computational Physics}
  \bibinfo{volume}{426} (\bibinfo{year}{2021}) \bibinfo{pages}{109935}.
\bibitem[{Rueda-Ram{\'\i}rez and Gassner(2021)}]{rueda2021subcell}
\bibinfo{author}{A.~M. Rueda-Ram{\'\i}rez}, \bibinfo{author}{G.~J. Gassner},
\newblock \bibinfo{title}{{A subcell finite volume positivity-preserving
  limiter for DGSEM discretizations of the Euler equations}},
\newblock \bibinfo{journal}{arXiv preprint arXiv:2102.06017}
  (\bibinfo{year}{2021}).
\bibitem[{Lin and Chan(2024)}]{lin2024high}
\bibinfo{author}{Y.~Lin}, \bibinfo{author}{J.~Chan},
\newblock \bibinfo{title}{High order entropy stable discontinuous galerkin
  spectral element methods through subcell limiting},
\newblock \bibinfo{journal}{J. Comput. Phys.} \bibinfo{volume}{498}
  (\bibinfo{year}{2024}) \bibinfo{pages}{112677}.
\bibitem[{Sonntag and Munz(2014)}]{sonntag2014shock}
\bibinfo{author}{M.~Sonntag}, \bibinfo{author}{C.-D. Munz},
\newblock \bibinfo{title}{{Shock capturing for discontinuous Galerkin methods
  using finite volume subcells}},
\newblock in: \bibinfo{booktitle}{Finite Volumes for Complex Applications
  VII-Elliptic, Parabolic and Hyperbolic Problems: FVCA 7, Berlin, June 2014},
  \bibinfo{organization}{Springer}, \bibinfo{year}{2014}, pp.
  \bibinfo{pages}{945--953}.
\bibitem[{Vilar(2019)}]{vilar2019posteriori}
\bibinfo{author}{F.~Vilar},
\newblock \bibinfo{title}{{A posteriori correction of high-order discontinuous
  Galerkin scheme through subcell finite volume formulation and flux
  reconstruction}},
\newblock \bibinfo{journal}{Journal of Computational Physics}
  \bibinfo{volume}{387} (\bibinfo{year}{2019}) \bibinfo{pages}{245--279}.
\bibitem[{Zhang and Shu(2010)}]{zhang2010maximumpreservinglimiter}
\bibinfo{author}{X.~Zhang}, \bibinfo{author}{C.-W. Shu},
\newblock \bibinfo{title}{On maximum-principle-satisfying high order schemes
  for scalar conservation laws},
\newblock \bibinfo{journal}{Journal of Computational Physics}
  \bibinfo{volume}{229} (\bibinfo{year}{2010}) \bibinfo{pages}{3091--3120}.
\bibitem[{Chen and Shu(2017)}]{chen2017entropy}
\bibinfo{author}{T.~Chen}, \bibinfo{author}{C.-W. Shu},
\newblock \bibinfo{title}{Entropy stable high order discontinuous {Galerkin}
  methods with suitable quadrature rules for hyperbolic conservation laws},
\newblock \bibinfo{journal}{Journal of Computational Physics}
  \bibinfo{volume}{345} (\bibinfo{year}{2017}) \bibinfo{pages}{427--461}.
\bibitem[{Cicchino and Nadarajah(2020)}]{Cicchino2020NewNorm}
\bibinfo{author}{A.~Cicchino}, \bibinfo{author}{S.~Nadarajah},
\newblock \bibinfo{title}{A new norm and stability condition for tensor product
  flux reconstruction schemes},
\newblock \bibinfo{journal}{Journal of Computational Physics}
  (\bibinfo{year}{2020}) \bibinfo{pages}{110025}.
\bibitem[{Sheshadri and Jameson(2016)}]{sheshadri2016stability}
\bibinfo{author}{A.~Sheshadri}, \bibinfo{author}{A.~Jameson},
\newblock \bibinfo{title}{On the stability of the flux reconstruction schemes
  on quadrilateral elements for the linear advection equation},
\newblock \bibinfo{journal}{Journal of Scientific Computing}
  \bibinfo{volume}{67} (\bibinfo{year}{2016}) \bibinfo{pages}{769--790}.
\bibitem[{Chan and Wilcox(2019)}]{chan2019discretely}
\bibinfo{author}{J.~Chan}, \bibinfo{author}{L.~C. Wilcox},
\newblock \bibinfo{title}{On discretely entropy stable weight-adjusted
  discontinuous {Galerkin} methods: {Curvilinear} meshes},
\newblock \bibinfo{journal}{Journal of Computational Physics}
  \bibinfo{volume}{378} (\bibinfo{year}{2019}) \bibinfo{pages}{366--393}.
\bibitem[{Chan(2019)}]{chan2019skew}
\bibinfo{author}{J.~Chan},
\newblock \bibinfo{title}{Skew-symmetric entropy stable modal discontinuous
  {Galerkin} formulations},
\newblock \bibinfo{journal}{Journal of Scientific Computing}
  \bibinfo{volume}{81} (\bibinfo{year}{2019}) \bibinfo{pages}{459--485}.
\bibitem[{Cicchino(2024)}]{cicchino2024weight}
\bibinfo{author}{A.~Cicchino},
\newblock \bibinfo{title}{Weight-adjusted nonlinearly stable flux
  reconstruction high-order methods for compressible flows in curvilinear
  coordinates}  (\bibinfo{year}{2024}).
\bibitem[{Zhang and Shu(2010)}]{ZHANG20108918}
\bibinfo{author}{X.~Zhang}, \bibinfo{author}{C.-W. Shu},
\newblock \bibinfo{title}{On positivity-preserving high order discontinuous
  galerkin schemes for compressible euler equations on rectangular meshes},
\newblock \bibinfo{journal}{Journal of Computational Physics}
  \bibinfo{volume}{229} (\bibinfo{year}{2010}) \bibinfo{pages}{8918--8934}.
\bibitem[{Wang et~al.(2012)Wang, Zhang, Shu, and Ning}]{WANG2012653}
\bibinfo{author}{C.~Wang}, \bibinfo{author}{X.~Zhang}, \bibinfo{author}{C.-W.
  Shu}, \bibinfo{author}{J.~Ning},
\newblock \bibinfo{title}{Robust high order discontinuous galerkin schemes for
  two-dimensional gaseous detonations},
\newblock \bibinfo{journal}{Journal of Computational Physics}
  \bibinfo{volume}{231} (\bibinfo{year}{2012}) \bibinfo{pages}{653--665}.
\bibitem[{Ranocha and Gassner(2021)}]{ranocha2021preventing}
\bibinfo{author}{H.~Ranocha}, \bibinfo{author}{G.~J. Gassner},
\newblock \bibinfo{title}{Preventing pressure oscillations does not fix local
  linear stability issues of entropy-based split-form high-order schemes},
\newblock \bibinfo{journal}{Communications on Applied Mathematics and
  Computation}  (\bibinfo{year}{2021}) \bibinfo{pages}{1--24}.
\bibitem[{Ranocha(2018)}]{ranocha2018comparison}
\bibinfo{author}{H.~Ranocha},
\newblock \bibinfo{title}{Comparison of some entropy conservative numerical
  fluxes for the euler equations},
\newblock \bibinfo{journal}{J. Sci. Comput.} \bibinfo{volume}{76}
  (\bibinfo{year}{2018}) \bibinfo{pages}{216--242}.
\bibitem[{Castonguay(2012)}]{castonguay2012high}
\bibinfo{author}{P.~Castonguay}, \bibinfo{title}{High-order energy stable flux
  reconstruction schemes for fluid flow simulations on unstructured grids},
  Ph.D. thesis, Citeseer, \bibinfo{year}{2012}.
\bibitem[{Shu and Osher(1988)}]{shu1988efficient}
\bibinfo{author}{C.-W. Shu}, \bibinfo{author}{S.~Osher},
\newblock \bibinfo{title}{Efficient implementation of essentially
  non-oscillatory shock-capturing schemes},
\newblock \bibinfo{journal}{Journal of computational physics}
  \bibinfo{volume}{77} (\bibinfo{year}{1988}) \bibinfo{pages}{439--471}.
\bibitem[{Ismail and Roe(2009)}]{ismail2009affordable}
\bibinfo{author}{F.~Ismail}, \bibinfo{author}{P.~L. Roe},
\newblock \bibinfo{title}{Affordable, entropy-consistent euler flux functions
  ii: Entropy production at shocks},
\newblock \bibinfo{journal}{Journal of Computational Physics}
  \bibinfo{volume}{228} (\bibinfo{year}{2009}) \bibinfo{pages}{5410--5436}.
\bibitem[{Kennedy and Gruber(2008)}]{KENNEDY20081676}
\bibinfo{author}{C.~A. Kennedy}, \bibinfo{author}{A.~Gruber},
\newblock \bibinfo{title}{Reduced aliasing formulations of the convective terms
  within the navier–stokes equations for a compressible fluid},
\newblock \bibinfo{journal}{Journal of Computational Physics}
  \bibinfo{volume}{227} (\bibinfo{year}{2008}) \bibinfo{pages}{1676--1700}.
\bibitem[{Chandrashekar(2013)}]{chandrashekar2013kinetic}
\bibinfo{author}{P.~Chandrashekar},
\newblock \bibinfo{title}{Kinetic energy preserving and entropy stable finite
  volume schemes for compressible euler and navier-stokes equations},
\newblock \bibinfo{journal}{Communications in Computational Physics}
  \bibinfo{volume}{14} (\bibinfo{year}{2013}) \bibinfo{pages}{1252--1286}.
\bibitem[{Johnsen et~al.(2010)Johnsen, Larsson, Bhagatwala, Cabot, Moin, Olson,
  Rawat, Shankar, Sjögreen, Yee, Zhong, and Lele}]{JOHNSEN20101213}
\bibinfo{author}{E.~Johnsen}, \bibinfo{author}{J.~Larsson},
  \bibinfo{author}{A.~V. Bhagatwala}, \bibinfo{author}{W.~H. Cabot},
  \bibinfo{author}{P.~Moin}, \bibinfo{author}{B.~J. Olson},
  \bibinfo{author}{P.~S. Rawat}, \bibinfo{author}{S.~K. Shankar},
  \bibinfo{author}{B.~Sjögreen}, \bibinfo{author}{H.~Yee},
  \bibinfo{author}{X.~Zhong}, \bibinfo{author}{S.~K. Lele},
\newblock \bibinfo{title}{Assessment of high-resolution methods for numerical
  simulations of compressible turbulence with shock waves},
\newblock \bibinfo{journal}{Journal of Computational Physics}
  \bibinfo{volume}{229} (\bibinfo{year}{2010}) \bibinfo{pages}{1213--1237}.
\bibitem[{Xu and Shu(2022)}]{xu2022third}
\bibinfo{author}{Z.~Xu}, \bibinfo{author}{C.-W. Shu},
\newblock \bibinfo{title}{Third order maximum-principle-satisfying and
  positivity-preserving lax-wendroff discontinuous galerkin methods for
  hyperbolic conservation laws},
\newblock \bibinfo{journal}{Journal of Computational Physics}
  \bibinfo{volume}{470} (\bibinfo{year}{2022}) \bibinfo{pages}{111591}.
\bibitem[{Rault et~al.(2003)Rault, Chiavassa, and Donat}]{rault2003shock}
\bibinfo{author}{A.~Rault}, \bibinfo{author}{G.~Chiavassa},
  \bibinfo{author}{R.~Donat},
\newblock \bibinfo{title}{Shock-vortex interactions at high mach numbers},
\newblock \bibinfo{journal}{Journal of Scientific Computing}
  \bibinfo{volume}{19} (\bibinfo{year}{2003}) \bibinfo{pages}{347--371}.
\bibitem[{Galbraith et~al.(2017)Galbraith, Murman, Kim, Persson, Fidkowski,
  Glasby, Hillewaert, and Ahrabi}]{HOW5}
\bibinfo{author}{M.~Galbraith}, \bibinfo{author}{S.~Murman},
  \bibinfo{author}{C.~Kim}, \bibinfo{author}{P.~Persson},
  \bibinfo{author}{K.~Fidkowski}, \bibinfo{author}{R.~Glasby},
  \bibinfo{author}{K.~Hillewaert}, \bibinfo{author}{B.~Ahrabi},
\newblock \bibinfo{title}{{5th international workshop on high-order CFD
  methods}},
\newblock in: \bibinfo{booktitle}{AIAA science and technology forum and
  exposition}, \bibinfo{year}{2017}.
\bibitem[{Hillier(1991)}]{hillier1991computation}
\bibinfo{author}{R.~Hillier},
\newblock \bibinfo{title}{Computation of shock wave diffraction at a ninety
  degrees convex edge},
\newblock \bibinfo{journal}{Shock waves} \bibinfo{volume}{1}
  (\bibinfo{year}{1991}) \bibinfo{pages}{89--98}.
\bibitem[{Vevek et~al.(2019)Vevek, Zang, and New}]{vevek2019alternative}
\bibinfo{author}{U.~Vevek}, \bibinfo{author}{B.~Zang}, \bibinfo{author}{T.~H.
  New},
\newblock \bibinfo{title}{On alternative setups of the double mach reflection
  problem},
\newblock \bibinfo{journal}{Journal of Scientific Computing}
  \bibinfo{volume}{78} (\bibinfo{year}{2019}) \bibinfo{pages}{1291--1303}.
\bibitem[{Ha et~al.(2005)Ha, Gardner, Gelb, and Shu}]{ha2005numerical}
\bibinfo{author}{Y.~Ha}, \bibinfo{author}{C.~L. Gardner},
  \bibinfo{author}{A.~Gelb}, \bibinfo{author}{C.-W. Shu},
\newblock \bibinfo{title}{Numerical simulation of high mach number
  astrophysical jets with radiative cooling},
\newblock \bibinfo{journal}{Journal of Scientific Computing}
  \bibinfo{volume}{24} (\bibinfo{year}{2005}) \bibinfo{pages}{29--44}.

\end{thebibliography}

\end{document}